\documentclass[a4paper,10pt]{article}
\usepackage{hyperref}
\hypersetup{colorlinks, citecolor=blue, filecolor=black, linkcolor=blue, urlcolor=blue}
\usepackage[top=2.54cm, bottom=2.54cm, outer=2.75cm, inner=2.75cm, headsep=14pt]{geometry} 

\newcommand{\act}{\mathop\mathrm{\rm act}}
\newcommand{\supp}{\mbox{\rm supp}}
\newcommand{\cl}{\mbox{\rm cl\,}}
\usepackage{mathpazo}
\newcommand{\qrint}{\mbox{\rm qri\,}}
\usepackage{enumerate}
\usepackage[shortlabels]{enumitem}
\usepackage[latin1]{inputenc}
\usepackage{amsthm}
\usepackage{amsmath}
\usepackage{amssymb}

\newtheorem{theorem}{Theorem}[section]
\newtheorem{assumption}[theorem]{Assumption}
\theoremstyle{definition}
\newtheorem{lemma}[theorem]{Lemma}
\newtheorem{corollary}[theorem]{Corollary}
\newtheorem{proposition}[theorem]{Proposition}
\newtheorem{remark}[theorem]{Remark}
\newtheorem{example}[theorem]{Example}

\newtheorem{definition}[theorem]{Definition}
\numberwithin{figure}{section}

\newcommand{\N}{\mathbb{N}}

\newcommand{\R}{\mathbb{R}}
\newcommand{\B}{\mathbb{B}}
\newcommand{\Rinf}{\mathbb{R}\cup \{+\infty \}}

\renewcommand{\AA}{A^{\hspace{-0.05cm}*} \hspace{-0.05cm} A}

\renewcommand{\S}{\mathbb{S}}
\newcommand{\T}{T_\lambda}
\newcommand{\spec}{\mbox{\rm spec}}

\newcommand{\eps}{\varepsilon}
\newcommand{\Loja}{\L ojasiewicz }
\newcommand{\doms}{\mbox{\rm dom}^* f}
\newcommand{\argmin}{\mbox{\rm argmin\,}}
\newcommand{\dom}{\mathop\mathrm{\rm dom}}
\newcommand{\dist}{\mbox{\rm dist\,}}

\DeclareMathOperator{\Ker}{Ker}
\newcommand{\rint}{\mbox{\rm ri\,}}
\newcommand{\srint}{\mbox{\rm sri\,}}
\newcommand{\rank}{\mbox{\rm rank\,}}
\newcommand{\prox}{\mbox{\rm prox}}
\newcommand{\proj}{\mbox{\rm proj}}
\newcommand{\sgn}{\mbox{\rm sgn}}
 
\newcommand{\kin}{{k\in\N}}
\newcommand{\nin}{{n\in\N}}

\renewcommand{\d}{{\rm d}}
\newcommand{\dt}{{\rm d}t}
\DeclareFontEncoding{FMS}{}{}
\DeclareFontSubstitution{FMS}{futm}{m}{n}
\DeclareFontEncoding{FMX}{}{}
\DeclareFontSubstitution{FMX}{futm}{m}{n}
\DeclareSymbolFont{fouriersymbols}{FMS}{futm}{m}{n}
\DeclareSymbolFont{fourierlargesymbols}{FMX}{futm}{m}{n}
\DeclareMathDelimiter{\VERT}{\mathord}{fouriersymbols}{152}{fourierlargesymbols}{147}

\usepackage[usenames, dvipsnames]{color}

\newcommand\firstpagefootnote[1]{%
  \begingroup
  \renewcommand\thefootnote{}\footnote{\hspace*{-1.8em}#1}%
  \addtocounter{footnote}{-1}%
  \endgroup
}

\begin{document}
\title{
\textsc{\Large Convergence  of the Forward-Backward algorithm: \\  Beyond the worst-case   with the help of geometry}
%
	}
\author{Guillaume Garrigos$^{1}$, Lorenzo Rosasco$^{2,3}$, and Silvia Villa$^4$
}

\date{\vspace*{0.5em}

\small
$\!^1$ LPSM, Universit\'e de Paris. 75205 Paris CEDEX 13, France.\\
$\!^2$ DIBRIS, Universit\`a degli Studi di Genova. Via Dodecaneso 35, 16146, Genova, Italy.\\
$\!^3$ LCSL, Istituto Italiano di Tecnologia and Massachusetts Institute of Technology. \\Bldg. 46-5155, 77 Massachusetts Avenue, Cambridge, MA 02139, USA.\\
$\!^4$ Dipartimento di Matematica, Universit\`a degli Studi di Genova. Via Dodecaneso 35, 16146, Genova, Italy.
\vspace*{-2em}
}

\maketitle

\firstpagefootnote{%
\textbf{Contact:} \quad 
G. Garrigos {\ttfamily{garrigos@lpsm.paris}} \quad
L. Rosasco {\ttfamily{lrosasco@mit.edu}} \quad
S. Villa {\ttfamily{silvia.villa@unige.it}}
\vspace*{1em}
}

\firstpagefootnote{%
\textbf{Acknowledgements:} 
This material is supported by the Center for Brains, Minds and Machines, funded by NSF STC award CCF-1231216, and the Air Force project FA9550-17-1-0390. L. Rosasco acknowledges the financial support of the Italian Ministry of Education, University and Research FIRB project RBFR12M3AC. S. Villa is supported by the INDAM GNAMPA research project 2017 Algoritmi di ottimizzazione ed equazioni di evoluzione ereditarie.
}

\begin{abstract}
We provide a comprehensive study of the convergence of the forward-backward algorithm under suitable geometric conditions, such as  conditioning or \Loja properties. 
These geometrical notions are usually local by nature, and may fail to describe  the fine geometry of objective functions relevant in inverse problems and signal processing, that have a nice behaviour on  manifolds, or  sets open with respect to a weak topology. 
Motivated by this observation, we revisit those geometric notions over arbitrary sets.
In turn, this allows us to present several new results as well as  collect in a unified view a variety of  results scattered in the literature.
Our contributions include the analysis of infinite dimensional convex minimization problems, showing the first \Loja inequality for a quadratic function associated to a compact operator, and the derivation of new linear rates for problems arising from  inverse problems with low-complexity priors.
Our approach allows to establish unexpected connections between geometry and a priori conditions in inverse problems, such as source conditions, or restricted isometry properties.
\end{abstract}

\bigskip

\section{Introduction}

Splitting algorithms based on first order descent methods are widely used to solve high dimensional convex optimization problems in  signal and image processing \cite{ComPes11}, compressed sensing \cite{DauDefDem04}, and machine learning \cite{MesRosSan10}. 
Their main advantage is  their simplicity and   complexity independent of the dimension of the problem. {The worst case convergence rates of these methods have been intensively investigated in the last  twenty years. The simplest example is the  gradient method applied to a smooth convex function, which is known to  converge in values as $o(n^{-1})$ \cite{DavYin14,Sal16}. Analogous results are known for the forward-backward splitting algorithm.  We refer to these results as {\em worst case} since no particular assumption is made on the objective function aside from convexity and  existence of a solution. Note that these rates are sharp, meaning that 
there are functions for which these rates are arbitrarily accurate. } Clearly such a large class of convex functions allows for functions with wild behaviors around the minimizers \cite{BolDanLeyMaz10}, behaviors that  might hardly  appear in practice.  It is then natural to ask whether improved rates can be proved under further regularity assumptions. 

\newpage

\noindent \textbf{Previous work on optimization rates with geometry.}
One classical geometrical assumption is strong convexity, which indeed guarantees linear convergence rates \cite{Gol62,SchLerBac11}.  
In practice, strong convexity is often too restrictive, and  one would wish to relax it, while retaining fast  rates. 
{A relaxation of this condition is given by} geometric conditions that,  roughly speaking, describe convex functions $f\in \Gamma_0(X)$ that behave  like 
\begin{equation}\label{eq:appetizer}
x\mapsto \dist^p( x , \argmin f), 
\end{equation}
for some $p \geq 1$ and on some subset $\Omega \subset X$, which is typically a neighborhood of the  minimizers and/or a sub-level set.
 The intuition behind this kind of assumption required on a neghborhood of the solution is clear: the bigger is $p$, the more the function is ``{\em flat}'' around its minimizers, which in turns means that a gradient descent algorithm will converge slowly. 
The idea of exploiting  geometric conditions to derive convergence rates has a long history dating back to \cite{Pol63,Roc76}, 
and plenty of similar convergence rates results have been derived  under different yet related geometrical properties. 

{The optimization community focused on  several different but related geometric assumptions, namely the $p$-conditioning, the $p$-metric subregularity and the $p$-\Loja properties (see Section \ref{S:Geometry} for their definitions).
The first\footnote{If we discard the ``classic'' strong convexity assumption.} result exploting geometry to derive fast convergence rates dates back to Polyak \cite[Theorem 4]{Pol63}, showing that the gradient method converges linearly (in terms of the values and iterates) when the objective function verifies the $2$-{\L}ojasiewicz inequality.
Improved  convergence rates for first-order descent methods were then obtained in \cite{Roc76},  considering notions slightly stronger than $p$-metric subregularity, and proving  finite convergence of the proximal algorithm for $p=1$, and linear convergence for $p=2$.
These results are  improved and extended in \cite{Luq84},  analyzing  for the first time convergence rates for the iterates of the proximal algorithm  using  metric subregularity for  general $p \in [1, + \infty[$.
The results in  \cite{Luq84} recover those in  \cite{Roc76}  (see also \cite{Spi85,Spi87}),
but also derive superlinear rates for $p \in \left]1,2\right[$, and sublinear rates for $p>2$. 
Roughly speaking, the results in \cite{Luq84} show that  the bigger is $p$ the slower is the algorithm. 
A related notion, nowadays called the Luo-Tseng  error bound condition, has been considered in the seminal paper \cite{LuoTse93}, and implies the linear convergence of several first order methods.
Recently, this condition has been shown to be equivalent to 2-conditioning \cite{DruLew16,LiPon17}.
In the early 90's, some attention was devoted to the study of $p$-conditioned functions, in  particular for $p=1$ (some authors call this property  superlinear conditioning, 
sharp growth or sharp minima property).
In this context,  \cite{Fer91,Lem92,BurFer93} showed that the proximal algorithm terminates after a finite number of iterations.
For $p=1$, Polyak \cite[Theorem 7.2.1]{Pol}  obtained the finite termination for the projected gradient method.
The $2$-conditioning was also used  to obtain linear rates for the proximal algorithm  in \cite{Li95}. 
In \cite{AttBol09},  it was  observed that the $p$-{\L}ojasiewicz property could be used to derive precise rates for the iterates of the  proximal algorithm. 
The authors obtain finite convergence when $p=1$, linear rates when $p \in \left]1,2\right]$, and sublinear rates when $p \in \left]2, + \infty\right[$. 
Similar results can be found in \cite{AttBolRedSou10,MerPie10}. 
Such convergence rates for the iterates have been extended to the forward-backward algorithm (and its alternating versions) in \cite{BolSabTeb13}, and similar rates also hold for the convergence of the values in \cite{ChoPesRep14,FraGarPey15}. 
More recently, various papers focused on conditions  equivalent (or stronger) to the $2$-conditioning  to derive linear rates \cite{Lev09,LiaFadPey14,DruMorNgh14,LiuWriReBitSri15,DruLew16,KarNutSch16}. 
Some effort has also been made to show that the \Loja property and conditioning are equivalent \cite{BolDanLeyMaz10,BolNguPeySut15}, and to relate it to other error bounds appearing in the literature \cite{KarNutSch16}.
See also \cite{NecNesGli15} for a refined analysis of linear rates for the projected gradient algorithm under conditions that interpolate between strong convexity and $2$-conditioning (see also Subsection \ref{SS: linear rates}). 
}

\medskip

{\noindent\textbf{A key observation.}
Our study  starts from a basic observation which allows a number of developments. Indeed, motivated by several relevant examples described in Section~\ref{S:Linear inverse problems}, we  require  condition~\eqref{eq:appetizer} to hold on an arbitrary set $\Omega$, which in general
is neither a neighborhood of the solution, nor a sublevel set. This extension allows to establish a connection with modeling assumptions considered in different contexts and  unveil their role in optimization.  As we explain below, modeling assumptions, such as source conditions in inverse problems \cite{EngHanNeu} or the restricted injectivity property in sparse recovery \cite{Can08}, correspond to conditioning assumptions on specific subsets. This ensures global  convergence rates  for the forward-backward algorithm that are faster compared to those given by a worst case analysis and indeed  often observed in practice. }

\medskip

{\noindent \textbf{Geometry and inverse problems.}
As a first example of the importance of considering arbitrary sets $\Omega$ to define geometrical properties, consider linear inverse problems $Ax = y$ for which the operator $A$ is an infinite dimensional compact operator, making the problem severely ill-posed. 
A common modeling assumption is to suppose that the minimal norm solution of the problem satisfies a \textit{source condition}, which can be seen as a measure of its regularity (see Section \ref{SS:least squares in Hilbert spaces} for a definition). Under this condition, 
it is shown that the sublinear rate of the gradient algorithm is faster than the worst case one \cite{EngHanNeu}. However, such a behavior 
cannot be apparently explained in terms of  classical geometrical conditions satisfied by the least squares function: indeed, it was shown in \cite{HarJen11} that such a least squares function cannot verify any \Loja inequality \eqref{eq:appetizer} in a neighborhood of its minimizers. On the contrary,thanks to the extension of the definition considered in this paper, we show that geometric assumptions are indeed
satisfied, but only on specific subsets.
More precisely, we show in  Theorem \ref{T:geometry of least squares} that the source condition guarantees that the least squares $\Vert Ax-y\Vert^2$ is $p$-\Loja ($p>2$) on a dense affine subspace having empty interior.
This allows therefore to explain the faster global rates of the gradient algorithm which are typically observed in this context.

As a second example, consider linear inverse problems with a low-complexity prior, such as sparse inverse problems.
For these problems, the restricted injectivity condition \cite{Can08} is a key modeling assumption to guarantee stable recovery: it means that, even if a linear measurement is corrupted by noise, we can hope to reconstruct an approximated solution by solving a regularized optimization problem.
In Section \ref{SS:sparse inverse problems}, we show that this assumption implies a $2$-conditioning of the problem over a (nonconvex) cone of sparse vectors.
Since this set is \text{active}, in the sense that it is reached by the algorithm after a finite time, it immediately gives us  asymptotic linear rate of the algorithm.
For problems with more general low-complexity priors the situation is similar: an active set will be identified by the iterates of the algorithm, and we show that restricted injectivity condition on the tangent cone to this active set induces a $2$-conditioning of the problem on this set.
Depending on the applications or on the hypothesis made on the problem, this set can be a low-dimensional manifold, or a set with less structure, and can be computed within the partial smoothness framework \cite{HarLew04} or the mirror stratification one \cite{FadMalPey18}.}

\medskip

{\noindent \textbf{Paper contents.}
Motivated by the estimation problems presented in Section~\ref{S:Linear inverse problems}, the goal of this paper is to provide  a comprehensive study of the convergence rates of the forward-backward algorithm for convex minimization problems satisfying geometric conditions \textit{on arbitrary sets}.  We  collect in a unified view a variety of  results scattered in the literature, and we extend them to this more general setting. In addition, we  derive several novel results along the way. The paper is organized as follows.}

{After reviewing and discussing worst-case convergence results for the  forward-backward algorithm in Section~\ref{S:classic FB},  we give in Section~\ref{S:Geometry} the definition of  different  geometric conditions for a proper convex lower semicontinuous function $f$: $p$-conditioning, $p$-metric subregularity, and $p$-\Loja property on general subsets $\Omega \subset X$, rather than sublevel sets or open sets, as typically done in the literature.  
We show that those geometrical notion are equivalent, provided that the set $\Omega$ is stable by the semigroup generated by $\partial f$  (see Proposition~\ref{P:equivalence geometrical notions}).
Since establishing $p$-conditioning of a function may be hard in general,  
we provide two sum rules for conditioned functions in Theorem \ref{T:sum rule polyhedral} and Theorem~\ref{T:sum rule}. 
The first one 
establishes that if a strictly convex function remains $p$-conditioned under linear perturbations, then it is also $p$-conditioned  under convex perturbation.
The second one gives conditions under which the sum of two conditioned functions are conditioned. It allows us to show in particular that the ROF model (minimization of the total variation and the Kullback-Leibler divergence) is $2$-conditioned on every bounded set.

Section~\ref{S:CV rates for Forward Backward}  exploits the $p$-{\L}ojasiewicz property on general sets to study   the convergence  of the forward-backward algorithm.
In Theorem \ref{T:CVKL discrete}, we recover and extend results from the literature, getting finite / superlinear / linear / sublinear convergence rates, depending on the value of $p \in [1,+\infty[$ to our more general setting.
Along the way, we extend the sharp superlinear rate known for the proximal method to the Forward-Backward algorithm. 
In addition, our 
approach allows to derive in a unified setting both nonasymptotic/global and asymptotic/local convergence results, see Corollaries~\ref{T:CV on invariant sets} and \ref{P:capture result for local omega}.  
We go beyond the classical analysis by introducing a $p$-{\L}ojasiewicz property with $p$ taking \textit{nonpositive} values.
This allows to study convex functions being bounded from below but with no minimizers, a case which has drawn little attention so far, but which can arise for instance in function approximation \cite{Dev86} or in statistical learning  theory \cite[Theorem 9]{DevRosVer06} (see also Section \ref{SS:least squares in Hilbert spaces}).
For such ill-posed problems, we derive new and sharp sublinear rates for the values in Theorem \ref{T:CVKL discrete rates negative p}, interpolating between $o(n^{-1})$ and $o(1)$.
We further show in Section \ref{SS: linear rates} that the $2$-conditioning is essentially equivalent to the linear convergence of the forward-backward algorithm, illustrating  the importance of this notion for convergence rate analysis.}

{In Section~\ref{S:Linear inverse problems}, we apply the aforementioned results to optimization problems arising from inverse problems, and discuss the interaction between geometry and modeling assumptions. The key results of this section are Theorem~\ref{T:geometry of least squares} and Theorem~\ref{T:injective Hessian implies 2 conditioning}. Theorem~\ref{T:geometry of least squares} establishes that classical source conditions in inverse problems guarantee the \Loja property on special sets, and therefore give better convergence rates of the gradient method with respect to worst case ones. Theorem~\ref{T:injective Hessian implies 2 conditioning} says that if we have an a priori assumption about the minimizer, which is assumed to belong to a set $C$, then a restricted  injectivity property of the Hessian of the smooth component of the objective function implies that $f$ is $2$-conditioned on this set $C$ around the minimizer. This  guarantees asymptotic linear rates for forward-backward when combined with Corollary \ref{T:partial smoothness rates}.}

\section{The forward-backward algorithm: notation and background}\label{S:classic FB}
\subsection{Notation and basic definitions}\label{SS:notations and materials}

We recall a few classic notions and introduce some notation.
Throughout the paper $X$ is a Hilbert space. 
Given $\Omega \subset X$, we note $\mbox{int}~\Omega$ and $\mbox{cl}~\Omega$ its interior and closure.
We say that $\Omega$ is a cone, if $\Omega = ]0,+\infty[\Omega$.
We note $\mbox{cone}(\Omega)$ (resp. $\mbox{span}(\Omega)$) the smallest cone (resp. linear subspace) in $X$ containing $\Omega$.
Let $x\in X$, $\delta\in \left]0,+\infty\right[$, and let $\B_X(x,\delta)$ and $\overline{\B}_X(x,\delta)$ 
denote respectively the open and closed  balls  of radius $\delta$ centered
at $x$. We also use $\B_X$ and $\overline{\B}_X$ to denote $\B_X(0,1)$ and $\overline{\B}_X(0,1)$, and $\mathbb{S}_X$ to denote the unit sphere $\overline{\B}_X\setminus \mathbb{B}_X$. 
The distance of  $x\in X$ from a set $\Omega\subset X$ is  $\dist(x,\Omega)=\inf\{\|x-y\| \colon y\in \Omega \}$, and $\Vert \Omega \Vert_\_$ stands for $\dist(0,\Omega)$, so, in particular $\Vert \emptyset \Vert_\_= +\infty$.
{If $\Omega$ is closed and convex, $\proj(x,\Omega)$ is the projection of $x$ onto $\Omega$, and the relative interior and the strong relative interior of $\Omega$ are respectively defined as \cite[Definition 6.9]{BauCom}: 
$
\rint \Omega = \{x \in \Omega \ | \ \mbox{cone}(C-x) = \mbox{span}(C-x) \}$, $\srint = \{x \in \Omega \ | \ \mbox{cone}(C-x) = \mbox{cl}~\mbox{span}(C-x) \}$.
Given a bounded linear operator $A$ between two Hilbert spaces, its \textit{spectrum}, noted $\spec(A)$, is the set of spectral values
$\lambda \in \mathbb{R}$ such that $A - \lambda I$ is not boundedly invertible.
We also note $\spec^*(A) := \spec(A) \setminus \{0\}$.
The set of \textit{singular values} of $A$, noted $\sigma(A)$, is defined as $\sigma(A) := \sqrt{\spec^*(AA^*)}$, and we note $\sigma_{inf}(A) := \inf \sigma(A)$.}
Let $\Gamma_0(X)$ be the class of convex, lower semi-continuous, and proper functions from $X$ to $\left]-\infty,+\infty\right]$.
For $f \in \Gamma_0(X)$ and $x \in X$, $\partial f(x) \subset X$ denotes the (Fenchel) subdifferential of $f$ at $x$ \cite[Definition 16.1]{BauCom},
and $\dom f$ (resp. $\dom \partial f$) denotes the effective domain of $f$ (resp. of $\partial f$). Moreover, $f^*$ is the Fenchel
conjugate of $f$, namely $f^*(v)=\sup_{x\in X} \langle x,v\rangle-f(x)$ for all $v \in X$.
We introduce the shorthand notation $\doms := \dom f \setminus \argmin f$.
We also introduce the following notation for the (strict) sublevel sets of $f \in \Gamma_0(X)$:
for every $r \in ]-\infty,+\infty]$, $[f<r]:=\{x \in X \ | \ f(x)  < r \}$.

\bigskip

The following assumption will be made throughout this paper. 
\begin{assumption} 
\label{ass:H}
Let $X$ be a Hilbert space,  $g\in \Gamma_0(X)$, and  $h\colon X\to \mathbb{R}$ be differentiable and convex, with  $L$-Lipschitz continuous gradient for some 
$L\in\left]0,+\infty\right[$ and  set $ f = g+h$.
\end{assumption}
Splitting methods, such as the forward-backward algorithm, are extremely popular for minimizing an objective function as in Assumption \ref{ass:H}.
To have an implementable procedure, we  implicitly assume 
that the proximal operator of $g$ can be easily computed (see e.g. \cite{ComPes11}):
\begin{equation}\label{D:prox}
(\forall \lambda > 0)(\forall x \in X)  \quad \prox_{\lambda g}(x)= \underset{u \in X}{\argmin} \left\{ g(u) + \frac{1}{2\lambda} \Vert u - x \Vert^2 \right\}.
\end{equation}
Remembering Assumption \ref{ass:H} is in force, we introduce the Forward-Backward (FB) map for $\lambda \in ]0,2L^{-1}[$:
\begin{equation}\label{D:Forward Backward mapping}
T_\lambda : x \in X \longmapsto T_\lambda x:=\prox_{\lambda g} (x - \lambda \nabla h(x)) \in X  , 
\end{equation}
so that the FB algorithm can be simply written as $x_{n+1}=T_\lambda x_n$.

\subsection{The Forward-Backward algorithm: worst-case analysis}

The following theorem collects known results about the convergence of the FB algorithm. 
This is a  ``worst-case'' analysis, in the sense that it holds for every  $f\in\Gamma_0(X)$ satisfying Assumption~\ref{ass:H}.
The main goal of Section~\ref{S:CV rates for Forward Backward}  is to show how these results can be improved taking
into account the geometry of $f$ at its infimum.

\begin{theorem}[Forward-Backward - convex case]\label{T:CV FB}
Suppose that Assumption~\ref{ass:H} is in force, and let $(x_n)_\nin$ be generated by the FB algorithm with $\lambda \in ]0,2L^{-1}[$.
Then:
\begin{enumerate}[i)]
	\item\label{T:CV FB:i} \textit{(Descent property)} The sequence $(f(x_n))_\nin$ is   decreasing, and converges to $\inf f$.
	\item\label{T:CV FB:ii} \textit{(F\'ejer property)} For all $\bar x \in \argmin f$, the sequence $\left(\Vert x_n - \bar x \Vert \right)_\nin$ is decreasing.
	\item\label{T:CV FB:iii} \textit{(Boundedness)} The sequence $(x_n)_\nin$ is bounded if and only if $\argmin f$ is nonempty.
\end{enumerate}
Suppose in addition that $f$ is bounded from below. Then
\begin{enumerate}[resume*]
	\item\label{T:CV FB:iv} \textit{(Subgradients convergence)} 
	The sequence $\left( \Vert \partial f(x_n) \Vert_\_ \right)_\nin$  converges decreasingly  to zero, with $
	\Vert \partial f(x_{n+1}) \Vert_\_^2 = O\left( f(x_n) - \inf f \right).$
\end{enumerate}
Moreover, if $\argmin f \neq \emptyset$, we have:
\begin{enumerate}[resume*]
	\item \label{T:CV FB:v} \textit{(Weak convergence)} The sequence $(x_n)_\nin$ converges weakly to a minimizer of $f$.
	\item \label{T:CV FB:vi} \textit{(Global rates for function values)} For all $\nin$, 
	\begin{equation*}
\hspace{-0.4cm}	f(x_n) - \inf f \leq {C} \frac{\dist(x_0,\argmin f)^2}{2\lambda n}, \text{ with } C=
	\begin{cases}
		1 & \text{ if } \lambda \leq L^{-1}, \\
		1+2(\lambda L -1)(2-\lambda L)^{-1} & \text{ otherwise}.
	\end{cases}
	\end{equation*}
	\item\label{T:CV FB:vii}\label{Asymptotic rates for values} \textit{(Asymptotic rates for function values)} When $n \to +\infty$, $f(x_n) - \inf f = o\left(n^{-1}\right).$
\end{enumerate}
\end{theorem}

\noindent Theorem \ref{T:CV FB} collects various convergence results on the FB algorithm.
Item \ref{T:CV FB:i} appears in \cite[Theorem 3.22]{Sal16} (see also \cite{Gul91}).
Item \ref{T:CV FB:ii} is a consequence of the nonexpansiveness of the FB map 
(see \eqref{D:Forward Backward mapping}) \cite[Lemma 3.2]{Lem96}.
Item \ref{T:CV FB:iii}, which is a consequence of Opial's Lemma \cite[Lem. 5.2]{Pey}, can  be found in \cite[Theorem 3.12]{Sal16}.
Item \ref{T:CV FB:iv} follows from Lemma~\ref{L:estimates for the forward-backward}.\ref{L:estimates for the forward-backward:gradients} in the Annex.
Item \ref{T:CV FB:v} is also a consequence of Opial's Lemma, see \cite[Proposition 3.1]{Lem96}.
Items \ref{T:CV FB:vi} and \ref{T:CV FB:vii} are proved in \cite[Theorem 3]{DavYin14} (see also \cite[Proposition 2]{BreLor08b} and \cite[Theorem 3.1]{BecTeb09}).

\begin{remark}[Sharpness of the results in the worst-case]\label{R:optimality of rates worst case}
The convergence results in Theorem \ref{T:CV FB} are  sharp, in the following sense.
First, the iterates may not converge strongly: see \cite{Bai78,Gul91} for a counterexample in $\Gamma_0(\ell^2(\N))$.
Even in finite dimension, no sublinear rates should be expected for the iterates. 
To see this, apply the proximal algorithm to the function  $x\in \R\mapsto f_p(x)=\vert x \vert^p$, whose unique minimizer is zero.
When $p\in \left]2,+\infty\right[$, there exists a constant $C_p >0$ depending on $(\Vert x_0 \Vert, \lambda, p)$ such that (see e.g. the discussion following \cite[Proposition 2.5]{MerPie10}, or Lemma \ref{L:lower bounds proximal}):
\begin{equation}\label{E:rates for norm to the p}
(\forall n \geq 1) \quad \vert x_n \vert \geq C_p n^{-1/(p-2)}, \quad\text{where} \lim_{p\to +\infty} \frac1{p-2}=0.
\end{equation}
The estimate \eqref{E:rates for norm to the p}  also provides a lower bound for the rates on the objective values:
\begin{equation}\label{E:rates for the norm to the p values}
f_p(x_n) - \inf f_p \geq C_p^p n^{-p/(p-2)}. 
\end{equation} 
The above lower bounds imply that the rate in Theorem \ref{T:CV FB:vii} cannot be improved into a rate $O(n^{-\delta})$, for some $\delta >1$, because we can always find a $p$ large enough verifying $p/(p-2) > \delta$.
It also means that no polynomial rates can be expected for $\Vert x^n - \bar x \Vert$.
This fact was also observed in \cite[Theorem 12]{DavYin14} on an infinite dimensional counterexample.
When $f$ is bounded from below, but has no minimizers, the values $f(x_n) - \inf f$ 
go to zero but no rates can be obtained in general.
To see this, consider for any  $\alpha >0$ the function $f_\alpha \in \Gamma_0(\mathbb{R})$ defined by 
\begin{equation}\label{E:counter example sublinear rates}
f_\alpha : \mathbb{R}\rightarrow ]-\infty,+\infty] \ : \ f_\alpha(x) = \vert x \vert^{-\alpha} \ \text{ if }  \ x <0, \ +\infty \ \text{ otherwise.}
\end{equation}

{If $(x_n)_\nin$ is  obtained by applying the proximal algorithm to this function, then (see Lemma \ref{L:lower bounds proximal}) there exists $C_\alpha>0$ such that:
\begin{equation}\label{E:rates for the norm to the -p values}
 f_\alpha(x_n) - \inf f_\alpha \geq C_\alpha^{-\alpha} n^{-\alpha/(2+\alpha)}, \
\text{ where }
\lim\limits_{\alpha \to 0} \frac{\alpha}{2+ \alpha} = 0 
\
\text{ and } 
\lim\limits_{\alpha \to +\infty} \frac{\alpha}{2+ \alpha} = 1. 
\end{equation}
Observe that this lower bound on the objective function values implies that the convergence for those functions is slower than the usual 
$O(n^{-1})$ rate obtained in Theorem \ref{T:CV FB}.\ref{T:CV FB:vi}. It also shows that no polynomial rates can be proven for the values when ${\rm{argmin}}~f =\emptyset$.}
\end{remark}

\section{Identifying the geometry of a function}\label{S:Geometry}
\subsection{Definitions}\label{SS:geometry definitions}
In this section we introduce the main geometrical concepts that will be used throughout the paper to derive precise rates for the FB method.
Roughly speaking, these notions characterize functions which behave like \eqref{eq:appetizer}  on an arbitrary set $\Omega \subset X$.
\begin{definition}
\label{D:geometric notions}
Let $p \in [1,+\infty[$, let $f \in \Gamma_0(X)$ with $\argmin f \neq \emptyset$, and $\Omega \subset X$.
We say that:
\begin{enumerate}[i)]
	\item  $f$ is \textit{$p$-conditioned} on $\Omega$ if there exists a constant $\gamma_{f,\Omega} >0$ such that:
\begin{equation*}\label{E:growth inequality}
\forall x \in  \Omega \cap \dom f , \quad \frac{\gamma_{f,\Omega}}{p} \dist(x,\argmin f )^p \leq  f(x) - \inf f.
\end{equation*}
	\item $\partial f$ is \textit{$p$-metrically subregular} on $\Omega$ if there exists a constant $\gamma_{\partial f,\Omega} >0$ such that:
	\begin{equation*}\label{E:metric subregularity}
\forall x \in  \Omega\cap \doms , \quad {\gamma_{\partial f,\Omega}} \dist(x,\argmin f )^{p-1} \leq \Vert \partial f(x) \Vert_\_ .
	\end{equation*}
	\item $f$ is \textit{$p$-{\L}ojasiewicz} on $\Omega$ if there exists a constant $c_{f,\Omega}>0$ such that:
\begin{equation*}
\label{E:Lojasiewicz inequality}
\forall x  \in \Omega \cap \doms ,   \quad (f(x) - \inf f)^{1- \frac{1}{p}} \leq c_{f,\Omega} \Vert \partial f(x) \Vert_\_. 
\end{equation*}
\end{enumerate}
We will refer to these notions as  \textit{global} if $\Omega=X$, and as \textit{local} if $\Omega= \B_X(\bar x ; \delta) \cap [f<r]$ for some $\bar x \in \argmin f,$  and $\delta \in ]0,+\infty]$, $r \in ] \inf f, +\infty]$.
\end{definition}

The notion of conditioning, introduced in \cite{Vai70,Zol78}, is a common tool in the optimization and regularization literature \cite{AttWet93,Pen96,Lem98,Zal,BolNguPeySut15}.
It is also called the  \textit{growth condition}  \cite{Pen96}, and it is strongly related to the notion of Tikhonov wellposedness \cite{DonZol93}.
The $p$-metric subregularity coincides with metric subregularity of the subdifferential at the origin, and it is less used, generally defined 
for $p=1$ or $2$  with $\Omega$ equal to a neighborhood of a specific minimizer \cite{DonLewRoc03,Lev09}. It is also called upper Lipschitz
continuity at zero of $\partial f^{-1}$ in \cite{CorJouZal97}, or inverse calmness \cite{DonRoc09}.
The {\L}ojasiewicz property goes back to \cite{Loj63}, and was initially designed as a tool to guarantee the convergence of 
trajectories for the gradient flow of analytic functions, before its recent use in convex and nonconvex optimization.
It is generally presented with a constant $\theta \in [0,1]$ which is equal, in our notation, to  $1-1/p$ \cite{Loj63,AbsMahAnd05,BolDanLew07,BolNguPeySut15}, 
or  $1/p$ \cite{MerPie10,HarJen11,FraGarPey15}.
In the remark below we explain the main difference between our definition and the one usually considered in the literature.
{\begin{remark}
There is a subtle but crucial difference in the terminology used in Definition~\ref{D:geometric notions}  with respect to the one commonly used for the {\L}ojasiewicz property. 
It is usually said that a function has the \L ojasiewicz property at $\bar{x}$ if there exist $\delta>0$, $c>0$, and $r>\inf f$ such that $f(x)-f(\bar{x})\leq c\|\partial f(x)\|_{- }$ holds on $\Omega=\B_X(\bar x ; \delta) \cap [f<r]$. If the latter property holds for every $\bar{x}\in S\subset X$, the function is said to have the  {\L}ojasiewicz property on $S$.
This is a different requirement with respect to the one in Definition~\ref{D:geometric notions}.
Indeed, we require the inequality to hold uniformly on $\Omega$, while the above definition must hold locally around every point of interest in a given set, and typically only allows for asymptotic convergence rates (see Corollary~\ref{P:capture result for local omega}). 
This change of viewpoint is motivated by the fact that for many \textit{convex} functions, we have more than just a local information about the geometry (see Sections \ref{SS:sum rule} and \ref{S:CV rates for Forward Backward}). 
More importantly, it is actually necessary for the analysis of the problems discussed in Section \ref{S:Linear inverse problems}, which motivated this paper. 
Beyond that, it also allows to understand in a unified framework both global (Corollary~\ref{T:CV on invariant sets}) and local 
(Corollary~\ref{P:capture result for local omega}) convergence rates.
\end{remark}}
The notions introduced in Definition~\ref{D:geometric notions} are closely related to each other. 
Indeed, for convex functions, $p$-conditioning implies metric subregularity, which implies
the {\L}ojasiewicz property.
Under some additional assumptions, it is possible to show that the reverse implications hold.
For instance, metric subregularity implies conditioning when $\Omega= \argmin f + \delta \B_X$, $\delta >0$ \cite[Theorem 4.3]{ZhaTre95}.
Similar results can also be found in \cite{AraGeo08,AzeCor14,DruMorNgh14,DruIof15}, and \cite[Theorem 5.2]{CorJouZal97} (for $\Omega=X$).
Also, it is shown in \cite[Theorem 5]{BolNguPeySut15} that the local {\L}ojasiewicz property implies local conditioning. 
The next result, proved in Annex~\ref{SS:equivalence geometric notions}, extends the mentioned ones, and states 
the equivalence between conditioning, metric subregularity, and {\L}ojasiewicz property on $\partial f$-invariant sets 
(see Definition~\ref{D:invariant set continuous dynamic} in Annex~\ref{SS:equivalence geometric notions}). 

\begin{proposition}
\label{P:equivalence geometrical notions}
Let $p \in [1,+\infty[$, let $\Omega \subset X$, and let $f \in \Gamma_0(X)$ be such that $\argmin f \neq \emptyset$. Consider the following properties:
\begin{enumerate}[i)]
\item
\label{P:equivalence geometrical notions:i}
 $f$ is $p$-conditioned on $\Omega$ ,
\item
\label{P:equivalence geometrical notions:ii}
 $\partial f$ is $p$-metrically subregular  on $\Omega$,
\item
\label{P:equivalence geometrical notions:iii}
 $f$ is $p$-{\L}ojasiewicz on $\Omega$.
\end{enumerate}
Then \ref{P:equivalence geometrical notions:i} $\implies$ \ref{P:equivalence geometrical notions:ii} $\implies$ \ref{P:equivalence geometrical notions:iii}.
One can  respectively take $\gamma_{\partial f,\Omega}=\gamma_{f,\Omega}/p$ and $c_{f,\Omega}= \gamma_{\partial f,\Omega}^{-1/p}$.
Assuming in addition that $\Omega$ is $\partial f$-invariant, we also have
\ref{P:equivalence geometrical notions:iii} $\implies$ \ref{P:equivalence geometrical notions:i} with $\gamma_{f,\Omega}=c_{f,\Omega}^{-p} p^{1-p}$.
\end{proposition}

The two next propositions show that these geometric notions are stronger when $p$ is smaller, and  are  meaningful only on sets containing minimizers (their proof follow directly from  Definition \ref{D:geometric notions} and are left to the reader).

\begin{proposition}
\label{P:hierarchy of conditionings with p}
Let $f \in \Gamma_0(X)$ be such that $\argmin f \neq \emptyset$, $\Omega \subset X$, and $p' \geq p \geq 1$.
\begin{enumerate}[i)]
	\item\label{P:hierarchy of conditionings with p:i} If $f$ is $p$-conditioned (resp. $\partial f$ is $p$-metrically subregular) on $\Omega$, then  $f$ is $p'$-conditioned (resp. $\partial f$ is $p'$-metrically subregular) on $\Omega \cap \delta\B_X$ for any $\delta \in ]0, + \infty[$.
	\item\label{P:hierarchy of conditionings with p:ii} If $f$ is $p$-{\L}ojasiewicz on $\Omega$, then $f$ is $p'$-{\L}ojasiewicz on $\Omega \cap [f<r]$ for any $r > \inf f$.
\end{enumerate}
\end{proposition}

\begin{proposition}
\label{P:trivial conditioning without minimizers}
Let $f \in \Gamma_0(X)$ be such that $\argmin f \neq \emptyset$.
If $\Omega \subset X$ is a weakly compact set  for which $\Omega \cap \argmin f = \emptyset$, then $f$ is $p$-conditioned on $\Omega$ for any $p \in [1,+\infty[$.
\end{proposition}

\subsection{Examples}\label{SS:geometry examples}
In this section, we collect some relevant examples.

\begin{example}[Uniformly convex functions]\label{Ex:uniformly convex functions}
Suppose that $f\in \Gamma_0(X)$ is uniformly convex of order $p \in [2,+\infty[$ \cite[Definition 10.7]{BauCom}. Then, there exists $\gamma >0$ such that \cite[Corollary 3.5.11.iv]{Zal}:
\begin{equation*}
(\forall (x_1,x_2) \in \dom \partial f^2) (\forall x^*_1 \in \partial f(x_1)) \quad f(x_2) - f(x_1) - \langle x_1^*,x_2 - x_1 \rangle \geq \frac{\gamma}{p}\Vert x_2 - x_1 \Vert^p.
\end{equation*}
Such function is globally $p$-conditioned, with $\gamma_{f,X}=\gamma$,
and globally $p$-{\L}ojasiewicz, 
with $c_{f,X}=(1-1/p)^{1-1/p} \gamma^{-1/p}$ (see Lemma \ref{L:Loja for uniformly convex}).
In the strongly convex case, when $p=2$, the $2$-{\L}ojasiewicz inequality holds with the constant $c_{f,X}=1/\sqrt{2\gamma}$, which is sharp.
Examples of uniformly convex functions of order $p$ are $x \mapsto \Vert x \Vert^p$ \cite[Example 10.16]{BauCom}.
\end{example}

\begin{example}[Least squares]\label{Ex: least squares}
Let $A : X \rightarrow Y$ be a nonzero bounded linear operator between Hilbert spaces, and $f(x)=(1/2)\Vert Ax - y \Vert^2$, for some $y \in Y$.
Then, the conditioning, metric subregularity, and {\L}ojasiewicz properties, with $p=2$ and $\Omega=X$, are  equivalent to verify on $\Ker A^\perp$, respectively:
\begin{equation*}
 \gamma_{f,X} \Vert x \Vert^2 \leq \langle A^*Ax,x \rangle, \ \ \gamma_{\partial f,X} \Vert x \Vert \leq \Vert  A^* Ax \Vert, \  \text{ and }  \ \langle A^* Ax,x \rangle \leq 2c_{f,X}^2 \Vert A^* Ax \Vert^2.
\end{equation*}
If $\sigma_{\inf}(A^*A) >0$ holds, one can see that the above inequalities hold with
\begin{equation*}
\gamma_{f,X}=\gamma_{\partial f,X} = 1/(2c_{f,X}^2) = \sigma_{\inf}(A^*A),
\end{equation*}
 meaning in particular that $f$ is globally $2$-conditioned.
Since $\sigma_{\inf}(A^*A) >0$ is equivalent for $R(A^*A)$ to be closed (see Proposition \ref{P:closed range singular values}), it is in particular always true when $Y$ has finite dimension. 
If instead $\sigma_{\inf}(A^*A) =0$ holds, \cite[Theorem 2.1]{HarJen11} shows that $f$ cannot satisfy any local $p$-{\L}ojasiewicz property, for any $p\geq 1$.
This is for instance the case for infinite dimensional compact operators.
Nevertheless, we will show in Section~\ref{S:Linear inverse problems}, that the least squares always satisfies a $p$-{\L}ojasiewicz property on the so-called regularity sets, for any  $p > 2$.
\end{example}

\begin{example}[Convex piecewise polynomials]\label{Ex:convex piecewise polynomials are conditioned}
A convex continuous function $f: \R^N \rightarrow \R$ is said to be \textit{convex piecewise polynomial} if $\R^N$ can be partitioned 
in a finite number of polyhedra $P_1,...,P_s$ such that for all $i\in \{1,...,s\}$, the restriction of $f$ to $P_i$ is a convex polynomial, of degree $d_i \in \N$.
The degree of $f$ is defined as $\deg (f):= \max \{d_i \ | \ i\in \{1,...,s \}\}$. Assume $\deg (f)>0$.
Convex   piecewise polynomial functions are conditioned \cite[Corollary 3.6]{Li13}.
More precisely, for all $r > \inf f$, $f$ is $p$-conditioned on its sublevel set $\Omega=[f<r]$, with $p=1+(\deg (f)-1)^N.$
In general, the constant $\gamma_{f,\Omega}$ (which depends on $r$) cannot be explicitly computed.
This result implies that polyhedral functions ($\deg (f)=1$) are $1$-conditioned (in agreement with \cite[Corollary 3.6]{BurFer93}), and that convex piecewise quadratic functions ($\deg (f)=2$) are $2$-conditioned (in agreement with \cite[Theorem 2.7]{Li95}).
More generally, convex semi-algebraic functions are locally $p$-conditioned \cite{BolDanLewShi07}.
\end{example}

\begin{example}[L1 regularized least squares]\label{Ex:LASSO is 2 conditioned}
Let $f(x)=\alpha \Vert x \Vert_1 + (1/2)\Vert Ax-y \Vert^2$, for some linear operator $A : \R^N \rightarrow \R^M$, $y \in \R^M$ and $\alpha >0$.
As observed in \cite[Section 3.2.1]{BolNguPeySut15}, $f$ is  convex piecewise polynomial  of degree $2$, thus it is $2$-conditioned on every nonempty level 
set $\Omega=[f<r]$. The computation of the conditioning constant $\gamma_{f,\Omega}$ is rather difficult.
In \cite[Lemma 10]{BolNguPeySut15} an estimate of $\gamma_{f,\Omega}$ is provided, by means of Hoffman's bound \cite{Hof52}. Extensions of this result to the infinite dimensional setting can be found in \cite{GarRosVil20}.
\end{example}

\begin{example}[Regularized problems]\label{Ex:general regularization}
Let $X$ be an Euclidean space, 
$f(x):= g(x) + h(Ax)$, where $A : X \rightarrow \R^M$ is a linear operator, $g \in \Gamma_0(X)$, and $h \in \Gamma_0(\R^M)$ is a strongly convex $C^{1,1}$ function, and $\argmin f \neq \emptyset$.
Then $f$ is $2$-conditioned on any level set $\Omega=[f<r]$, for $r > \inf f$, if
\begin{enumerate}[i)]
	\item $g(x) = \Vert x \Vert_p$ with $p \in \left]1,2\right]$,  (see \cite[Corollary 2]{ZhoZhaSo15}),
	\item $g(x)=\Vert x \Vert_p^p$ with $p \in \left]1,2\right]$, (use \cite[Theorem 4.2]{DruLew16}; the details are left to the reader as an exercise, and can be checked in the Appendix),
	\item $g(x)=\Vert x \Vert_*$ is the nuclear norm of the matrix $x \in X$, provided the following qualification condition holds\footnote{We mention that this result was originally announced in \cite[Theorem 3.1]{HouZhoSoLuo13} without the qualification condition, but then corrected in \cite[Proposition 12 \& following remarks]{ZhoSo15}, in which the authors show that such condition is necessary.} (see \cite{ZhoSo15}): $\exists \bar x \in \argmin f$ such that $   -A^*\nabla h(A \bar x) \in \rint \partial \Vert \cdot \Vert_* (\bar x)$.
	\item $g$ is polyhedral (see \cite[Proposition 6]{ZhoSo15}).
\end{enumerate}
Note that in \cite{ZhoSo15,ZhoZhaSo15}, the authors do not prove directly that the functions are $2$-conditioned, but that they verify the so-called Luo-Tseng error bound, that is known to be equivalent to $2$-conditioning on sublevel sets \cite[Corollary 3.6]{DruLew16}.
Note also that in items ii-iv),  the strong convexity and $C^{1,1}$ assumptions on $h$ can be weakened (see \cite{ZhoSo15} and \cite[Theorem 4.2]{DruLew16}).
\end{example} 

\begin{example}[Distance to an intersection]
\label{Ex:distance to intersection}
Let $C,D$ be two closed convex sets in $X$ such that $C\cap D\neq\varnothing$, and for which the intersection is 
sufficiently regular, i.e. $0 \in \srint (C - D)$.  Let $f(\cdot)=\max\{\dist(\cdot,C),\dist(\cdot,D)\}$.
Clearly, $f \in \Gamma_0(X)$, and $\argmin f = C\cap D$.
Then $f$ is $1$-conditioned on bounded sets \cite[Theorem 4.3]{BauBor93}.
Let $p \in \left[1, + \infty\right[$. From  $\Vert \cdot \Vert_\infty \leq \Vert \cdot \Vert_p$, it follows that the function $x \mapsto \dist(x,C)^p + \dist(x,D)^p$ is $p$-conditioned on bounded sets.
The regularity condition $0 \in \srint (C-D)$ is not necessary if the two sets are polyhedral, as proved by Hoffman \cite{Hof52}.
\end{example}

{\begin{example}[Minimum of {\L}ojasiewicz functions]
If $f = \min_{i=1, \dots , m} f_i$, with $f_i \in \Gamma_0(\mathbb{R}^N)$ being continuous on its domain, and locally $p$-{\L}ojasiewicz at $\bar x \in {\rm{argmin}}~f$, then $f$ is locally $p$-{\L}ojasiewicz at $\bar x$ \cite[Theorem 3.1]{LiPon17}.
It is important to notice that this result do not need the $f_i$'s to be convex.
\end{example}}

The next section presents new sum rules for conditioned functions.

\subsection{A sum rule for $p$-conditioned functions}\label{SS:sum rule}
{Since verifying conditioning directly with the definition can be difficult, it is very
useful to establish which basic operations preserve conditioning. 
In this section we present two new
sum rules for conditioned functions in a setting where  $f = g+h\circ A$, where $g$ and $h$ are convex and $A$ is a bounded linear operator.
Theorem~\ref{T:sum rule polyhedral} states that if $g$ strictly convex and $p$-conditioned up to linear perturbations then also $f$ is  $p$-conditioned. 
Theorem~\ref{T:sum rule} provides an alternative where the assumption of strict convexity of $g$ is replaced by a stable conditioning assumption on $h$, 
which we formalise in the next definition, inspired by the terminology used in \cite{PolRoc96,DruMorNgh14,DruLew16}.}

{\begin{definition}
Let $f \in \Gamma_0(X)$, $\Omega \subset X$, and $p \in [1,+\infty[$. 
We say that $f$ is $p$-\textit{tilt-conditioned} if, for every $u \in X$, the tilted function $f + \langle u, \cdot \rangle$ has no minimizers, or is $p$-conditioned on $\Omega$.
\end{definition}}
{Note that a similar notion is already present in the literature: if $f$ is $p$-tilt-conditioned (in our sense) on every compact set, then it is \textit{firmly convex} in the sense of \cite[Definition 4.1]{DruLew16}.
}

\begin{example}[Tilt-conditioned functions]\label{R:tilt conditioned functions}
Many conditioned functions relevant for inverse problems are also tilt-conditioned:
\begin{itemize}
    \item The $1$-norm $\Vert \cdot \Vert_1$, and more generally every polyhedral function, are $1$-tilt-conditioned on Euclidean spaces \cite[Cor. 3.6]{BurFer93}.
    \item Convex piecewise polynomials of degree 2 are $2$-tilt-conditioned on their sublevel sets. 
    This is due to Example \ref{Ex:convex piecewise polynomials are conditioned} and the fact that this class of functions is stable up to linear perturbations.
    \item For the same reasons as above, $p$-uniformly convex functions are $p$-tilt-conditioned on $X$, for $p\geq 2$.
    \item {If $KL(x_1;x_2)$ denotes the Kullback-Leibler divergence between two vectors in $]0,+\infty[^N$, then the divergence $KL(x_1; \cdot)$ is $2$-tilt-conditioned on bounded sets.
    This result is new, and its proof can be found in Lemma \ref{L:kullback libler tilt conditioning}.}
    \item The nuclear norm is $2$-tilt-conditioned on bounded sets \cite[Proposition 11]{ZhoSo15}.
    \item See  \cite[Section 4]{DruLew16} for more examples and properties of $2$-tilt-conditioned functions on compact sets.
\end{itemize}
\end{example}

{In this first theorem, we show that if a strictly convex function remains conditioned up to linear perturbations, then it is also stable up to \textit{convex} perturbations:
\begin{theorem}[Sum rule involving a strictly convex tilt-conditioned function]\label{T:sum rule polyhedral}
Let $f=g+h\circ A$, where $g \in \Gamma_0(X)$, let $Y$ be a Hilbert space, $h \in \Gamma_0(Y)$ and $A:X \rightarrow Y$  a bounded linear operator.
Suppose that ${\rm{argmin}}~f \neq \emptyset$.
Let $\Omega\subset X$, and assume that:
\begin{enumerate}[label=\alph*)]
    \item\label{T:sum rule polyhedral:compatibility} the nondegeneracy condition $0 \in \srint \left(  \dom h - A(\dom g) \right)$ holds,
    \item\label{T:sum rule polyhedral:strict convex} 
    $g$ is strictly convex on its domain,
    \item\label{T:sum rule polyhedral:tilt conditioning} $g$ is $p$-tilt conditioned on $\Omega$ for some $p\in \left[1,+\infty\right[$.
\end{enumerate}
Then, $f$ is $p$-conditioned on $\Omega$. We have $\gamma_{f,\Omega} = \gamma_{\tilde g, \Omega}$, where $\tilde g=g + \langle \cdot, u \rangle$, for some $u \in X$.
\end{theorem}
\begin{proof}
Let $\bar x \in {\rm{argmin}}~f$;  Fermat's rule implies that $0 \in \partial f(\bar x)$.
Using assumption \ref{T:sum rule polyhedral:compatibility} with  \cite[Thm. 16.47]{BauCom}, we can write $0 \in \partial g(\bar x) + A^* \partial h(A\bar x)$.
Let $\bar v \in -\partial h(A\bar x)$ be such that $0 \in \partial g(\bar x) - A^* \bar v$, i.e., $\bar x \in \partial g^*(A^* \bar v)$.
Let $x \in \Omega \cap \dom f$, and set $\tilde g=g - \langle A^*\bar v, \cdot \rangle$.
Using the fact that linear forms are continuous, we can use again Fermat's rule together with a sum rule \cite[Thm. 3.30]{Pey} to write
\begin{equation}\label{srcsc1}
    u \in {\rm{argmin}}~\tilde g 
    \Leftrightarrow 
    0 \in \partial (g - \langle A^*\bar v, \cdot \rangle)(u) 
    =
    \partial g(u) - A^* \bar v
    \Leftrightarrow 
    A^*\bar v \in \partial g(u)
    \Leftrightarrow u \in \partial g^*(A^* \bar v),
\end{equation}
meaning that $\argmin \tilde g = \partial g^*(A^* \bar v) \neq \emptyset$.
It follows then from assumption \ref{T:sum rule polyhedral:tilt conditioning} that $\tilde g$ is $p$-conditioned on $\Omega$.
Moreover, because $g$ is strictly convex, we have $\partial g^*(A^* \bar v) = \{ \bar x \}$ \cite[Prop. 16.37.i]{BauCom}, and ${\rm{argmin}}~f = \{\bar x \}$ \cite[Cor 11.9]{BauCom}.
These facts mean that ${\rm{argmin}}~\tilde g = {\rm{argmin}}~f$.
We can now write the conditioning of $\tilde g$  evaluated at $x$, together with the convexity of $h$ (remember that $- \bar v \in \partial h(A \bar x)$):
\begin{eqnarray*}
g(x) & \geq & g(\bar x) + \langle A^* \bar v, x - \bar x \rangle + (\gamma_{\tilde g, \Omega}/p) \dist^p(x,{\rm{argmin}}~f), \\
h(Ax) & \geq & h(A\bar x) + \langle -\bar v, Ax - A\bar x  \rangle.
\end{eqnarray*}
Observe that we are allowed to use the conditioning of $\tilde g$ at $x$, because $x \in \dom f \subset \dom g = \dom \tilde g$.
Summing these two last inequalities gives
\begin{equation*}
 f(x) - \inf f \geq \frac{\gamma_{f,\Omega}}{p}  \dist^p(x,{\rm{argmin}}~f) ,
\end{equation*}
with $\gamma_{f,\Omega}:= \gamma_{\tilde g, \Omega}$, which concludes the proof.
\end{proof}}

{\begin{remark}[On the nondegeneracy condition \ref{T:sum rule polyhedral:compatibility} of Theorem \ref{T:sum rule polyhedral}]\label{R:compatibility condition}
This condition is very mild, and is satisfied under any of the following sufficient conditions (we note $\bar x$ a minimizer of $f$):
\begin{itemize}
    \item $h$ is continuous at $A\bar x$ (see \cite[Prop. 16.27 \& Prop. 6.19.vii]{BauCom}).
    \item $h$ has a full domain.
    \item $\dim Y <+ \infty$, $\bar x \in \qrint \dom g$ and $A\bar x \in \rint \dom h$ (see \cite[Def. 6.9 \& Prop. 6.19.ix]{BauCom}). 
    These inclusions hold for instance if $g$ and $h$ have open domains.
\end{itemize}
\end{remark}}
{Theorem \ref{T:sum rule polyhedral} is useful, but proves to be impractical when $g$ is not strictly convex, which typically happens when $g$ corresponds to some low-complexity-inducing regularizer used in inverse problems ($\ell^1$ norm, group lasso,  nuclear norm, total variation, etc).
The next theorem provides a setting for those functions; in exchange for the strict convexity of $g$, we will require $h$ to also be tilt-conditioned, and to some strong qualification condition to hold.}
{\begin{theorem}[Sum rule for tilt-conditioned functions]\label{T:sum rule}
Let $f=g+h\circ A$, where $g \in \Gamma_0(X)$, $h\in \Gamma_0(Y)$ and $A:X \rightarrow Y$ is a bounded linear operator with closed range.
Suppose that ${\rm{argmin}}~f \neq \emptyset$, and let $\Omega\subset X$.
If $\psi \in \Gamma_0(Y)$ denotes the corresponding Fenchel-Rockafellar dual problem $\psi(v)=g^*(A^*v) + h^*(-v)$, and
\begin{enumerate}[label=\alph*)]
    \item\label{T:sum rule:compatibility} the nondegeneracy condition $0 \in \srint \left(  \dom h - A(\dom g) \right)$ holds,
\end{enumerate}
then ${\rm{argmin}}~\psi \neq \emptyset$. Moreover, if
\begin{enumerate}[label=\alph*),resume]
    \item\label{T:sum rule:qualification condition} there is $\bar v \in {\rm{argmin}}~\psi$ for which the following qualification conditions are satisfied:
\begin{eqnarray}
	0 & \in & \srint (\partial g^*(A^* \bar v) - A^{-1} \partial h^* ( - \bar v) ), \label{D:QC2} \\
 0 & \in & \srint \left(R(A) - \partial h^*(-\bar v) \right), \label{D:QC1}
\end{eqnarray}    
    \item\label{T:sum rule:conditioning} $g$ is $p_1$-tilt-conditioned on $\Omega$, and $h$ is $p_2$-tilt-conditioned on $A\Omega$ for some $p_1,p_2 \geq 1$,
\end{enumerate}
then $f$ is $p$-conditioned on every bounded subset of $\Omega$, with $p:=\max\{p_1,p_2\}$.
\end{theorem}}

\begin{proof}
{The beginning of this proof starts as in the proof of Theorem \ref{T:sum rule polyhedral}: we use the nondegeneracy assumption \ref{T:sum rule:compatibility} with  \cite[Thm. 16.47]{BauCom} to get some $\bar x \in {\rm{argmin}}~f$ and $\bar v \in - \partial h(A\bar x)$ such that $\bar x \in \partial g^*(A^* \bar v)$.
So the condition \cite[Thm. 19.1.iii]{BauCom} is verified, meaning that strong duality holds (in the sense that $\inf f = -\inf \psi$).
This allows to use \cite[Cor. 19.2]{BauCom} to obtain
\begin{equation}\label{src1}
\bar x \in \argmin f  =  \partial g^*(A^*\bar v) \cap A^{-1}\partial h^*(- \bar v). 
\end{equation}}
We can use again \cite[Cor. 19.2]{BauCom}, this time on the dual problem, to also obtain
\begin{equation*}\label{src1.5}
\bar v \in {\rm{argmin}}~\psi  = -\partial h(A\bar x) \cap {A^*}^{-1} \partial g(\bar x).
\end{equation*}
The above equality allows us to assume, without loss of generality,  that $\bar v$ is the element of ${\rm{argmin}}~\psi$ satisfying \ref{T:sum rule:qualification condition}.
So, it remains to prove that, for all $\delta >0$, there exists $\gamma>0$ such that:
\begin{equation}\label{src2}
(\forall x \in \Omega \cap \delta \B_X \cap \dom f)\quad f(x) - \inf f \geq \gamma \dist^p(x,\partial g^*(A^*\bar v) \cap A^{-1}\partial h^*(- \bar v)).
\end{equation}

Fix $\delta >0$, let $x \in \Omega_\delta:=\Omega \cap \delta \B_X \cap \dom f$, and set $\tilde g=g - \langle A^*\bar v, \cdot \rangle$ and $\tilde h = h + \langle \bar v, \cdot \rangle$.
Setting $p = \max\{p_1,p_2\}$, we see from assumption \ref{T:sum rule:conditioning} and Proposition \ref{P:hierarchy of conditionings with p} that $\tilde g$ and $\tilde h$ are  $p$-conditioned on the bounded sets $\Omega_\delta$ and $A\Omega_\delta$, respectively.
Using the same arguments as in \eqref{srcsc1}, we obtain that  $\argmin \tilde g = \partial g^*(A^*\bar v) \ni \bar x$ and  $\argmin \tilde h= \partial h^*(-\bar{v}) \ni A \bar x$.
Therefore, the conditioning of $\tilde g$ (resp. $\tilde h$) evaluated at $x\in \dom f \subset \dom g = \dom \tilde g$ (resp. $Ax\in A \dom f \subset \dom h = \dom \tilde h$) writes as
\begin{eqnarray*}
g(x) & \geq & g(\bar x) + \langle A^* \bar v, x - \bar x \rangle + (\gamma_{\tilde g, \Omega_\delta}/p) \dist^p(x,\partial g^* (A^* \bar v)), \\
h(Ax) & \geq & h(A\bar x) + \langle - \bar v, Ax - A\bar x  \rangle + (\gamma_{\tilde h,A \Omega_\delta} /p) \dist^p(Ax, \partial h^*(-\bar v)).
\end{eqnarray*}
Summing these two last inequalities gives, 
\begin{equation}\label{src3}
 f(x) - \inf f \geq C_1  \left( \dist^p(x,\partial g^* (A^* \bar v)) + \dist^p(Ax, \partial h^*(-\bar v)) \right),
\end{equation}
with $C_1=p^{-1}\min\{ \gamma_{\tilde g, \Omega_\delta}, \gamma_{\tilde h, A\Omega_\delta}\}$.
Since $\Vert \cdot \Vert_\infty \leq \Vert \cdot \Vert_p$ on $\R^2$, we deduce that
\begin{equation*}\label{src4}
 f(x) - \inf f \geq C_1  \max\left\{ \dist(x,\partial g^* (A^* \bar v)) , \dist(Ax, \partial h^*(-\bar v)) \right\}^p,
\end{equation*}
It remains to lower bound the right hand side by the distance to $\argmin f$.
By Example \ref{Ex:distance to intersection}, thanks to the qualification condition \eqref{D:QC2} and the fact that $\Omega_\delta$ is bounded, we derive from \eqref{src1} 
that there exists $C_2 > 0$ independent of $x$ such that
\begin{equation}\label{src5}
 \dist(x,\argmin f) \leq C_2 \max\{ \dist(x, \partial g^*(A^* \bar v)),\dist (x,A^{-1} \partial h^*(- \bar v)) \}.
\end{equation}
Define $y:=\proj(Ax,R(A)\cap \partial h^*(- \bar v)))$, which is well defined since we assumed $R(A)$ to be closed.
Let $\phi_y \in \Gamma_0(X)$ be defined by $\phi_y(u):=(1/2)\Vert Au - y \Vert^2$.
Since $y \in R(A)$, necessarily $\inf \phi_y=0$, so we deduce from Example \ref{Ex: least squares} that
\begin{equation}\label{src6}
(\forall u \in X) \quad \phi_y(u) \geq (\sigma_{\inf}(A^*A)/2) \dist^2(u, \argmin \phi_y).
\end{equation}
On the one hand, we have $\argmin \phi_y = A^{-1}y \subset A^{-1} \partial h^*(- \bar v)$.
On the other hand, the definition of $y$ implies $\phi_y(x)=(1/2)\dist^2(Ax,R(A) \cap  \partial h^*(- \bar v))$.
Thus, it follows from \eqref{src6}  that
\begin{equation*}\label{src7}
 \dist(Ax,R(A) \cap  \partial h^*(- \bar v)) \geq \sigma_{\inf}(A) \ \dist(x, A^{-1} \partial h^*(- \bar v))).
\end{equation*}
Since this is true for any  $x \in \Omega_\delta$, we can combine it with \eqref{src5} to get for all $x \in \Omega_\delta$
\begin{equation}\label{src8}
 \dist(x,\argmin f) \leq C_3 \max\{ \dist(x, \partial g^*(A^* \bar v)), \dist(Ax,R(A) \cap  \partial h^*(- \bar v))  \},
\end{equation}
with $C_3=C_2 \max\{1, \sigma_{\inf}(A)^{-1}\}$.
To end the proof, use the qualification condition \eqref{D:QC1} with Example \ref{Ex:distance to intersection} again to get some $C_4 > 0$ such that for all $x \in  \Omega_\delta$, 
\begin{eqnarray}
 \dist(Ax,R(A) \cap  \partial h^*(- \bar v)) & \leq & C_4 \max \{ \dist(Ax, R(A)), \dist(Ax, \partial h^*(-\bar v)) \} \label{src9} \\
 & = & C_4 \dist(Ax, \partial h^*(-\bar v)). \notag
\end{eqnarray}
The above inequality, combined with \eqref{src8} and \eqref{src2}, concludes the proof.
\end{proof}
{\begin{remark}[On the qualification conditions]\label{R:sum rule qualification condition}
When $g$ is not strictly convex, the conclusion of Theorem \ref{T:sum rule} may not hold if the qualification conditions \eqref{D:QC2} and \eqref{D:QC1} are removed, as proved in \cite[Section 4.4.4]{ZhoSo15} with $g = \Vert \cdot \Vert_*$.
Let us give some sufficient conditions for \eqref{D:QC2} and \eqref{D:QC1} to hold:
\begin{itemize}
    \item If $X$ and $Y$ have finite dimension,  \ref{T:sum rule:qualification condition} is equivalent to 
    \begin{equation*}
        0 \in \rint A \partial g^*(A^*\bar v) - \rint \partial h^*(-\bar v).
    \end{equation*}
    To prove this, use \cite[Cor. 6.15]{BauCom} and \cite[Thm. 6.7]{Roc} to see that the above condition is equivalent to \eqref{D:QC2}, which implies \eqref{D:QC1}.
    This condition is for instance satisfied if $0 \in \rint \partial \psi(\bar v)$ and $0 \in \rint \dom g^* + A^* (\rint \dom h^*)$ (see \cite[Thm. 16.47]{BauCom}).
    Those are the two conditions needed in \cite[Theorem 4.2]{DruLew16}.
    \item If $X$ and $Y$ have finite dimension and $h$ is strictly convex, then a sufficient condition for \ref{T:sum rule:qualification condition} is $\bar x \in \rint \partial g^*(A^* \bar v)$ \cite[Prop. 18.9]{BauCom}.
    \item If $X$ and $Y$ have finite dimension, $g$ is polyhedral and $h$ is strictly convex, then assumption \ref{T:sum rule:qualification condition} is not needed. 
    As pointed out in \cite[Cor. 4.3]{DruLew16}, this is due to the fact that the subdifferentials of $h^*$ and $g^*$ are polyhedral, which allows the use of Hoffman's bound \cite{Hof52} instead of \cite[Theorem 4.3]{BauBor93} in the proof.
\end{itemize}
\end{remark}}
{\begin{remark}[On the closedness of the range]\label{R:sum rul closed range}
In Theorem \ref{T:sum rule polyhedral} we assume $R(A)$ to be closed.
To see how important this hypothesis is in infinite dimension, take $g =0$ (which is not strictly convex), $h=\Vert \cdot \Vert^2$ and $A$ an  operator with a nonclosed range.
Then, for this example, the qualification conditions cannot be satisfied.
Indeed, even if \eqref{D:QC2} is automatically satisfied (because $\partial g^*(0) = X$), condition \eqref{D:QC1} reduces to $0 \in \srint R(A)$, which is equivalent by definition to $R(A) = \cl R(A)$, which is impossible.
Worse, even if we could get rid of this qualification condition, and if the conclusion of the theorem were true, we would obtain that $x \mapsto \Vert Ax \Vert^2$ is $2$-conditioned on bounded sets, which was proven to be impossible in \cite[Theorem 2.1]{HarJen11} (combine it with Proposition \ref{P:equivalence geometrical notions}).
\end{remark}}

{\begin{remark}[Previous results]\label{R:sum rule previous results}
Our results can be seen as extensions and refinements of arguments from \cite{DruLew16}, where the authors introduce the ideas of exploiting the $2$-conditioning of tilted functions on compact sets, together with the description of ${\rm{argmin}}~f$ as an intersection \eqref{src1}. 
Theorem \ref{T:sum rule} improves on \cite[Thm. 4.2]{DruLew16} and \cite[Cor. 4.3]{DruLew16} which require the ${\rm{argmin}}~f$ to be bounded, and $h$ to be in $C^1$ with $\dom h=Y$ (we only ask for a compatibility condition which is satisfied if $h$ is continuous at $A \bar x$, see Remark \ref{R:compatibility condition}).
As far as we know, Theorem \ref{T:sum rule polyhedral} is the first sum rule of this kind with such weak assumptions on $g$.
\end{remark}
To illustrate the interest of these sum rules, we provide a new result for regularized inverse problems where the loss function is the Kullback-Leibler divergence, and the regularizer is a polyhedral function, such as the $\ell^1$ norm, or the Total Variation, which are commonly used in the signal and image processing literature.
\begin{proposition}
Let $f(x) = g(x) + KL(y;Ax)$, where $g \in \Gamma_0(\mathbb{R}^N)$ is polyhedral, $A \in \mathcal{M}_{M,N}(\mathbb{R})$, and $y \in ]0,+\infty[^M$.
If ${\rm{argmin}}~f \neq \emptyset$, then $f$ is $2$-conditioned on bounded sets.
\end{proposition}
\begin{proof}
We just have to verify the hypotheses of Theorem \ref{T:sum rule}, by noting $h:=KL(y; \ \cdot \ )$.
First, the nondegeneracy condition \ref{T:sum rule:compatibility} is verified because $\dom h$ is open, and $h$ is continuous on its domain (see Remark \ref{R:compatibility condition}).
Second, the qualification conditions \ref{T:sum rule:qualification condition} are not needed because we are in a finite dimensional setting, $g$ is polyhedral and $h$ is strictly convex (see Remark \ref{R:sum rule qualification condition}).
Finally, $g$ being polyhedral implies that it is globally $1$-tilt-conditioned (see \cite[Corollary 3.6]{BurFer93}), and we prove in Lemma \ref{L:kullback libler tilt conditioning} that $h$ is $2$-tilt-conditioned on bounded sets, so \ref{T:sum rule:conditioning} is verified.
\end{proof}}
\section{Sharp convergence rates for the Forward-Backward algorithm}
\label{S:CV rates for Forward Backward}
In this section, we present sharp convergence results for the forward-backward algorithm applied to 
$p$-{\L}ojasiewicz functions on a subset $\Omega$, building on the ideas in \cite{AttBolSva13}.
We extend the analysis to the case where $\Omega$ is an arbitrary set, which will allow us to deal
with infinite dimensional inverse problems (see Section \ref{SS:least squares in Hilbert spaces}), 
or structured problems for which all the information is encoded in a manifold (see Section \ref{SS:sparse inverse problems}).
We also provide explicit rates of convergence, for both the iterates and the values.
The proofs of Section \ref{SS:Convergence with Lojasiewicz} are left in the Annex \ref{SS:Annex proof section 4}.

\subsection{Refined analysis with $p$-{\L}ojasiewicz functions}
\label{SS:Convergence with Lojasiewicz}

\begin{theorem}[Strong convergence and rates, $p \geq 1$]\label{T:CVKL discrete}
Suppose that Assumption~\ref{ass:H} is in force, and that $f$ is bounded from below.
Let  $(x_n)_\nin$ be generated by the FB algorithm.
Assume that:
\begin{enumerate}[a)]
	\item \emph{(Localization)} for all $\nin$, $x_n \in \Omega \subset X$,
	\item \emph{(Geometry)} $f$ is $p$-\Loja on $\Omega$, for some $p \geq 1$.
\end{enumerate}
Then the sequence $(x_n)_\nin$ has finite length in $X$, meaning that $\sum_{\nin} \Vert x_{n+1} - x_n \Vert < +\infty$,
and converges strongly to some $x_\infty \in \argmin f \neq \emptyset$. Moreover, there exist some constants $C_p,C_p'>0$ 
with explicit expressions (see equations~\eqref{e:cp} and \eqref{e:cpp}), such that the following convergence rates 
hold, depending on the value of $p$, and of $\kappa:= \lambda(2-\lambda L) [2c_{f,\Omega}^2]^{-1}$:
\begin{enumerate}[i)]
\item If $p=1$, then  $x_n=x_\infty$  for every $n\geq (f(x_0)-\inf f)/\kappa$.
\item  If $p \in ]1,2[$, the convergence is superlinear: for all $\nin$,
	\begin{equation*}
	f(x_{n+1}) - \inf f \leq \left( \frac{f(x_n) - \inf f}{\kappa} \right)^{\frac{p}{2(p-1)}} \quad\text{ and }\quad \Vert x_{n+1} - x_\infty \Vert \leq  C_p (f(x_n)-\inf f)^{1/2},
	\end{equation*}
\item If $p=2$, the convergence is linear: for all $\nin$,
	\begin{equation*}
	\hspace*{-0.5cm} f(x_{n+1}) - \inf f \leq \frac{1}{1+\kappa} (f(x_n) - \inf f) \text{ and } \Vert x_{n+1} - x_\infty \Vert \leq C_2 (f(x_0) - \inf f)^{1/2} \left( {1+\kappa} \right)^{-n/2}.
	\end{equation*}
\item If $p\in ]2,+\infty[$, the convergence is sublinear: for all $\nin$,
	\begin{equation*}
	f(x_n) - \inf f \leq( C_p')^{p/(p-2)}{n^{-\frac{p}{p-2}}} \quad \text{ and }\quad  \Vert x_{n+1} - x_\infty \Vert \leq  C_p(C_p')^{1/(p-2)}{n^{-\frac{1}{p-2}}}.
	\end{equation*}
\end{enumerate}
\end{theorem}
Note that the rates range from the finite termination, for $p=1$, to the worst-case rates seen in Theorem \ref{T:CV FB}, when $p$ tends to $+\infty$.
The bigger is $p$, the more the function is ill-conditioned, in the sense that the rates of its values become closer to $o(n^{-1})$, and the rates of its iterates become arbitrarily slow.
\begin{remark}[Related work]
Theorem~\ref{T:CVKL discrete} collects known and new results. 
We present a simple proof of this theorem, focusing on the analysis of a real sequence satisfying \eqref{cfl7} (see \cite[Theorem 3.2]{ChoPesRep14} or 
\cite[Theorem 3.4]{FraGarPey15} for previous results). 
The superlinear rates in ii), which were known for the proximal point algorithm \cite{Luq84}, are new  for the Forward-Backward algorithm.
Moreover, the case $p=2$ was giving R-linear rates for the values in \cite{ChoPesRep14,FraGarPey15}, while we prove here Q-linear rates.
Also, the quantification of the number of steps in the case $p=1$ involving $\kappa$ is new.
\end{remark}

\begin{remark}[On the sharpness of the rates I]
\label{R:optimality of rates Lojasiewicz}
Let $f=\Vert \cdot \Vert^p$.   
According to \eqref{E:rates for norm to the p} and \eqref{E:rates for the norm to the p values}, the order of the sublinear rates for the forward-backward algorithm that we obtain for both iterates and values are sharp when $p \in ]2, + \infty[$, see Remark \ref{R:optimality of rates worst case}.
When $p=2$, we see that the proximal algorithm verifies $x_{n+1}=(1+2\lambda)^{-1} x_n$, and the algorithm converges linearly.
Finally, when $p \in \left]1,2 \right[$, the order of superlinearity that we obtain is not sharp, since  
for this function  the proximal algorithm has a Q-superlinear rate of order $(p-1)^{-1}$.
It is shown in \cite[Theorem 3.1]{Luq84} that $\dist(x_n,\argmin f)$ converges with this order for the proximal algorithm.
For this, the author uses the stronger notion of metric subregularity, and we will extend this result in Theorem \ref{T:superlinear rates with conditioning} to the FB algorithm. 
\end{remark}

\begin{remark}[Best stepsize and condition number]
When $p \in \left[1,2\right]$, we directly see that the bigger is $\kappa$, the better are the constants in the rates for the values.
This is  true also for $p>2$, by looking in the proof of Theorem \ref{T:CVKL discrete} to the definition of the constant $C_p'$.
The constant $\kappa$ is maximal when we take $\lambda=L^{-1}$, in which case $\kappa=(L2c_{f,\Omega}^2)^{-1}$.
When $f$ is a $\gamma$-strongly convex function,  $\kappa=\gamma/L$ is  the condition number of $f$ (see Example \ref{Ex:uniformly convex functions}) .
So $(L2c_{f,\Omega}^2)^{-1}$ can be seen as a generalized condition number, 
extending this notion from strongly convex functions to $p$-{\L}ojasiewicz ones.
\end{remark}

In Theorem \ref{T:CVKL discrete} the $p$-{\L}ojasiewicz assumption with $p\in \left[1,+\infty\right[$ implies that the $\argmin f$ is nonempty.
In what follows we will derive convergence rates for the objective function values, even in the case where  
$f$ is bounded from below but has no minimizers. Such results are of interest for instance in function approximation theory, where the goal is to find the best 
approximation of a  target function within a specified function class \cite{Dev86}. 
Since in general the considered classes are not closed in the ambient space, the minimizer of the error does not exist, but convergence rates in objective function values are
useful. 
A similar problem appears also in supervised statistical learning theory, where some convergence results can still be obtained are available 
(see e.g. \cite[Theorem 9]{DevRosVer06} and \cite[Theorem A.1]{DevCapRos05}).

We show below that the $p$-{\L}ojasiewicz notion can be extended to \textit{nonpositive} values of $p$, which allows to describe thegeometry of problems without minimizers. 
Based on this new definition, we then derive sharp convergence rates for the objective function values. 

\begin{definition}
Let $p \in \left]-\infty,0\right[$, let $f \in \Gamma_0(X)$ be bounded from below, and let $\Omega \subset X$.
We say that $f$ is \textit{$p$-{\L}ojasiewicz} on $\Omega$ if $ \exists c_{f,\Omega}>0$ such that the {\L}ojasiewicz inequality holds:
\begin{equation*}
\forall x  \in \Omega \cap \doms,   \quad (f(x) - \inf f)^{1- \frac{1}{p}} \leq c_{f,\Omega} \Vert \partial f(x) \Vert_\_. 
\end{equation*}
\end{definition}

\noindent Similarly to the case $p \geq 1$, where this property describes the behavior of $f$ around its minimizers, here it describes the decay of $f(x)$ 
when $\Vert x \Vert$ goes to  $+\infty$.
This assumption leads to convergence rates, interpolating between $o(1)$ and $o(n^{-1})$, depending on the value of $p<0$.
We will see in Section \ref{SS:least squares in Hilbert spaces} that this result applies to ill-posed linear problems involving a compact operator between infinite dimensional spaces.

\begin{theorem}[Rates of convergence, $p<0$]\label{T:CVKL discrete rates negative p}
Let $f\in \Gamma_0(X)$ be bounded from below and satisfying Assumption \ref{ass:H},  $(x_n)_\nin$ be generated by the FB algorithm.
Assume that:
\begin{enumerate}[a)]
	\item \emph{(Localization)} for all $\nin$, $x_n \in \Omega \subset X$,
	\item \emph{(Geometry)} $f$ is $p$-\Loja on $\Omega$, for some $p <0$.
\end{enumerate}
Then the values converge sublinearly (with $C_p'$ defined as in \eqref{e:cp}):
\begin{equation*}
\forall n \in \N, \ f(x_n) - \inf f \leq C_p'^{\frac{p}{p-2}}{n^{\frac{p}{2-p}}}.
\end{equation*}
\end{theorem}

\begin{remark}[On the sharpness of the rates II]\label{R:optimality of rates Lojasiewicz negatif}
The rates obtained in Theorem \ref{T:CVKL discrete rates negative p} are sharp.
Indeed, the function defined in \eqref{E:counter example sublinear rates} is $p$-{\L}ojasiewicz on $\R$ with $p=-\alpha$, and our  rates match the lower bounds obtained in Remark \ref{R:optimality of rates worst case}. 
\end{remark}

Theorem \ref{T:CVKL discrete rates negative p}, together with Theorem \ref{T:CVKL discrete}, give a complete (and sharp) picture of the asymptotic behavior of the FB algorithm.
In fact, looking at the proofs of the mentioned results, we see that the  only properties of forward-backward algorithm that are used are \eqref{cflH1} and \eqref{cflH2}.
We can then extend the previous theorems to a broader class of first-order descent methods, which encompasses block coordinate descent methods, and/or variable metric extensions of the FB algorithm \cite{AttBolSva13,BolSabTeb13,FraGarPey15}.

\begin{theorem}[General first-order descent method]\label{T:CV general descent}
The statements of Theorems \ref{T:CVKL discrete} and  \ref{T:CVKL discrete rates negative p} remain true if the sequence $(x_n)_\nin$ is generated by any algorithm satisfying:
\begin{eqnarray}
(\exists a > 0) & \quad a \Vert x_{n+1} - x_n \Vert^2 \leq f(x_{n+1}) - f(x_n)  \label{H1}\\
(\exists b > 0) & \quad \Vert \partial f(x_{n+1}) \Vert_\_ \leq b \Vert x_{n+1} - x_n \Vert. \label{H2}
\end{eqnarray}
In that case the constant appearing in Theorem \ref{T:CVKL discrete} becomes $\kappa:= ab^{-2}c_{f,\Omega}^{-2}$.
\end{theorem}
\subsection{How to localize the sequence of iterates}\label{SS:localization}
One of the two assumptions we do in Theorems \ref{T:CVKL discrete} and \ref{T:CVKL discrete rates negative p} is
that the sequence  belongs to a set $\Omega$ on which the geometry of $f$ is known.
We discuss here some possible choices.
One first simple case is when $\Omega$ remains invariant under the action of $\T$ (see also Annex \ref{SS:equivalence geometric notions}).

\begin{definition}\label{D:invariant sets FB}
We  say that $\Omega \subset X$ is \textit{FB-invariant} if for all $\lambda \in ]0,2L^{-1}[$, 
$ \T \Omega \subset \Omega$.
\end{definition}

\begin{example}[FB-invariant sets]
Theorem~\ref{T:CV FB}\ref{T:CV FB:i}-\ref{T:CV FB:ii} and Lemma \ref{L:estimates for the forward-backward}.\ref{L:estimates for the forward-backward:gradients} imply that  these sets are FB-invariant (as well as any of their intersection):
\begin{itemize}
	\item $\B_X(\bar x,\delta)$ and $\overline{\B}_X(\bar x,\delta)$ for every $\bar x \in \argmin f$, and for every $\delta \in \left]0, +\infty\right]$,
	\item $[f<r]$ for every $r>\inf f$,
	\item $\{x \in X \ | \ \Vert \partial f(x) \Vert_\_ < M\}$ and $\{x \in X \ | \ \Vert \partial f(x) \Vert_\_ \leq M\}$, for every $M \in ]0,+\infty]$,
	\item  $\Omega=\{x_n\}_\nin$ if $(x_n)_\nin$ is generated by the FB algorithm. 
\end{itemize}
\end{example}

Assuming that $\Omega$ is FB-invariant,  the localization property becomes a simple assumption on the initialization of the algorithm.
The proof of the next corollary is immediate:

\begin{corollary}[Geometry on stable sets gives global rates]\label{T:CV on invariant sets}
Let $f\in \Gamma_0(X)$ be bounded from below and satisfying Assumption~\ref{ass:H}, and  $(x_n)_\nin$ be generated by the FB algorithm.
Assume that $\Omega \subset X$ is FB-invariant and that:
\begin{enumerate}[a)]
	\item \emph{(Initialization)}  $x_0 \in \Omega$,
	\item \emph{(Geometry)} $f$ is $p$-\Loja on $\Omega$, for some $p \in ]-\infty,0[\cup[1,+\infty[$.
\end{enumerate}
Then the results of Theorems \ref{T:CVKL discrete} and \ref{T:CVKL discrete rates negative p} apply for the sequence $(x_{n})_\nin$.
\end{corollary}

In some cases, it is possible to remove the assumption $x_0 \in \Omega$, to the price of having only \textit{asymptotic rates}.
Indeed, it suffices to prove that the sequence will enter in $\Omega$ at a certain iteration, which is the argument used in \cite{AttBolSva13,FraGarPey15}, in a non-convex setting.
This happens for instance with the local level sets, under a slight compactness assumption (see below).

\begin{corollary}[Local geometry gives asymptotical rates]\label{P:capture result for local omega}
Let $f\in \Gamma_0(X)$ be such that $\argmin f \neq \emptyset$ and satisfying Assumption~\ref{ass:H}.
Let $(x_n)_\nin$ be generated by the FB algorithm and assume that:
\begin{enumerate}[a)]
	\item \emph{(Compactness)} $(x_n)_\nin$ admits a  subsequence strongly converging to $\bar x$ in $X$,
	\item \emph{(Local geometry)} for some $p\in\left[1,+\infty\right[$:
	\begin{equation*}\label{E:local Loja}
 (\exists (\delta,r)\in \left]0,+\infty\right]) \text{ such that $f$ is $p$-{\L}ojasiewicz on } \B_X(\bar x,\delta) \cap [f<r+\inf f].
\end{equation*}
\end{enumerate}
Then there exists $n_0 \in \N$ such that the rates of Theorem \ref{T:CVKL discrete} apply for the sequence $(x_{n_0+n})_\nin$.
\end{corollary}

\begin{proof}
Let $(x_{n_k})_{k\in\mathbb{N}}$ be a subsequence strongly converging to some $x_\infty $, which belongs to $\argmin f$ according to Theorem \ref{T:CV FB}.
Therefore, $f$ is $p$-{\L}ojasiewicz on $\Omega:=\B_X(x_\infty,\delta) \cap [f<r+\inf f]$, for some $(\delta,r)\in \left]0,+\infty\right]$.
Since $x_{n_k}\to x_\infty$ and $f(x_{n_k}) \downarrow \inf f$, there exists  $K\in \N$ such that $x_{n_K} \in \Omega$.
Since $\Omega$ is FB-invariant, we conclude that $(x_n)_{n \geq N} \subset \Omega$.
\end{proof}

\begin{remark}[On the compactness assumption]
The compactness assumption made in Corollary \ref{P:capture result for local omega} is always satisfied in finite dimension.
Indeed Theorem \ref{T:CV FB} guarantees that the sequence is bounded  under the assumption that $\argmin f \neq \emptyset$.
If $X$ has infinite dimension, this assumption can be verified provided that $f$ has compact level sets, due to the decreasing property of $f(x_n)$.
\end{remark}

The property that a sequence $(x_n)_{\nin}$ generated by an algorithm reaches a set  of interest $\Omega$ after a finite number of iterations, is usually called identifiability, or finite identification of $\Omega$  \cite{Wri93,Lew02,HarLew04}, and $\Omega$ is therefore called an \textit{active set}.
For instance, the so-called \textit{active manifolds} can be identified in finite time, under the assumption that $f$ is {partially smooth} with respect to this manifold  \cite{HarLew04,HarLew07}.
{An alternative approach, recently introduced in \cite{FadMalPey18}, shows that the \textit{strata} of mirror-stratifiable functions are identifiable.
We will use this notion of active strata to derive another asymptotic convergence result.

Before introducing the notion of mirror-stratifiability, we recall that
a set $M \subset \mathbb{R}^N$ is said to be stratified by $\{M_i\}_{i=1}^s \subset M$ if this family is a finite partition $\sqcup M_i = M$ such that $M_i \cap \cl M_j \neq \varnothing \Leftrightarrow M_i \subset \cl M_j$.
The latter inclusion endows the family of strata with an order relation $M_i \preceq M_j \Leftrightarrow M_i \subset \cl M_j$.
Given a point $x \in M$, it will be useful to note $M_x$ the unique strata which contains $x$.
{
\begin{definition}[Mirror-stratifiable function] We say that a function $f \in \Gamma_0(\mathbb{R}^N)$ is \textit{mirror-stratifiable} if
\begin{enumerate}[label=\alph*)]
    \item $\dom \partial f$ (resp. $\dom \partial f^*$) is stratified by $\{M_i\}_{i=1}^s$ (resp. $\{M_i^*\}_{i=1}^s$),
    \item the map $J_f : M \longmapsto \bigcup\limits_{x \in M} \rint \partial f(x)$ realizes a bijection between $\{M_i\}_{i=1}^s$ and $\{M_i^*\}_{i=1}^s$,
    \item the map $J_f$ is decreasing, in the sense that $M_i \preceq M_j \Leftrightarrow J_f(M_j) \preceq J_f(M_i)$.
\end{enumerate}
\end{definition}
Both notions appear naturally in most sparsity-based inverse problems such as the $1$-norm, group-lasso norm, nuclear norm, or the total variation, or any polyhedral function, see  \cite{FadMalPey18} for more details and many examples.

\begin{corollary}\label{T:partial smoothness rates}
Suppose that Assumption~\ref{ass:H} is in force, that $X=\R^N$, and let $(x_n)_\nin$ be the sequence generated 
by the FB algorithm converging to some $\bar x \in \argmin f$. Assume that:
\begin{enumerate}[a)]
    \item \label{T:mirror stratification rates:mirror} $g$ is mirror-stratifiable, and we define $C_{\bar x} := \cup \{ M \ | \ M_{\bar x} \preceq M \preceq J_g^{-1} (M^*_{-\nabla h(\bar x)}) \}$,  
	\item\label{T:partial smoothness rates:loja} $f$ is $p$-{\L}ojasiewicz on $C_{\bar x} \cap \B_X(\bar x,\delta)$ for some $\delta \in ]0, +\infty]$ and $p\in[1,+\infty[$.
\end{enumerate}
Then there exists $n_0 \in \N$ such that the rates of Theorem   \ref{T:CVKL discrete}  apply for the sequence $(x_{n_0+n})_\nin$. 
Note that $C_{\bar x} = M_{\bar x}$ holds whenever $0 \in \rint \partial f(\bar x)$.
\end{corollary}
\begin{proof}

It follows from \cite[Theorem 4]{FadMalPey18} that there exists $n_0 \in \N$ for which $x_{n_0+n}\in C_{\bar x}$ for every $\nin$.
Since $(x_n)_{\nin}$ converges to $\bar x$, we can assume that $n_0$ is such that $x_{n_0+n}\in C_{\bar x} \cap \B_X(\bar x,\delta)$ for every $\nin$.
This, together with \ref{T:partial smoothness rates:loja}, allows to apply Theorem \ref{T:CVKL discrete} to the sequence $(x_{n_0+n})_\nin$.
The equality $C_{\bar x} = M_{\bar x}$ follows directly from the bijectivity of $J_g$, and the fact that $-\nabla h(\bar x) \in \rint \partial g(\bar x)$.
\end{proof}
The reader not familiar with the notion of mirror-stratifiability might wonder what is the active set $C_{\bar x}$ appearing in Corollary \ref{T:partial smoothness rates}.
Here are a few example of interest:
\begin{example}\label{Ex:mirror active sets}
We keep here the notations of Corollary \ref{T:partial smoothness rates}:
\begin{itemize}
    \item If $g(x) = \Vert x \Vert_1$, we can choose a stratification based on sets with prescribed support, which gives
    \begin{equation}\label{Ex:mirror active sets:L1 support}
        C_{\bar x} = \{x \in \mathbb{R}^N \ | \ \supp(\bar x) \subset \supp(x) \subset \act(-\nabla h(\bar x)) \},
    \end{equation}
    where $\supp(x)$ is the support of $x$, and $\act(x^*) = \{i \ | \ \vert x_i \vert =1 \}$ is the set of active indices of $x^*$ in $[-1,1]^N$. 
    Some authors call $\act(-\nabla h(\bar x))$ the extended support of $\bar x$. In the case that $0 \in \rint \partial f(\bar x)$, we have $\supp(\bar x) = \act(-\nabla h(\bar x))$.
\item If $g(x) = \Vert x \Vert_*$ is the nuclear norm, we can choose a stratification based on sets of matrices with prescribed rank, which gives
    \begin{equation}\label{Ex:mirror active sets:nuclear}
        C_{\bar x} = \{x \in \mathcal{M}_{M,N}(\mathbb{R}) \ | \ \rank(\bar x) \leq \rank(x) \leq \#\act(\sigma(-\nabla h(\bar x))) \},
    \end{equation}
    where $\sigma(x^*)$ denotes the set of singular values of the matrix $x^*$. 
    If $0 \in \rint \partial f(\bar x)$, we have $\rank(\bar x) = \#\act(\sigma(-\nabla h(\bar x)))$.
\end{itemize}
\end{example}
\begin{remark}[Partial smoothness]
Even if there is no direct relation between mirror stratification and partial smoothness, all the above mentioned functions are both mirror-stratifiable and partially smooth, and it would be immediate to provide an analogue 
result to Corollary \ref{T:partial smoothness rates} for partially smooth functions.
Note that when using the identification theorems for partially smooth functions, it is necessary to assume the qualification condition $0 \in \rint \partial f(\bar x)$ to hold.
In this case, the active manifold coincide with the active set $C_{\bar x} = M_{\bar x}$ for most practical cases (polyhedral functions, spectral norms), meaning that those cases are already covered by Corollary \ref{T:partial smoothness rates}.
\end{remark}
\begin{remark}[On the assumptions]
Note that our assumptions do not require or imply that $f$ has unique minimizer; we only require $f$ to be \Loja on the active set.
In Section \ref{SS:sparse inverse problems}, we will show how this geometrical assumption can be guaranteed, provided that $\nabla^2 h(\bar x)$ is injective when restricted to the tangent cone of the active set.
In \cite[Thm. 3.7]{LiPon17} the authors provide a sufficient condition for the \Loja inequality to hold locally when $g$ is a partially smooth function.
\end{remark}
}}

\subsection{Linear rates of convergence for the Forward-Backward algorithm}
\label{SS: linear rates}

In this Section  we give more insights on the linear rates for the FB algorithm.
According to Theorem \ref{T:CVKL discrete},  $f(x_n)-\inf f$ and $\Vert x_n - x_\infty \Vert$  converge linearly when a $2$-{\L}ojasiewicz property is verified.
Another decreasing quantity of interest is $\dist(x_n,\argmin f)$, and its Q-linear convergence is equivalent to asking that the forward-backward map $T_\lambda$ satisfies
\begin{equation}\label{D:Linear rates of distance minimizers}
(\exists \eps_{f,\Omega}\in ]0,1[) (\forall x \in \Omega \cap \dom f)\quad \dist(\T x,\argmin f) \leq \eps_{f,\Omega} \dist(x,\argmin f).
\end{equation}
If such property holds on a set $\Omega$ containing  $(x_n)_\nin$, the sequence $(\dist(x_n,\argmin f))_\nin$ will converge Q-linearly.
In fact, it is possible to show that \eqref{D:Linear rates of distance minimizers} is \textit{equivalent} to the $2$-conditioning of $f$ on 
$\Omega$, provided this set is FB-invariant (see Definition \ref{D:invariant sets FB}).
This fact has been observed in \cite{NecNesGli15} for the projected gradient method, with $\Omega=X$ and $\lambda=L^{-1}$, and 
below we extend the  argument to our more general setting.

\begin{proposition}[Linear rates and $2$-conditioning]\label{P:equivalence linear rates and 2-conditioning}
Suppose that Assumption~\ref{ass:H} is in force and assume that $\argmin f \neq \emptyset$. Let $\Omega \subset X$ and $\lambda \in ]0,2 L^{-1}[$.
\begin{enumerate}[i)]
	\item
\label{P:equivalence linear rates and 2-conditioning:i}
	 \hspace{-0.045cm}If $f$ verifies \eqref{D:Linear rates of distance minimizers} on $\Omega$, then it is $2$-conditioned on $\Omega$ with  $\gamma_{f,\Omega} = \lambda^{-1}(2 - \lambda L) (1-\eps_{f,\Omega})^{2}$.
	\item
	\label{P:equivalence linear rates and 2-conditioning:ii}
	 \hspace{-0.045cm}If $f$ is $2$-conditioned on $\T \Omega$, then it verifies \eqref{D:Linear rates of distance minimizers} on $\Omega$ with $\eps_{f,\Omega}=(1+\lambda \gamma_{f,\Omega})^{-1/2}$ for 
	stepsizes $\lambda\in \left ]0,L^{-1}\right]$.
\end{enumerate}
Then, on FB-invariant sets, the $2$-conditioning is equivalent to \eqref{D:Linear rates of distance minimizers},  for stepsizes $\lambda\in \left ]0,L^{-1}\right]$.
\end{proposition}

\begin{proof}
Let $S=\argmin f$, and let $x \in \Omega$.
It follows from the triangular inequality that
\begin{equation}\label{elt1}
\dist(x,S) \leq \Vert x - \proj(\T x,S) \Vert \leq \Vert x - \T x \Vert + \dist(\T x,S).
\end{equation}
Lemma~\ref{L:estimates for the forward-backward}.\ref{L:estimates for the forward-backward:descent} implies that
\begin{equation}\label{elt2}
\Vert x - \T x \Vert^2 \leq 2\lambda(2-\lambda L)^{-1} (f(x)-\inf f) 
\end{equation}
\noindent For item \ref{P:equivalence linear rates and 2-conditioning:i}, combine  \eqref{D:Linear rates of distance minimizers}, \eqref{elt1},  and \eqref{elt2}:
\begin{equation*}
(1-\eps_{f,\Omega})^2\dist(x;S)^2 \leq \Vert \T x - x \Vert^2 \leq 2\lambda(2-\lambda L)^{-1} (f(x) - \inf f).
\end{equation*}

\noindent For item \ref{P:equivalence linear rates and 2-conditioning:ii}, Lemma \ref{L:estimates for the forward-backward}.\ref{L:estimates for the forward-backward:descent} with $u=\proj(x;S)$, and the fact that $\lambda \leq 1/L$ implies
\begin{equation*}
\Vert \T x - \proj(x;S) \Vert^2 \leq \dist(x;S)^2 - 2 \lambda(f(\T x)-\inf f).
\end{equation*}
Then, since $f$ is $2$-conditioned on $\T \Omega \ni \T x$,  we can conclude from
\begin{equation*}
\dist(\T x;S)^2 \leq \Vert \T x - \proj(x;S) \Vert^2 \leq \dist(x;S)^2 - \lambda \gamma_{f,\Omega} \dist(\T x;S)^2. \tag*{\qedhere}
\end{equation*} 
\end{proof}

Let us assume that $f$ is a $\gamma$-strongly convex function, with $\gamma>0$ as in Example \ref{Ex:uniformly convex functions}, and let $\bar x$ be its unique minimizer.
Let $(x_n)_\nin$ be generated by the FB algorithm, for which we take $\lambda=1/L$, and define the condition number of $f$ as $\kappa:=\gamma/L$.
We compare the different linear rates that we can get for $\Vert x_n - \bar x \Vert$ by using different theorems, relying on more or less strong assumptions.
Using that $f$ is $2$-{\L}ojasiewicz (with $c_{f,X}=(2\gamma)^{-1/2}$, see Example \ref{Ex:uniformly convex functions}), 
 Theorem \ref{T:CVKL discrete} yields R-linear rates of the form
\begin{equation*}
\Vert x_n - \bar x \Vert \leq C\eps^n, \ C>0,
\end{equation*}
where 
$ \eps = {1}/{\sqrt{1 +\kappa}}.$
If instead we exploit  $2$-conditioning (recall that in general this is a stronger notion than 2-\Loja, Proposition \ref{P:equivalence geometrical notions}), we obtain Q-linear 
rates from Proposition \ref{P:equivalence linear rates and 2-conditioning}  with exactly the same constant $\eps$.
If we use directly the strong convexity of $f$, we obtain in this case Q-linear rates with 
$\eps = 1- \kappa$ (see e.g. \cite[Proposition 3]{SchLerBac11}).
So, the more information we use, the better rates we derive.
In \cite{NecNesGli15}, the authors investigate different notions belonging between  strong convexity and the $2$-conditioning.
For instance, under an assumption of ``quasi strong convexity", they obtain
$\eps = \sqrt{{(1-\kappa})/(1+\kappa)},$
which is smaller than $ (1+\kappa)^{-1/2}$, but not as good as $ 1- \kappa$.
In conclusion, two aspects are crucial in the linear convergence of forward-backward.
First, to have Q-linear rates for the iterates, it is necessary and sufficient to require the $2$-conditioning of the function, 
due to the equivalence result of Proposition \ref{P:equivalence linear rates and 2-conditioning}.
Second, just assuming $2$-conditioning is not a guarantee of having a \textit{fast computation} of the solution, since linear rates can be arbitrarily slow on any finite number of iterations.
Indeed two constants play a key role: the condition number $\kappa$, which is directly related to $\gamma_{f,\Omega}$ (some extra assumptions on $f$ could improve the value of $\gamma_{f,\Omega}$,  see e.g. the discussion in Subsection \ref{SS:sparse inverse problems}), and  $\eps$ (see also \cite{NecNesGli15}).

\subsection{Superlinear rates and finite termination}
{In this section, we refine the convergence analysis for the case $p \in \left]1,2\right[$, replacing the $p$-{\L}ojasiewicz property with $p$-metric subregularity (or $p$-conditioning). 
As discussed in  Remark \ref{R:optimality of rates Lojasiewicz}, the order of superlinear convergence that we derive for the FB algorithm in the case $p\in]1,2[$ is not sharp. 
In Theorem \ref{T:superlinear rates with conditioning}, using $p$-metric subregularity (or $p$-conditioning) instead of $p$-\Loja, we derive  better (and indeed sharp, see Remark~ \ref{R:optimality of rates Lojasiewicz}) superlinear rates. 
Keep in mind these three notions are only equivalent via   Proposition~\ref{P:equivalence geometrical notions} if $\Omega$ verifies a stability condition.
The  proof of  Theorem \ref{T:superlinear rates with conditioning} below follows directly from the next lemma, which is a partial analogue of Proposition~\ref{P:equivalence linear rates and 2-conditioning}-\ref{P:equivalence linear rates and 2-conditioning:ii}.}
\begin{lemma}\label{L:error bound for superlinear rates}
Suppose that Assumption~\ref{ass:H} is in force and assume that $\argmin f \neq \emptyset$. 
\begin{enumerate}[i)]
	\item\label{L:error bound for superlinear rates:i}  If $\partial f$ is $p$-metrically subregular on $\Omega \subset X$, then for all $p\in \left]1,2\right[$, and $x \in \doms$:
	\begin{equation*}
	\T x \in \Omega \Rightarrow \dist(\T x,\argmin f)^{p-1} \leq {2}/({\lambda \gamma_{\partial f,\Omega}})^{-1} \dist(x,\argmin f).
	\end{equation*}
	\item\label{L:error bound for superlinear rates:ii} If $f$ is $p$-conditioned on $\Omega$, then   for all $p\in \left]1,2\right[$, and $x \in \doms$:
	\begin{equation*}
	(x,\T x) \in \Omega^2 \Rightarrow (f(\T x) - \inf f)^{p-1} \leq \left( {p}/{\gamma_{f,\Omega}} \right)^2 \left({2}/{\lambda} \right)^p (f(x) - \inf f).
	\end{equation*}
\end{enumerate}
\end{lemma}
\begin{proof}
Let $S=\argmin f$. Lemma \ref{L:estimates for the forward-backward}.\ref{L:estimates for the forward-backward:gradients}, the triangular inequality, and Theorem \ref{T:CV FB}-\ref{T:CV FB:ii} yield
\begin{equation}\label{ebs1}
\lambda \Vert \partial f(\T x) \Vert_\_ \leq  \Vert \T x - x \Vert \leq \Vert \T x - \proj(x,S) \Vert + \Vert \proj(x,S) - x \Vert \leq 2 \dist(x,S).
\end{equation}
For \ref{L:error bound for superlinear rates:i}, use the hypothesis with \eqref{ebs1} to derive $\gamma_{\partial f,\Omega}  \dist(\T x, S)^{p-1} \leq (2/\lambda) \dist(x,S).$
For  \ref{L:error bound for superlinear rates:ii}, use the $p$-{\L}ojasiewicz inequality  via Proposition \ref{P:equivalence geometrical notions} 
, together with \eqref{ebs1} and the $p$-conditioning:
\begin{equation*}
(f(\T x) - \inf f)^{p-1}\leq (p/\gamma_{f,\Omega}) \Vert \partial f(\T x ) \Vert_\_^p   \leq (p/\gamma_{f,\Omega})^2(2/\lambda)^p (f(x)-\inf f) \tag*{\qedhere}
\end{equation*}
\end{proof}

\begin{theorem}\label{T:superlinear rates with conditioning}
Assume that $p \in ]1,2[$ and that the hypotheses of Theorem \ref{T:CVKL discrete} hold.
If the $p$-{\L}ojasiewicz hypothesis  is replaced by $p$-metric subregularity (resp. $p$-conditioning), then $\dist(x_n,\argmin f)$ (resp. $f(x_n) - \inf f$)  Q-superlinearly converges with order $(p-1)^{-1}$.
\end{theorem}

We now discuss the relevance of these fast rates when $f$ is $p$-{\L}ojasiewicz with $p \in \left[1,2\right[$.
While superlinear rates are well-known for the proximal algorithm applied to sharp functions, it is not observed for the gradient method.
The apparent contradiction between this result and practice is in fact related to a quite intuitive fact, stated in the following Proposition:  the more a function is smooth, the less it can be sharp.
This means that the gradient algorithm cannot be applied to $p$-{\L}ojasiewicz function, with $p <2$, because it is incompatible with $\nabla f$ being Lipschitz continuous. 
A similar statement, under  different assumptions, can be found in \cite[Proposition 2.8]{BegBolJen15}.

\begin{proposition}\label{P:conditioning and Holder}
Let $f \in \Gamma_0(X)$ be such that $\dom f$ {has a nonempty interior.} Assume $f$ to be differentiable on $\Omega$, where $\Omega \subset X$ is convex and such that\footnote{Note that $\proj(\Omega;\argmin f) \subset \Omega$ holds when $\Omega=\B_X(\bar x,\delta) \cap [f<r]$, for $\bar x \in {\rm{argmin}}f$, because $\proj(\cdot;\argmin f)$ is nonexpansive.} $\proj(\Omega;\argmin f) \subsetneq \Omega$. 
Assume that $f$ is $p$-conditioned on $\Omega$, and that $\nabla f$ is $\alpha$-H\"older continuous on $\Omega$, i.e.
\begin{equation*}
(\exists L_{\nabla f,\Omega, \alpha}>0)(\exists\alpha>0)(\forall (x,y)\in \Omega^2)\quad \ \Vert \nabla f(x) - \nabla f(y) \Vert \leq L_{\nabla f,\Omega,\alpha} \Vert x - y \Vert^\alpha.
\end{equation*}
Then $p \in [\alpha + 1, +\infty[$. In the case that $p=\alpha+1$, we have moreover that $\gamma_{f,\Omega} \leq L_{\nabla f, \Omega,\alpha}$.
\end{proposition}

\begin{proof}
Let $x \in \Omega \cap \dom^* f$, and $\bar x := \proj(x, \argmin f)$. Then $\bar x \in \Omega$ and $\bar x \neq x$.
For all $t \in \left]0,1\right]$, let $x_t:=t x +(1-t) \bar x$. Then $x_t \in \Omega \setminus\argmin f$ and $\bar x = \proj(x_t, \argmin f)$.
From the $p$-conditioning assumption and the Descent Lemma \ref{L:descent lemma} applied at $(\bar x,x_t) \in \Omega^2$, we see that:
\begin{equation}\label{cah1}
(\forall t \in \left]0,1\right])\quad 0<\frac{\gamma_{f,\Omega}}{p} \Vert x_t - \bar x \Vert^p \leq f(x_t) - f(\bar x) \leq \frac{L_{\nabla f,\Omega}}{\alpha+1} \Vert x_t - \bar x \Vert^{\alpha + 1}.
\end{equation}
If we suppose that $p< \alpha +1$, then by passing to the limit for $t \to 0$, we get $\gamma_{f,\Omega}/p \leq 0$
which is impossible. So $p\geq \alpha +1$, and if equality holds, $\gamma_{f,\Omega} \leq L_{\nabla f, \Omega}$ follows from \eqref{cah1}.
\end{proof}

As a consequence of Proposition \ref{P:conditioning and Holder}, we should not expect more than linear rates for the gradient method applied to a $C^{1,1}$ convex function.
Such a result cannot be extended straightforwardly to the Forward-backward algorithm.
For instance,  the function $f(x)=\Vert x \Vert^2 + \Vert x \Vert$ has a nontrivial smooth term in its decomposition, but is still sharp at its minimizer.

\section{Linear inverse problems: from modeling assumptions to convergence rates}
\label{S:Linear inverse problems}
{Throughout this section, $X$ and $Y$ are Hilbert spaces and 
 $A : X \longrightarrow Y$ is a bounded linear operator. $X$ is called the parameter space and $Y$ is the data space.
Given the linear  inverse problem $Ax=y$, for some $y \in Y$,  we are interested in the (possibly regularized) convex optimization problem
\begin{equation}\label{D:least squares inverse problems}
\min\limits_{x \in X} f(x) = g(x) + D( Ax; y), 
\end{equation}
where $g \in \Gamma_0(X)$ and $h = D(\cdot,y) \in \Gamma_0(Y)$.
The goal of this section is to show that typical modeling assumptions made in the inverse problem literature
can be interpreted as geometric assumptions on \eqref{D:least squares inverse problems}, which are often not local, in the sense of Definition \ref{D:geometric notions}.}
First, we show that the classical source conditions are equivalent to a \Loja  condition on suitable subsets, that we call  source sets. 
Second, we show that the restricted isometry
property, which is the key for exact recovery in sparsity based regularization, induces a 2-conditioning of the problem over a cone of sparse vectors, which is identified in finite time by the algorithm. 
{This result extends to general inverse problems with mirror-stratifiable regularizing functions, for which the restricted isometry property entails a 2-conditioning of the problem over an active set (introduced in Corollary \ref{T:partial smoothness rates}). 
\subsection{\Loja property of quadratic functions via source conditions in Hilbert spaces}\label{SS:least squares in Hilbert spaces}
All across this Section \ref{SS:least squares in Hilbert spaces}, we assume that $A : X \rightarrow Y$ is a bounded linear operator, that $y \in Y$, and that $f(x):= (1/2)\Vert Ax-y \Vert^2$ is the associated least squares function.
We will also note $y^\dagger := \proj(y, \mbox{cl}~R(A))$, and, whenever $\argmin f \neq \emptyset$, we will note $x^\dagger := \proj(0, \argmin f)$, which verifies $Ax^\dagger = y^\dagger$.
\subsubsection{Elements of linear algebra}
Before going further into the topic, let us recall some basic (but not necessarily well-known) facts about bounded linear operators in Hilbert spaces.
A first important difference with the finite-dimensional setting is that the set of minimizers of $f$ can be empty:
\begin{proposition}[{\cite[Theorem 3.1.1]{Gro77}}]\label{P:quadratic empty argmin}
Let $A : X \rightarrow Y$ be a bounded linear operator, $y \in Y$ and $f(x) := \Vert Ax - y \Vert^2/2$.
Then $\argmin f \neq \emptyset \Leftrightarrow y \in R( A) + R(A)^\perp \Leftrightarrow y^\dagger \in R(A)$.
\end{proposition}
\noindent We see that $\argmin f \neq \emptyset$ is guaranteed when $R(A)$ is closed, which for instance cannot happen for compact operators with  infinite-dimensional range \cite[Theorem 3.1.3]{Gro77}.
Observe that the closedness of $R(A)$ can be checked by means of its  singular values:
\begin{proposition}\label{P:closed range singular values}
Let $A : X \rightarrow Y$ be a bounded linear operator. 
Then $R(A) $ is closed if and only if $\sigma_{inf}(A) > 0$.
\end{proposition}
\begin{proof}
Use the fact that $R(A) = R((AA^*)^{1/2})$ \cite[Proposition 2.18]{EngHanNeu} together with \cite[Remark 2.3]{HarJen11} and the fact that $\spec((AA^*)^{1/2}) = \spec(AA^*)^{1/2}$ \cite[\textsection 32 Theorem 3]{Hel69}.
\end{proof}
}

\subsubsection{Known results about the Landweber algorithm}
The quadratic function $f$ can be minimized by means of a gradient method, defined as
\begin{equation}\label{D:Landweber}
(\forall \nin)\quad x_{n+1} = x_n - \lambda A^*(Ax_n - y), \text{ with } x_0 \in X \text{ and }  \lambda\in\left]0,2\| A^*A \|^{-1}\right[.
\end{equation}
A vast literature is devoted to this algorithm, which is often called in this context the Landweber algorithm.
It is well-known that whenever $\argmin f \neq \emptyset$, the sequence $(x_n)_\nin$ generated by the
Landweber algorithm converges strongly to the projection of $x_0$ onto  $\argmin f$ (see e.g. \cite[Theorem 6.1]{EngHanNeu}, or \cite[Theorem 3.3.2]{Gro77} for varying stepsizes).
When the range $R(A)$ is closed, the algorithm behaves exactly as in finite dimensions: both iterates and values converge linearly, see Example~\ref{Ex: least squares} and Theorem~\ref{T:CVKL discrete}.
If the $R(A)$ is not closed, instead, the rates for $\Vert x_n - \bar x_0 \Vert$ can be arbitrarily slow without additional assumptions  \cite[Theorem 12]{DavYin14}. Moreover, \cite[Theorem 2.1]{HarJen11} shows that \textit{no local {\L}ojasiewicz property} can be satisfied by such quadratic function when $R(A)$ is not closed. This could suggest that it is not possible to rely on geometrical assumptions to obtain convergence rates. Nevertheless, as we will see below, this is not true.
Indeed, in the inverse problem literature, this worst-case scenario is avoided by making an extra assumption on the problem.
For instance, if the following \textit{source condition} is verified
\begin{equation}\label{D:source condition}
(\exists \mu \in \left]0,+\infty\right[) \quad  x^\dagger \in R(A^*A)^\mu,
\end{equation}
the Landweber algorithm initialized with $x_0=0$ is known \cite{EngHanNeu} to have the rates
\begin{equation}\label{E:rates for Landweber and source condition}
f(x_n) - \inf f = O (n^{-(1+2\mu) }), \text{ and } \Vert x_n -  x^\dagger \Vert = O(n^{-\mu }).
\end{equation}
Also, when $\argmin f = \emptyset$, a source condition in $Y$ can be made:
 \begin{equation}\label{D:source condition in Y}
(\exists \nu \in \left]0,+\infty\right[) \quad  y^\dagger \in R(AA^*)^\nu,
\end{equation}
so that the Landweber algorithm initialized with $x_0=0$ verifies \cite[Theorem 2.10]{DevRosVer06}:
\begin{equation}\label{E:rates for Landweber and source condition on Y}
f(x_n) - \inf f = O (n^{-2\nu }).
\end{equation}
{The source condition \eqref{D:source condition in Y} can be understood in light of Proposition \ref{P:quadratic empty argmin}.
Indeed, this proposition says that the problem is well posed (in the sense that $\argmin f \neq \emptyset$) when $y^\dagger \in R(A)$.
So it is reasonable to think that the ``deeper'' $y^\dagger$ is in $R(A)$, and the easier the problem is.
In the ill-posed case $y^\dagger \in \cl R(A)\setminus R(A)$, we could also imagine that the ``further away'' $y^\dagger$ is from $R(A)$, and the more difficult the problem is.
Estimating the location of $y^\dagger$ can be done thanks to the spaces $R(AA^*)^\nu$, because they form a sequence of nonincreasing dense subsets of $\cl R(A)$ (see Lemma \ref{L:powers self adjoint operators} and \cite[Proposition 2.8]{EngHanNeu}):
\begin{equation*}
    \mbox{cl}~R(A) = \mbox{cl}~\bigcup\limits_{\nu > 0} R(AA^*)^\nu 
    \quad \text{ and } \quad
    R(AA^*)^{1/2} = R(A).
\end{equation*}}
The aim of this section is to highlight how the rates \eqref{E:rates for Landweber and source condition} and \eqref{E:rates for Landweber and source condition on Y} can be simply explained using the results of Section~\ref{S:CV rates for Forward Backward}.
We show that the source conditions \eqref{D:source condition} and \eqref{D:source condition in Y} are equivalent to assume that the initialization $x_0$ of the algorithm belongs to a so-called \textit{source set}.
Our main result in this section consists  in showing that the function $f$ satisfies a {\L}ojasiewicz inequality on these source sets, which are FB-invariant.
As a by-product of Corollary \ref{T:CV on invariant sets}, we will obtain a new and simple geometrical interpretation of the rates in \eqref{E:rates for Landweber and source condition} and \eqref{E:rates for Landweber and source condition on Y}.

\subsubsection{Regularity spaces and source sets}

{
\begin{definition}[Regularity space and source set]\label{D:regularity spaces and source sets}
\begin{enumerate}
    \item Given $(\nu,\delta) \in \left]0 , + \infty\right[\times \left]0 , + \infty\right[ $, the data regularity space and the data source set are respectively defined as: 
\begin{equation*}
    Y_\nu := y^\dagger + R(AA^*)^\nu
    \quad \text{ and } \quad
    Y_{\nu,\delta}  :=  y^\dagger  + \{ (AA^*)^\nu \omega \ | \  \omega  \in \cl R(A), \Vert \omega \Vert \leq \delta \}.
\end{equation*}
    \item Given $(\mu,\delta) \in \left]-1/2 , + \infty\right[\times \left]0 , + \infty\right[ $, the  regularity space and the source set are respectively defined as: 
\begin{equation*}
X_\mu   :=  A^{-1} Y_{\mu+1/2} \quad \text{ and } \quad X_{\mu,\delta} (y) :=A^{-1} Y_{\mu+1/2,\delta},
\end{equation*}
where $A^{-1}$ denotes the preimage of a set under the application $A$.
\end{enumerate}
\end{definition}

\begin{proposition}\label{P:Regularity spaces in terms of mu}
~~
\begin{enumerate}[i)]
    \item\label{P:Regularity spaces in terms of mu:i}   $\argmin f = \emptyset$ if and only if $X_\mu=\emptyset$ for all $\mu \in \left[0, + \infty\right[$.
	\item\label{P:Regularity spaces in terms of mu:ii}  $\argmin f \neq \emptyset$ if and only if $X_\mu=X$ for all $\mu \in \left]-1/2,0\right]$.
	\item\label{P:Regularity spaces in terms of mu:iii} Assume $R(A)$ is closed. Then $X_\mu=X$ for all $\mu \in \left]-1/2,+\infty\right[$.
\end{enumerate}
\end{proposition}
\begin{proof}
Given any $x \in X$, observe that $x \in X_0$ is, by definition, equivalent to $Ax \in Y_{1/2}$.
Since $R(A)= R({AA^*})^{1/2}$, the latter is equivalent to $Ax \in y^\dagger + R(A)$.
We can then easily deduce, using also Proposition \ref{P:quadratic empty argmin}, that $
    X_0 = X
    \Leftrightarrow 
    X_0 \neq \emptyset 
    \Leftrightarrow 
    y^\dagger  \in R(A)
    \Leftrightarrow 
    \argmin f \neq \emptyset.$
For items \ref{P:Regularity spaces in terms of mu:i}
and \ref{P:Regularity spaces in terms of mu:ii}, the claim follows directly from the nonincreasingness of $\{X_\mu\}_{-1/2< \mu < + \infty}$.
For item \ref{P:Regularity spaces in terms of mu:iii}, observe that for all $\nu >0$,  $\spec((AA^*)^\nu) = \spec(AA^*)^\nu$ \cite[\textsection 32 Theorem 3]{Hel69}.
As a consequence of Proposition \ref{P:closed range singular values}, we deduce that $R (AA^*)^\nu$ is closed, 
and therefore $R(AA^*)^\nu=R(A)$ (see Lemma \ref{L:powers self adjoint operators} in the Annex).
In particular, $Y_\nu=y^\dagger + R(A)$ for all $\nu \in \left]0,+\infty\right[$, and the result follows from item \ref{P:Regularity spaces in terms of mu:ii}.
\end{proof}

For well-posed problems, for which $\argmin f \neq \emptyset$ (and $x^\dagger$ exists), the source sets can be  expressed with a simpler expression (the proof is left in the Annex):

\begin{lemma}[Source sets for well-posed problems]\label{L:regularity trough data space}
Assume that $\argmin f \neq \emptyset$.
Then, for all $\mu >0$ and $\delta >0$:
\begin{equation*}
    X_\mu  =  \{ x^\dagger \} + \Ker A + R(A^*A)^\mu
    \ \text{ and } \
    X_{\mu,\delta} =  \{ x^\dagger \} + \Ker A + \{ (A^*A)^\mu w \ | \  w  \in \Ker A^\perp, \Vert w \Vert \leq \delta \} .
\end{equation*}
\end{lemma}
}

\begin{remark}
Given that $x^\dagger \in \ker A^\perp$, we see that the classical  conditions in \eqref{D:source condition} and \eqref{D:source condition in Y} are  equivalent, with our notations, to $0 \in X_\mu$ and $0 \in X_{\nu -1/2}$.	
This means in particular that \eqref{D:source condition} is just a particular case of \eqref{D:source condition in Y}.
\end{remark}

{
\begin{remark}[Source sets as balls]\label{R:source sets as balls}
Assume that $A$ is injective and $y\in R(A)$.
For all $\mu>0$, $R(A^*A)^\mu$ is a dense subspace of $X$ (Lemma \ref{L:powers self adjoint operators}), and we can endow it with the norm induced by the unbounded operator $(A^*A)^{-\mu}$, defined as $\Vert x \Vert_\mu:= 
    \inf \{ \Vert w \Vert \ | \ w \in X \text{ and } 
    x = (A^*A)^\mu w \}$.
Then, we see that the source sets $X_{\mu,\delta}$ are nothing but balls centered at the 
solution $ x^\dagger$, with respect to this norm:
\begin{equation*}
    X_{\mu,\delta} = \{ x \in X \ | \ \Vert x - x^\dagger \Vert_\mu \leq \delta \},
\end{equation*}
while $X_\mu$ is the affine space spanned by these balls.
By doing an analogy with the following example, the  reader can think about this norm $\Vert \cdot \Vert_\mu$ in $X$ as if it was a Sobolev norm in an $L^2$ space.
Note that these balls may have an empty interior with respect to the topology of $X$.
\end{remark}
}

\begin{example}[Regularity spaces as Sobolev spaces]
Assume that $X$ is the space of zero mean $L^2$-functions on $[0,2\pi]$:
\begin{equation*}
X=\left\{ \varphi \in L^2([0,2\pi]), \ \int_0^{2\pi} \varphi(t) \ \dt =0 \right\}.
\end{equation*}
If $A$ is the linear integration operator defined on $X$, then  $R (A^*A)^\mu$ coincides  with the Sobolev space $H^{2\mu}([0,2\pi]) \cap X$ \cite[Theorem 6.4]{Hoh02}, so that the regularity space is here
\begin{equation*}
X_\mu = \{x^\dagger \} + H^{2\mu}([0,2\pi]) \cap X.
\end{equation*}
\end{example}

\subsubsection{Properties of quadratic functions on source sets}
Here is the main result of this section: on each source set $X_{\mu,\delta}$, the least squares 
functional $f$ is $p$-{\L}ojasiewicz with $p=2+\mu^{-1}$.
\begin{theorem}[Geometry of least squares on source sets]\label{T:geometry of least squares}
Let $\mu \in \left]-1/2,0\right[\cup \left]0,+\infty\right[$ and $\delta \in \left]0,+\infty\right[$.
Then $f(x)=\frac{1}{2}\Vert Ax-y\Vert^2$ is $p$-{\L}ojasiewicz on $\Omega=X_{\mu,\delta}$, with 
\begin{equation}\label{T:geometry of least squares:constants}
\text{$p=2 + \mu^{-1}$ and $c_{f,\Omega}=2^{-(\mu +1)/(2\mu + 1)} \delta^{1/(1+2\mu)}$.}
\end{equation}
Moreover, these two constants are sharp.
\end{theorem}
\begin{proof}
Let $x \in X_{\mu ,\delta}$ and remind that $y^\dagger =\proj(y, \cl R(A))$.
From Definition \ref{D:regularity spaces and source sets} and the definition of $y^\dagger$, we get
\begin{eqnarray}
Ax= y^\dagger + (AA^*)^{\mu + 1/2} \omega, \text{ where } \omega \in {\ker A^*}^\perp \text{ with } \Vert \omega \Vert \leq \delta, \label{klp1} \\
f(x) - \inf f = (1/2)\Vert Ax - y^\dagger \Vert^2 \text{ and } \Vert \nabla f(x) \Vert = \Vert A^*(Ax - y^\dagger) \Vert.\label{klp1.5}
\end{eqnarray}
We first prove that $f$ verifies the {\L}ojasiewicz inequality by using the interpolation inequality (see Lemma \ref{L:interpolation inequality} in 
the Annex) with $\alpha=\mu+ (1/2)$ and $\beta= \mu +1$, together with \eqref{klp1}:
\begin{equation}\label{klp2}
\hspace{-0.18cm}\Vert Ax- y^\dagger \Vert =
 \Vert (AA^*)^{\mu + 1/2} \omega \Vert \leq 
 \Vert (AA^*)^{1+\mu} \omega \Vert^\frac{2\mu+1}{2\mu +2} \Vert \omega \Vert^\frac{1}{2\mu+2} \leq
 \Vert (AA^*)^{1+\mu} \omega \Vert^\frac{2\mu+1}{2\mu +2} \delta^\frac{1}{2\mu+2}.
\end{equation}
We  use \eqref{klp1} in the right member of \eqref{klp2}, to write
\begin{equation}\label{klp2.5}
\hspace{-0.09cm}\Vert (AA^*)^{1+\mu} \omega \Vert^2= \Vert (AA^*)^{1/2}(AA^*)^{\mu+1/2} \omega\Vert^2 = \Vert (AA^*)^{1/2} (Ax- y^\dagger) \Vert^2 = \Vert A^*(Ax- y^\dagger) \Vert^2.
\end{equation}
By combining \eqref{klp1}, \eqref{klp1.5}, \eqref{klp2} and \eqref{klp2.5}, we obtain the following inequality
\begin{equation*}\label{klp3}
f(x) - \inf f = (1/2) \Vert Ax - y^\dagger \Vert^2 \leq (1/2) \delta^\frac{1}{\mu +1} \Vert A^*(Ax - y^\dagger) \Vert^\frac{2\mu+1}{\mu +1}
= (1/2) \delta^\frac{1}{\mu +1} \Vert \nabla f(x) \Vert^\frac{2\mu+1}{\mu +1}.
\end{equation*}
Then the desired \Loja inequality holds by taking $p:=2+\mu^{-1}$.
Now we verify that the obtained constants in \eqref{T:geometry of least squares:constants} are sharp.
For this, let $X=\ell^2(\N)$, and let $(e_k)_\kin \subset X$ be its canonical basis. 
Let $(\sigma_k)_\kin$ be a {strictly} positive sequence converging to zero, and define $A : X \longrightarrow X$ as follows: $\forall x=(x_k)_\kin \in X, Ax:= \sum_\kin \sigma_k x_k e_k$. 
Let $f(x)=(1/2)\Vert Ax \Vert^2$, $y= 0$, and let us assume that $f$ is $p$-{\L}ojasiewicz on $X_{\mu,\delta}$ for some $p\geq 1$:
\begin{equation}\label{ocg1}
(\forall x \in X_{\mu,\delta}) \quad [(1/2)\Vert Ax \Vert^2]^{1-(1/p)} \leq c_{f,X_{\mu,\delta}} \Vert \AA x \Vert.
\end{equation}
Let $v^k:=\delta \sigma_k^{2\mu} e_k \in X_{\mu,\delta}$, which satisfies $\Vert \AA v^k \Vert =\delta \sigma_k^{2 +2\mu} $, and deduce from \eqref{ocg1} that
\begin{equation}\label{ocg2}
(\forall \kin)\quad 2^{(1-p)/p} \delta^{(2p-2)/p} \leq c_{f,X_{\mu,\delta}} \sigma_k^{(1/p)(4\mu + 2) - 2\mu} \delta.
\end{equation}
It follows from $\sigma_k \to 0$ that  $4\mu - 2\mu p +2\leq 0$, which is equivalent to  $\mu p \geq 2\mu + 1$.
If $\mu >0$, it means that $ p \geq 2 + \mu^{-1} >0$, which is a regime in which the smallest is $p$, the better.
If $\mu \in ]-1/2,0[$, then $ p \leq 2 + \mu^{-1} < 0$, which is a regime in which the largest is $p$, the better.
In both cases we see that $2 + \mu^{-1}$ is the best possible exponent.
Moreover, when $p=2 + \mu^{-1}$, \eqref{ocg2} becomes
$ 2^{-\frac{1+\mu}{1+2\mu}} \delta^{\frac{1}{1+2\mu}} \leq c_{f,X_{\mu,\delta}},$
which implies the sharpness of the constant  obtained in  \eqref{T:geometry of least squares:constants}.
\end{proof}
\begin{remark}
The result of Theorem \ref{T:geometry of least squares} contrasts with \cite[Theorem 2.1]{HarJen11}, in which the authors show that 
no local {\L}ojasiewicz property can be satisfied by a quadratic function when $R(A)$ is not closed.
The key difference here is that we look at the {\L}ojasiewicz property on specific \textit{dense} sets with empty interior (see Remark \ref{R:source sets as balls}).
\end{remark}
Let us now verify that the source sets are invariant under the action of the Landweber algorithm \eqref{D:Landweber}. 
As mentioned at the beginning of the section, the Landweber algorithm is the gradient decent algorithm applied to a quadratic function, and therefore it is an instance of the FB algorithm. We can thus apply the convergence rates of Section \ref{S:CV rates for Forward Backward} once we prove that the source sets are invariant. 
\begin{proposition}[Invariance of source sets]\label{P:stability of source condition}
For all $(\mu,\delta) \in \left]-1/2,\infty\right[\times \left]0,+\infty\right[^2$, the source set $X_{\mu,\delta}$ is FB-invariant.
\end{proposition}
\begin{proof}
Let $x \in  X_{\mu,\delta}$, $\lambda\in\left]0,2/\| \AA \|\right[$, and let us prove that 
$\T x= x - \lambda A^*(Ax - y)$
belongs to $X_{\mu,\delta}$.
By using Lemma~\ref{L:regularity trough data space}, we deduce that 
$Ax= y^\dagger + (AA^*)^\nu \omega$,  $\nu:= \mu+1/2$, and $\omega \in \cl R(A)$ with $\Vert \omega \Vert \leq \delta$.
Since $A^*(Ax-y)=A^*(Ax-y^\dagger)$, this implies that
\begin{equation*}
A\T x  = Ax - \lambda AA^*(Ax- y^\dagger) = y^\dagger + (AA^*)^\nu (I - \lambda AA^*)\omega 
\end{equation*}
The above equality shows that $\T x \in X_{\mu}$. 
It remains only to prove that $\hat{T}_\lambda \omega:= (I - \lambda AA^*) \omega$ verifies $\hat{T}_\lambda \omega \in \cl R(A)$ and $\Vert \hat{T}_\lambda \omega \Vert \leq \delta$.
The condition $\hat{T}_\lambda \omega \in \cl R(A)$ immediately follows from $\omega \in \cl R(A)$ and $AA^*\omega \in R(A)$.
Next, observe that $\hat{T}_\lambda \omega$ is obtained by applying a gradient descent step to $\omega$ with respect to the function $u \mapsto (1/2)\Vert A^*u \Vert^2$.
Since this function has zero as a minimizer, and is differentiable with a $\| \AA \|$-Lipschitz gradient,  the Fej\'er property (see Theorem~\ref{T:CV FB}-\ref{T:CV FB:ii}) implies that $\Vert \hat{T}_\lambda \omega \Vert \leq \Vert \omega \Vert \leq \delta$.
\end{proof}
Next we  combine all the results of this section to derive convergence rates of the Landweber algorithm under source conditions from \Loja conditions. 
\begin{corollary}[Convergence rates for Landweber algorithm]
\label{Cor:coLan}
Let $(x_n)_\nin$ be a sequence generated by the Landweber algorithm \eqref{D:Landweber}.
Assume that for some $\mu \in \left ]-1/2,+\infty\right[$, the source condition $x_0 \in X_\mu$ is satisfied. Then:
\begin{enumerate}[i)]
	\item \label{Cor:coLan:i} $f(x_n) - \inf f =O( {n^{-(1+2\mu)}})$,
	\item  \label{Cor:coLan:ii} If $\mu > 0$, then  $\Vert x_n - \bar x_0 \Vert = O(n^{-\mu})$, where $\bar x_0:= \proj(x_0,\argmin f)$.
\end{enumerate}
\end{corollary}

\begin{proof}
For item \ref{Cor:coLan:i}, the source condition together with Proposition \ref{P:stability of source condition} imply $(x_n)_\nin \subset X_{\mu,\delta}$ for some $\delta >0$.
If $\mu \neq 0$,  we derive from Theorem \ref{T:geometry of least squares} that $f$ is $2+\mu^{-1}$-{\L}ojasiewicz on $X_{\mu,\delta}$.
Depending on the sign of $2+\mu^{-1}$, the rates on $f(x_n) - \inf f$ follow from Theorems~\ref{T:CVKL discrete} and \ref{T:CVKL discrete rates negative p}. 
If $\mu = 0$, then the source condition and Proposition \ref{P:Regularity spaces in terms of mu} ensures that $y^\dagger \in R(A)$, meaning that $\argmin f \neq \emptyset$, so the rate $O(n^{-1})$ follows from Theorem~\ref{T:CV FB}. 
For item \ref{Cor:coLan:ii}, the convergence and rates on the iterates follows from Theorem~\ref{T:CVKL discrete}.
To show that the limit of the sequence (let us note it $x_\infty$) is $\bar x_0$, it is enough to verify that $x_\infty -x_0 \in \ker A^\perp$, since $\argmin f$ is an affine space parallel to $\ker A$.
Because of the definition of the algorithm, it is easy to show by recurrence that $x_n - x_0 \in R(A^*)$. 
This being true for all $n \in \N$, we can pass to the limit and deduce that $x_\infty - x_0 \in \cl R(A^*) = \ker A^\perp$.
\end{proof}

\subsection{Sparsity based regularization, partial smoothness, and restricted injectivity}\label{SS:sparse inverse problems}

{
In this section we turn to the general case of optimization problems coming from a regularized inverse problem \eqref{D:least squares inverse problems}.
In particular, we focus on the case where $\nabla^2 h$ verifies a restricted injectivity condition at a solution, a situation which typically arises when $g$ is mirror-stratifiable, and typical modeling assumptions from the inverse problems/compressed sensing literature hold.
In this setting we will be able to derive the 2-conditioning of the objective function in \eqref{D:least squares inverse problems}.
In what follows, we will use the notation $\mathcal{S}_+(X)$ to refer to the set of bounded selfadjoint positive linear operators on $X$.

\subsubsection{Coercive linear operators on a cone}
\begin{definition}
We say that $K \subset X$ is a cone if it is a union of rays: $[0,+\infty[ K \subset K$.
\end{definition}

Note that we do not require a cone to be convex. 
This is important for certain applications in which we have geometrical information about a function over a union of linear spaces, see for instance \eqref{D:RIP partial} in the context of sparse regularization problems.

\begin{definition}\label{D:coercive operator}
Let $S\in \mathcal{S}_+(X)$, let $\gamma\in\left]0,+\infty\right[$, and let $K \subset X$ be a cone.
We say that $S$ is $\gamma$-coercive on $K$ if, for all $d \in K$, $\langle Sd, d \rangle \geq \gamma \Vert d \Vert^2$.
\end{definition}

\begin{example}[coercivity for positive symmetric matrices]
\label{Ex:elliptic operator facts}
A matrix $S \in \mathcal{S}_+(\mathbb{R}^N)$ is coercive on a closed cone $K \subset \mathbb{R}^N$ 
if and only if $S$ is injective when restricted on $K$ (see Proposition \ref{P: matrix coercive on injective cone} for a proof):
\begin{equation*}
    K \cap \mbox{Ker}~S = \{0 \}.
\end{equation*}
\end{example}

\begin{example}
Any operator $S \in \mathcal{S}_+(X)$ is $\sigma_{\inf}(S)$-coercive on
${\mbox{Ker}}~S^\perp$ (see e.g. the proof in \cite[Thm. 4]{CraGoc20}).
In particular, if $S$ is positive definite then it is $\sigma_{\inf}(S)$-coercive on $X$.  
\end{example}
}

In the next proposition we relate the coercivity of the Hessian of a function $f$ on a cone to the $2$-conditioning of $f$ on this cone.
This relation can be seen as a weakened analogue of the well known fact (see \cite[Prop. 10.8 \& 17.7.(iii)]{BauCom}) that, for $f \in C^2(X)$ : 
\begin{center}
$f$ is $\gamma$-strongly convex $\Leftrightarrow$ ($\forall x \in X$) $\nabla^2 f(x)$ is $\gamma$-coercive on $X$.
\end{center}
Strong convexity is a global notion, which requires the function to have a \textit{positive definite} quadratic-like geometry at each $x\in X$.
On the contrary, the $2$-conditioning requires the function to have a \textit{positive} quadratic-like geometry, on a given set $\Omega$.
We now state our result (its proof is left in the Annex \ref{SS:Annex proof section 5.2}). For similar results, see also \cite[Section 3.3.1]{BonSha00} and  \cite{DruMorNgh14}.

\begin{proposition}[Coercivity of the Hessian implies $2$-conditioning]
\label{P:ellipticity implies conditioning}
Let $f=g+h$ with $g,h \in \Gamma_0(X)$ and $\argmin f \neq \emptyset$. Assume that $h$ is of class $C^2$ in a neighbourhood of $\bar x \in \argmin f$,
and that $\nabla^2 h( \bar x)$ is $\gamma$-coercive on a closed cone $K\subset X$.
Then,
\begin{equation*}
\text{$(\forall \gamma' \in ]0,\gamma[)$ $(\exists \delta \in ]0,+\infty])$ s.t. $f$ is $2$-conditioned on $\Omega:=\bar x + (K \cap \delta \B_X )$ with  $\gamma_{f,\Omega}=\gamma'$, }
\end{equation*}
and $\Omega \cap \argmin f = \{ \bar x \}$.
If $h\in C^2(X)$ and $\nabla^2 h$ is $L$-Lipschitz, we can take $\delta=\frac{\gamma - \gamma'}{L}$.
\end{proposition}

{
\subsubsection{Conditioning on prox-regular sets via restricted injectivity of the Hessian}
Let us define some useful tools from variational analysis.
The notion of reached set (or set with positive reach) was introduced by Federer \cite[Def. 4.1]{Fed59}, and later extended to  \textit{prox-regularity} (see Proposition \ref{P:curvature manifold} and \cite{RocWet}).

\begin{definition}
Let $C \subset \mathbb{R}^N$.
The (Bouligand) tangent cone to $C$ at $\bar x \in C$ is defined as 
\begin{equation*}
    T_C(\bar x) :=  \{d \in \mathbb{R}^N \ | \ (\exists t_n \downarrow 0)(\exists d_n \rightarrow d)  \quad \bar x + t_n d_n \in C \}.
\end{equation*}
The normal cone to $C$ at $\bar x$ is $N_C(\bar x) := \{ \eta \in \mathbb{R}^N \  | \ (\forall d \in T_C(\bar x)) \ \langle \eta,d  \rangle \leq 0 \}$.
\end{definition}

\begin{definition}\label{D:prox regular set}
Let $C \subset \mathbb{R}^N$, and $\rho >0$.
We say that $C$ is $\rho$-reached at $\bar x \in C$, if it is locally closed at $\bar x$, and verifies
\begin{equation*}
    (\forall \eta \in N_{C}(\bar x) \cap \mathbb{S}_{\mathbb{R}^N}) \quad \mathbb{B}(\bar x + \frac{1}{\rho} \eta, \frac{1}{\rho}) \cap {C} = \emptyset.
\end{equation*}
We say that $C$ is prox-regular at $\bar x$ if there exists $\rho>0$ and a closed neighbourhood $U$ of $\bar x$ such that $C\cap U$ is $\rho$-reached at any $x \in U$.
We say further that $C$ is prox-regular if it is prox-regular at every $\bar x \in C$.
\end{definition}
Convex sets, and in particular affine spaces, are prox-regular. Manifolds of class $C^2$ are locally prox-regular (see Proposition \ref{P:curvature manifold}).

We now provide the result at the core of this section, which says that if a minimizer $\bar x$ belongs to some prox-regular set, and if the Hessian $\nabla^2 h(\bar x)$ is injective when restricted to the tangent cone of this set, then $f$ is $2$-conditioned on this set around $\bar x$.
This will guarantee asymptotic linear rates when combined with Corollary \ref{T:partial smoothness rates}.

\begin{theorem}[Injective Hessian on tangent cone implies $2$-conditioning]\label{T:injective Hessian implies 2 conditioning}
Let $g,h \in \Gamma_0(\mathbb{R}^N)$, and $f=g+h$.
Assume that there exists some $\bar x \in \argmin f$ such that:
\begin{enumerate}[label=\alph*)]
    \item\label{T:injective Hessian implies 2 conditioning:minimizer localisation}  $\bar x$ belongs to some $C\subset \mathbb{R}^N$ which is $\rho$-reached at $\bar x$,
    \item\label{T:injective Hessian implies 2 conditioning:C2} $h$ is of class $C^2$ in a neighbourhood of $\bar x$,
    \item\label{T:injective Hessian implies 2 conditioning:injectivity}  $\mbox{Ker}~\nabla^2 f(\bar x)$ is $\gamma$-coercive on  $T_C(\bar x)$.
\end{enumerate}
Then ${\rm{argmin}}~f_{|C} = \{ \bar x \}$, and for every $\gamma' \in ]0,\gamma[$, there exists $\delta \in ]0,+\infty]$ such that $f$ is $2$-conditioned on $\Omega := C \cap \mathbb{B}(\bar x, \delta)$, with $\gamma_{f,\Omega}=\gamma'$.
If we assume moreover that $\nabla^2 h$ is $L$-Lipschitz continuous, then we can take $\delta = \frac{2(\gamma - \gamma')}{2L + \rho \Vert \nabla^2 h(\bar x) \Vert}$.
\end{theorem}

\begin{proof}
Let $K:= T_C(\bar x)$.
Using Proposition \ref{P:ellipticity lifted from tangent to manifold}, we see that for every $\gamma' < \gamma$ there exists a $\theta \in ]0, \frac{\pi}{2}[$ such that the enlarged cone $K_\theta$ (see Definition~\ref{def:enlcone}) contains $(C - \bar x) \cap \mathbb{B}(0, \delta)$ for $\delta >0$ small enough, and such that $\nabla^2 h(\bar x)$ is $\gamma'$-coercive on $K_\theta$.
The conclusion of the claim follows from Proposition \ref{P:ellipticity implies conditioning} applied to $h$ and $K_\theta$.
Under the additional assumption that  $\nabla^2 h$ is $L$-Lipschitz, take any $\gamma' \in ]0, \gamma[$, and let $\gamma'' := \alpha \gamma + (1- \alpha) \gamma'$, with $\alpha = 2L/(2L + \rho \Vert \nabla^2 h(\bar x) \Vert)$.
Using again Proposition \ref{P:ellipticity lifted from tangent to manifold}, we obtain that $\nabla^2 h(\bar x)$ is $\gamma''$-coercive on some cone $K_\theta$, with $\bar x + K_\theta \supset C\cap \mathbb{B}(\bar x, \delta_1)$ and $\delta_1 = 2(\gamma - \gamma'')/ (\rho\Vert \nabla^2 h(\bar x) \Vert)$.
Then, Proposition \ref{P:ellipticity implies conditioning} shows that $f$ is $2$-conditioned on $\Omega = \bar x + K_\theta \cap \mathbb{B}(\bar x, \delta_2)$, with $\delta_2 = (\gamma''-\gamma')/L$ and $\gamma_{f,\Omega} = \gamma'$.
The conclusion follows by seeing that  $\delta_1 = \delta_2$ with our choice of $\gamma''$.
\end{proof}

Theorem \ref{T:injective Hessian implies 2 conditioning} can be used in combination with Corollary \ref{T:partial smoothness rates}: in this case we obtain that the restricted injectivity of the Hessian on the tangent cone to the active set $C_{\bar x}$ guarantees asymptotic linear rates.
In the example below, we detail what our assumptions mean for the examples in Example \ref{Ex:mirror active sets}.

\begin{example}\label{Ex:mirror active sets tangent proxregular}~~
\begin{itemize}
    \item If $g(x) = \Vert x \Vert_1$, the active set \eqref{Ex:mirror active sets:L1 support} is an open and dense subset of the vector space $X_I = \{x \in \mathbb{R}^N \ | \ \supp(x) \subset I \}$ with $I = \act(-\nabla h(\bar x))$.
    It is therefore $\rho$-reached for every $\rho>0$, and $T_{C_{\bar x}}(\bar x) = X_I$.
\item If $g(x) = \Vert x \Vert_*$, let $r=\#\act(\sigma(-\nabla h(\bar x)))$ and let $\mathcal{M}_r$ be the manifold of matrices with rank equal to $r$. 
If $0 \in \rint \partial f(\bar x)$, the active set $C_{\bar x}$ (see \eqref{Ex:mirror active sets:nuclear}) is equal to $\mathcal{M}_{r}$.
In particular, it is prox-regular (see Proposition \ref{P:curvature manifold}), and an expression for its tangent space can be found in \cite[Example 2.2]{LewMal08}.
More generally, $C_{\bar x}$ is locally prox-regular at $\bar x$ if $\rank(\bar x) = r$.
To see this, use the same arguments as in \cite[Prop. 3.1]{Luk13}: the fact that the singular values depend continuously on the matrix allows to find a neighbourhood $U$ of $\bar x$ where the matrices have a rank greater or equal to $r$.
This means that $C_{\bar x}\cap U=\mathcal{M}_r \cap U$, which is  prox-regular.
\end{itemize}

\end{example}
\begin{remark}[Related results with partial smoothness]
While our results are new in the setting of mirror-stratifiable functions (where no condition $0 \in \rint \partial f(\bar x)$ is required), they intersect with existing results when $g$ is partially smooth with respect to an active manifold $\mathcal{M}$.
It is shown in \cite{LiaFadPey14} that the $\gamma$-coercivity of $\nabla^2 h(\bar x)$ on the tangent space $T_\mathcal{M}(\bar x)$ guarantees asymptotic linear rates.
We recover a similar result by combining Theorem \cite[Theorem 5.3]{HarLew04} with Theorem \ref{T:injective Hessian implies 2 conditioning} and Theorem \ref{T:CV FB}.
For a fixed stepsize $\lambda =1/L$, \cite[Thm. 3.1]{LiaFadPey14} predicts a Q-linear rate arbitrarily close to $\sqrt{2(1-\kappa)}$ (where $\kappa = \gamma/L$) provided that $\kappa \geq 1/2$.
Instead, our results predict a R-linear rate arbitrarily close to $(1+ (\kappa/4))^{-1/2}$, without condition on $\kappa$. 
Note that our constant is worse (resp. better) than $\sqrt{2(1-\kappa)}$ when $\kappa$ is close to $1$ (resp. $1/2$).
Note also that the partial smoothness of $g$ together with \cite[Theorem 6.2.ii)]{HarLew04} ensures that $f$ is $2$-conditioned on a neighbourhood $\Omega$ of the solution, with $\gamma_{f,\Omega} = \gamma'$, meaning that we can use Proposition \ref{P:equivalence linear rates and 2-conditioning} to obtain Q-linear rates arbitrarily close to $(1+ \kappa)^{-1/2}$.
\end{remark}

}
\subsubsection{Application to low-complexity inverse problems}

Consider $f \in \Gamma_0(\R^N)$ be defined by, for every $x\in\R^N$,  $f(x) = \alpha \Vert x \Vert_1 + (1/2)\Vert Ax - y \Vert^2$.
$f$ is the sum of a smooth function, with Hessian equal to $A^*A$, and a nonsmooth function $ \alpha \Vert x \Vert_1$.
Example \ref{Ex:LASSO is 2 conditioned}  ensures that $f$ is locally $2$-conditioned on its sublevel sets \textit{without any assumption} on $A$. 
This means, according to Theorem \ref{T:CVKL discrete}, that for any $r>\inf f$, and any $x_0 \in [f<r]$, there exists a  constant $\eps \in ]0,1[$ such that the iterative soft-thresholding initialized at $x_0$ verifies
$f(x_{n+1}) - \inf f \leq \eps (f(x_n)  - \inf f)$.
Nevertheless, expressing the $2$-conditioning constant, or $\eps$, in terms of the components of the problems is far to be easy \cite{BolNguPeySut15}.
One way to recover a meaningful constant is to exploit modeling assumptions which are usually made to ensure the stability and recovery of the inverse problem $Ax=y$.

Suppose that we are given the sequence generated by the iterative soft-thresholding, which converges to a minimizer of $f$, $x_n \rightarrow \bar x$.
It is known that, after some iterations, the support of the sequence is stable \cite{LiaFadPey17,GarRosVil20}:
\begin{equation*}
(\exists I \subset \{1,\ldots,N \}) (\exists n_0 \in N)  (\forall n \geq n_0) \quad \supp(x_n) \subset I.
\end{equation*}
In particular, if the qualification condition $0 \in \rint \partial f(\bar x)$ holds, we can take $I=\supp(\bar x)$ \cite[Prop. 3.6]{LiaFadPey17}.
To estimate the rates of convergence for the sequence, it is then sufficient to make a \textit{restricted injectivity} 
assumption on the matrix $A$, depending on the knowledge we have on $I$.

In the case we have access to $I$, suppose that on the space $X_I:=\{x \in \R^N \ | \ \supp(x) \subset I \}$ the matrix $A$ is injective, i.e. $\Ker A \cap X_I = \{0 \}$ holds.
Then, there exists a constant $\gamma_I >0$ such that $A^*A$ is $\gamma$-coercive on $X_I$ (see Example \ref{Ex:elliptic operator facts}), which implies via Proposition \ref{P:ellipticity implies conditioning} that $f$ is $2$-conditioned on $X_I$, with $\gamma_{f,X_I}=\gamma_I$. 
We deduce then that, asymptotically, the rates are governed by $\eps= (1+\gamma_I \Vert A^*A \Vert^{-1})^{-1}$.
It might happen that instead of knowing $I$, we have only access to a partial information via the sparsity level $s:= \vert I \vert$.
We can then follow the same reasoning with the (nonconvex) cone $K_s:= \{x \in \R^N \ | \ \vert \supp(x) \vert \leq s \}$ instead of $X_I$.
In that case, the constant $\gamma_s$ of coercivity of $A^*A$ on $K_s$ is defined by
\begin{equation}
\label{D:RIP partial}
(\forall x \in K_s)\quad \gamma_s \Vert x \Vert^2 \leq \Vert Ax \Vert^2,
\end{equation}
and guarantees linear rates governed by $\eps= (1+\gamma_s \Vert A^*A \Vert^{-1})^{-1}$, using again Proposition \ref{P:ellipticity implies conditioning}.
Such assumption is classical in sparsity based regularization, and it is related to the so-called Restricted Isometry Property \cite{Can08},
to ensure uniqueness of the minimizer and 
guarantee the robustness or recovery \cite{VaiPeyFad14,ChaRecParWil12}.
Observe that while the computation of $\gamma_s$ remains impracticable \cite{BanDobMixSaw13}, it is \textit{meaningful} with respect to the properties of our problem, and, more importantly, can be estimated when the matrix $A$ is random \cite[Section 9]{FouRau}.
Of course, this whole discussion can be extended to other regularized inverse problems, in particular if $\Vert \cdot \Vert_1$ is replaced by a mirror-stratifiable function.
In this case we will use Theorem \ref{T:injective Hessian implies 2 conditioning} instead of Proposition \ref{P:ellipticity implies conditioning} to derive linear rates.

\section{Conclusion and perspectives}

In this paper, we dicussed in details how geometry can be used to improve the rates of the FB method, or more general first-order descent schemes.
We characterized the geometry, using tools that  are  often encountered in practice, like the $p$-conditioning, and we provided a new sum rule for it.
In Figure \ref{F:spectra of rates} we recall the various  rates obtained for the FB method, 
 from the worst case scenario (no minimizers, no assumptions) to the best one (sharp functions).

\begin{figure}[h]
\begin{center}
\begin{tabular}{|c|c|c|c|}
\hline 
 & $f(x_n) - \inf f$ & $\Vert x_n - x_\infty \Vert$ 	 \\ 
\hline 
$\inf f > - \infty$ & $o(1)$ & ---  \\ 
\hline 
$p \in \left]-\infty,0\right[$ & $O(n^{p/(2-p)})$ & ---   \\ 
\hline 
$\argmin f \neq \emptyset$ & $o(n^{-1})$ & decreasing, $o(1)$ in finite dimension  \\ 
\hline 
$p \in \left]2, + \infty\right[$ & $O(n^{-p/(p-2)})$ & $O(n^{-1/(p-2)})$  \\ 
\hline 
$p=2$ & Q-linear with $\eps=1/(1+\kappa)$ & R-linear with $\eps=1/(1+\kappa)$  \\ 
\hline 
$p\in \left]1,2\right[$ & Q-superlinear of order $1/(p-1)$ & R-superlinear of order $1/(p-1)$  \\ 
\hline 
$p=1$ & finite & finite  \\ 
\hline 
\end{tabular} 
\end{center}

\caption{Convergence rates of the FB algorithm for locally $p$-{\L}ojasiewicz functions (with the constant $\kappa$ defined in Theorem \ref{T:CVKL discrete}).}
\label{F:spectra of rates}
\end{figure}

\noindent We also  have discussed how  those refined results can be obtained by decoupling the geometrical information we have on the function and the localization of the sequence we are looking at.
This geometry-based analysis reduces then the gap between theory and practice, where the observed rates are often better than the ones resulting from a worst case analysis. It moreover shows that linear rates are tightly linked to $2$-conditioned function.
In addition, we showed how our analysis can be specialized to the inverse problems setting, and allows to explain typical modeling assumptions in this context, such as source conditions and restricted injectivity property. 
It is worth noting that the geometrical information such as conditioning or {\L}ojasiewicz property can be exploited to derive sharper convergence rates  for a broader class of functions and/or algorithms than just forward-backward algorithm \cite{AttBolSva13}.
We also emphasize that convexity plays no role in the proofs of Theorems \ref{T:CVKL discrete} and \ref{T:CVKL discrete rates negative p}. Indeed, some of these results were already known for non-convex functions \cite{BolSabTeb13,ChoPesRep14,FraGarPey15}. 
One of the challenges in the future is to have quantitative results concerning the geometry of classes of nonconvex functions.
For instance, what can be said about ``simple'' nonconvex piecewise polynomial functions (see \cite{LiMorPha15} for an answer about maximum of finitely many polynomials)?
Can we estimate the {\L}ojasiewicz exponent of semialgebraic functions, depending on the degree of the polynomials defining their graph?
Finally, a last challenge is the application of such geometrical tools to derive precise rates for nondescent methods.
First results  in this direction, using  $2$-conditioning are known for  inertial methods \cite{NecNesGli15,LiaFadPey16} or stochastic gradient methods \cite{KarNutSch16}.
It would be  of interest to understand the behavior of these algorithms for other geometries.

\appendix

\section{Appendix}

\subsection{Worst case analysis: proofs of Section~\ref{S:classic FB} }
\label{SS:proof worst case}

{
The following Lemma contains a detailed proof for the lower bound  \eqref{E:rates for the norm to the -p values} in Example \ref{R:optimality of rates worst case}, which can also be applied to \eqref{E:rates for the norm to the p values} by using a symmetry argument.

\begin{lemma}[Lower bounds for the proximal algorithm]\label{L:lower bounds proximal}
Let $p \in ]-\infty,0[ \cup ]2,+\infty[$, and let $f_p \in \Gamma_0(\R)$ be the function defined by
\begin{equation*}
\text{ if } p<0, \, f_p(x) = 
\begin{cases}
\vert x \vert^p  & \text{ if } x <0,\\
+\infty & \text{ if } x \geq 0,
\end{cases}
\quad
\text{ and if } p>2, \ f_p(x) =
\begin{cases}
0 & \text{ if } x <0, \\
\vert x \vert^p & \text{ if } x \geq 0.
\end{cases}
\end{equation*}
If $x_0 \in \dom f \setminus \argmin f$, and $x_{n+1} = \mbox{prox}_{\lambda f}(x_n)$, then for all $n \geq 1$:
\begin{equation*}
f_p(x_n) - \inf f_p \geq C_p^p n^{\frac{p}{2-p}}
\quad \text{ with } \quad
C_p= \left( \vert x_0\vert ^{2-p} + p (p-2) \lambda \right)^{\frac{1}{2-p}}.
\end{equation*}
\end{lemma}

\begin{proof}
Note that $\dom f_p$ is an open interval, and that $f_p$ is infinitely derivable there.
We can then see that $f_p$, $f_p'$ and $f_p''$ are non-negative.
In particular, we deduce that $f_p$ and $f'_p$ are non-decreasing on $\dom f$.

Let us now take  some $x_0 \in \dom f\setminus \argmin f$,
and consider the following continuous trajectory
\begin{equation*}
(\forall t \geq 0) \quad x(t) := \sgn(p) \left( \vert x_0\vert ^{2-p} + p (p-2) t \right)^{\frac{1}{2-p}}.
\end{equation*}
It is a simple exercise to verify that $x(\cdot)$ is a solution of this  differential equation:
\begin{equation*}
x(0)=x_0, \quad \dot{x}(t) + f'(x(t)) = 0, \quad x(t) \in \dom f_p.
\end{equation*}
The main step towards proving our lower bound is to show, by induction, that for every $n \in \N$, $x_n \geq x(n\lambda)$. 
This is clearly true for $n=0$, so, let us assume now that this is true for $n \in \N$, and show that this implies $x_{n+1} \geq x((n+1)\lambda)$.
Start by writing
\begin{equation*}
x((n+1)\lambda) = 
x(n\lambda) + \int_{n\lambda}^{(n+1)\lambda} \dot x(t)  ~ dt = 
x(n\lambda) + \int_{n\lambda}^{(n+1)\lambda} (-f_p' \circ x)(t) ~  dt.
\end{equation*}
On the one hand, $f_p'$ is non-negative on $\dom f$, and $\dot x(t) = -f_p'(x(t))$, which means that $x(\cdot)$ is increasing.
On the other hand, $f_p'$ is non-decreasing, which means that $(-f_p' \circ x)$ is increasing.
This fact, together with our induction assumption, allows us to write
\begin{eqnarray*}
x((n+1)\lambda) & \leq &
x_n + \int_{n\lambda}^{(n+1)\lambda} (-f_p' \circ x)((n+1)\lambda)) ~ dt
= x_n - \lambda f_p'(x((n+1)\lambda)),\\
\Leftrightarrow \quad x((n+1)\lambda) + \lambda f_p'(x((n+1)\lambda)) & \leq & x_n.
\end{eqnarray*}
Consider now the function $\phi : \dom f_p \to ]0,+\infty[$ defined by $\phi(t) = t+ \lambda f_p'(t)$.
It is clearly increasing and bijective on its image, so its inverse $\phi^{-1}$ is also increasing.
We observe moreover that, by definition, the proximal sequence satisfies $x_{n+1} = \phi^{-1}(x_n)$.
This allows us to write
\begin{equation*}
\phi(x((n+1)\lambda)) \leq 
 x_n 
 \quad
 \Leftrightarrow 
 \quad
x((n+1)\lambda) \leq \phi^{-1}(x_n) = x_{n+1}.
\end{equation*}
This ends the proof of the induction argument.

Observe 
that, given non-negative numbers $a,b>0$, the following inequality holds
\begin{equation*}
(\forall n \geq 1) \quad \sgn(p)(a+bn)^{\frac{1}{2-p}} \geq \sgn(p)(a+b)^{\frac{1}{2-p}} n^{\frac{1}{2-p}}.
\end{equation*}
This means that, for all $n \geq 1$,
\begin{equation*}
x_n \geq 
\sgn(p) \left( \vert x_0\vert ^{2-p} + p (p-2) \lambda n \right)^{\frac{1}{2-p}}
\geq 
\sgn(p) \left( \vert x_0\vert ^{2-p} + p (p-2) \lambda \right)^{\frac{1}{2-p}} n^{\frac{1}{2-p}} = 
\sgn(p) C_p n^{\frac{1}{2-p}}.
\end{equation*}
Passing this inequality through $f_p$ (which is non-decreasing) yields the desired result.
\end{proof}
}

\subsection{Proofs of Section~\ref{S:Geometry}}
\label{SS:equivalence geometric notions}

\subsubsection{Invariant sets and proofs of Section \ref{SS:geometry definitions}}

We provide here a result concerning the equivalence between all the notions in Definition \ref{D:geometric notions}, for a  large class of sets $\Omega \subset X$.
The sets $\Omega$ we will consider are directly related to the gradient flow induced by $\partial f$.
Given $u_0 \in \dom f$, it is known\footnote{See \cite[Thm 3.1]{Bre} when $u_0 \in \dom \partial f$, and \cite[Thm. 3.2]{Bre} with \cite[Cor. 16.39]{BauCom} when $u_0 \in \cl \dom f$.} that there exists a unique absolutely continuous trajectory noted $u(\cdot;u_0) : [0,+\infty[ \longrightarrow X$, called the steepest descent trajectory,  which satisfies:
\begin{equation}\label{D:SD}
(\text{for a.e. } t>0) \quad \frac{\d}{\dt} u(t;u_0) + \partial f(u(t;u_0)) \ni 0, \text{ and } u(0;u_0) = u_0.
\end{equation}
Following \cite{Bre}, we introduce the notion of invariant sets for the flow of $\partial f$:
\begin{definition}
\label{D:invariant set continuous dynamic}
A set $\Omega \subset X$ is \textit{$\partial f$-invariant} if for any  $x \in \Omega \cap \dom \partial f$ and a.e. $t>0$, $u(t;x) \in \Omega$ holds.
\end{definition}

\noindent In other words, $\Omega$ is said to be $\partial f$-invariant if any steepest descent trajectory starting in $\Omega$ remains therein.
It is straightforward to see that the intersection of two $\partial f$-invariant sets is still $\partial f$-invariant.

\begin{example}\label{Ex:invariant sets for continuous flow}
An easy way to construct a $\partial f$-invariant set is to consider the sublevel set of a \textit{Lyapunov} function $\psi : X \rightarrow \Rinf$ for the gradient flow induced by $\partial f$.
A function is said to be Lyapunov if for any $x \in \dom f$,  $\psi (u(\cdot;x)) : [0, +\infty [ \rightarrow \R$ is decreasing.
Classical examples of this kind are:
\begin{itemize}
	\item $\Omega=X$, which is $[\psi < 1]$ with $\psi = 0$.
	\item $\Omega = [f<r]$ for $r >\inf f$, which is $[\psi <r]$ with $\psi =f$ (see \cite[Thm. 3.2.17]{Bre}).
	\item $\Omega = \B(\bar x, \delta)$  for $\bar x \in \argmin f$, $\delta>0$, which is $[\psi < \delta]$ with $\psi(x)=\Vert x-\bar x \Vert$ (see \cite[Thm. 3.1.7]{Bre}).
	\item $\Omega = \{x \in X \ | \ \Vert \partial f(x) \Vert_\_ < M \}$ for $M>0$, which is $[\psi < M]$ with $\psi(x)=\Vert \partial f(x) \Vert_\_$ (see \cite[Thm. 3.1.6]{Bre}).
\end{itemize}
See \cite[Section IV.4]{Bre} for more details on the subject, as well as \cite{Bre70,LadLak74}.
It is also a good exercise to verify that the source sets considered in Proposition \ref{P:stability of source condition} are $\partial f$-invariant.
\end{example}

We next prove Proposition \ref{P:equivalence geometrical notions},  stating the equivalence between conditioning, metric subregularity and {\L}ojasiewicz on $\partial f$-invariant sets.
The proof   is based on an argument used in \cite[Theorem 5]{BolNguPeySut15}, which relies essentially on the following convergence rate property for the 
continuous steepest descent dynamic \eqref{D:SD}.

\begin{proof}[Proof of Proposition \ref{P:equivalence geometrical notions}]
Convexity of $f$ and the Cauchy-Schwartz inequality imply
\begin{equation*}
(\forall x \in \dom f) \quad f(x) - \inf f \leq \Vert \partial f(x) \Vert_\_ \dist(x, \argmin f ),
\end{equation*}
and so \ref{P:equivalence geometrical notions:i} $\implies$  \ref{P:equivalence geometrical notions:ii} $\implies$ 
\ref{P:equivalence geometrical notions:iii}.
Next, we just have to prove that the {\L}ojasiewicz property implies the conditioning one.
So let us assume that $f$ is $p$-{\L}ojasiewicz on $\Omega$, which is $\partial f$-invariant, and fix $x \in \Omega \cap \dom^* f$.
Define, for all $t \geq 0$, $\varphi(t):=(pc_{f,\Omega})^{-1}t^{1/p}$, which is derivable on $]0,+\infty[$, and for all $u\in \dom f$, $r(u)=f(u)-\inf f$.
Let us lighten the notations by noting $u(\cdot)$ instead of $u(\cdot;x)$, so that $u(0)=x$.
{Because we will need to distinguish the case in which the trajectory converges in finite time, we introduce $T:= \inf\{t \geq 0 \ | \ u(t) \in {\rm{argmin}}~f \} \in [0, + \infty]$.
Since $x \in \dom^*f$ and $u(\cdot)$ is continuous, we see that $T >0$.
For every $t \in [0,T[$, we have $u(t) \notin {\rm{argmin}}~f$, so $u(t) \in \Omega \cap \dom^* f$ and $r(u(t)) \neq 0$.
If $T < + \infty$, we also have for every $t > T$ that $u(t)=u(T)$ and $\dot u(t) =0$.}
Now, we write:
\begin{eqnarray*}
(\forall t \in ]0,T[) \quad \varphi(r(x)) \geq \varphi(r(x)) - \varphi (r(u(t))) & = & \displaystyle \int_t^0 (\varphi \circ r \circ u)'(\tau) \ \d \tau
= \int_t^0 \varphi' (( r \circ u )(\tau))\cdot (r \circ u)' (\tau) \ \d \tau.
\end{eqnarray*}
But $\frac{\d}{\d \tau} (r \circ u) (\tau) =- \Vert \dot u(\tau) \Vert^2 =  - \Vert \partial f(u(\tau)) \Vert_\_^2$ (see \cite{Bre}), so that the above equality becomes
\begin{equation}\label{rcs1}
(\forall t \in ]0,T[) \quad  \varphi(r(x)) \geq \int_0^t   \varphi' (( r \circ u )(\tau)) \Vert \partial f(u(\tau)) \Vert_\_^2 \ \d \tau.
\end{equation}
Since we assume $\Omega$ to be $\partial f$-invariant, we can apply the {\L}ojasiewicz inequality at $u(\tau) \in \Omega \cap \dom^* f$ for all $\tau \in ]0,t[$, which can be rewritten in this case as $1 \leq \varphi'(r(u(\tau))) \Vert \partial f(u(\tau)) \Vert_\_.$
This applied to \eqref{rcs1} gives us:
\begin{equation}\label{rcs3}
(\forall t \in ]0,T[) \quad  \varphi(r(x)) \geq \int_0^t \Vert \dot u(\tau ) \Vert \ \d \tau.
\end{equation}
From \eqref{rcs3} and the definition of $T$, we see that $\int_0^{+ \infty} \Vert \dot u(\tau ) \Vert \ \d \tau \leq \varphi(r(x)) < +\infty$, meaning that the trajectory $u(\cdot)$ has finite length.
As a consequence, it converges strongly to some $\bar u$ when $t$ tends to $+\infty$.
Finally, we use on \eqref{rcs3} the fact that $\Vert u(0) - u(t) \Vert \leq  \int_0^t \Vert \dot u(\tau ) \Vert \ \d \tau$, together with the fact that $\bar u \in \argmin f$ (see \cite[Thm. 3.11]{Bre}) to conclude that
\[
\frac{1}{pc_{f,\Omega}} \dist(x, \argmin f) \leq \frac{1}{pc_{f,\Omega}} \Vert x - \bar u \Vert \leq (f(x) - \inf f)^{1/p}. \tag*{\qedhere} 
\]
\end{proof}

\begin{proof}[Proof of Proposition~\ref{P:hierarchy of conditionings with p}]
\ref{P:hierarchy of conditionings with p:i}: let $S:=\argmin f \neq \emptyset$.
Given $\delta >0$, there exists $M \in ]0, + \infty[$ such that
\[
\sup \{  \dist(x,S) \ | \ {x \in \Omega \cap \delta \B_X} \} \leq M
\]
Since $f$ is $p$-conditioned on $\Omega$, we deduce that:
\begin{equation*}
(\forall x \in \Omega \cap \delta \B_X) \quad
\dist(x,S)^{p'} = \dist(x,S)^{p}\dist(x,S)^{p'-p} \leq (pM^{p'-p}/\gamma_{f,\Omega}) (f(x) - \inf f),
\end{equation*}
meaning that $f$ is $p'$-conditioned on $\Omega \cap \delta \B_X$. 

\ref{P:hierarchy of conditionings with p:ii}: the proof follows the same lines as in \ref{P:hierarchy of conditionings with p:i}.
\end{proof}

\begin{proof}[Proof of Proposition~\ref{P:trivial conditioning without minimizers}]
Assume by contradiction that there exists a sequence $(z^n)_\nin \subset \Omega$ such that 
\begin{equation}\label{jkl1}
n^{-1} \dist^p(z^n,\argmin f) > f(z^n) - \inf f.
\end{equation}
Since $\Omega$ is weakly compact, we can assume without loss of generality that $z^n$ weakly converges to some $z^\infty \in \Omega$ when $n \to + \infty$.
Then, it follows from \eqref{jkl1}, the boundedness of $(z^n)_\nin \subset \Omega$ and the weak lower semi-continuity of $f$ that $f(z^\infty) - \inf f \leq 0$, meaning that $z^\infty \in \argmin f$, contradicting $\Omega \cap \argmin f = \emptyset$.
\end{proof}

\subsubsection{Proofs of Section \ref{SS:geometry examples}}

{
\begin{lemma}[The \Loja constant for uniformly convex functions]\label{L:Loja for uniformly convex}
Let $f \in \Gamma_0(X)$ be uniformly convex, of order $p \geq 2$, with constant $\gamma$.
Then $f$ is $p$-\Loja on $X$, with $c_{f,X}=q^{-1/q}\gamma^{-1/p}$, where $1=(1/p)+(1/q)$.
\end{lemma}

\begin{proof}
Let $x \in \dom \partial f$, $\bar x \in \argmin f$, and $x^* \in \partial f(x)$.
By definition of uniformly convex functions
\begin{equation}\label{lcu1}
f(x) - \inf f = \sup\limits_{u \in X} \  f(x) - f(u)  \leq - \inf\limits_{u \in X}\left( \langle x^* , u - x \rangle + (\gamma/p) \Vert u - x \Vert^p \right).
\end{equation}
The right member of the above inequality involves a strictly convex optimization problem, whose unique optimal value $\bar u$ can be determined by using Fermat's rule: 
\begin{equation*}
0= x^* + \gamma \Vert \bar u - x \Vert^{p-2} (\bar u-x) \Leftrightarrow \bar u = x - \gamma^{-1/(p-1)} \Vert x^* \Vert^{(2-p)/(p-1)} x^*.
\end{equation*}
Injecting this optimal value in \eqref{lcu1} gives, after rearranging the terms,
\begin{equation*}
f(x) - \inf f \leq (1 - 1/p) \gamma^{-1/(p-1)} \Vert x^* \Vert^{p/(p-1)} ,
\end{equation*}
and, since $x^*$ is arbitrary in $\partial f (x)$, the result follows after passing this inequality to the power $1- 1/p$.
\end{proof}
}

\begin{proof}[Proof of Example \ref{Ex:general regularization}.ii)]
To prove the claim, it is enough to verify the three conditions of \cite[Theorem 4.2]{DruLew16}.
The first condition (boundedness of $\argmin f$) is guaranteed by the fact that $f$ is coercive.
Indeed, $h$ is strongly convex, therefore bounded from below, and $g$ is itself coercive.
The second condition (dual qualification conditions) follows immediately from the fact that both $h^*$ and $g^*$, and are continuously differentiable.
To see this, observe that in this example $g^*$ is (up to a constant) $\Vert \cdot \Vert_q^q$, where $q$ is the conjugate number of $p$: $(1/p) + (1/q) = 1$.
Moreover, $h$ being strongly convex means that $h^*$ is also continuously differentiable, with $\dom h^* = \mathbb{R}^M$.
The third condition (firm convexity) is easy to check for $h$ because it is strongly convex; for $g$ the proof is left in the following Lemma.
We can then apply \cite[Theorem 4.2]{DruLew16}, which ensures that $f$ is $2$-conditioned on every compact set.
Using again the fact that $f$ is coercive, and therefore has bounded sublevel sets, we conclude that $f$ is $2$-conditioned on every sublevel set.
\end{proof}

\subsubsection{Proofs of Section \ref{SS:sum rule}}

{
\begin{lemma}[$p$-powers are $2$-tilt conditioned when \text{$p \in ]1,2]$}]\label{L:p powers tilt 2 conditioned}
Let $p \in ]1,2]$, $u \in \mathbb{R}^N$, and $f : \mathbb{R}^N \rightarrow \mathbb{R}$ be defined as $f(x)=\frac{1}{p}\Vert x \Vert_p^p - \langle u,x \rangle$.
Then $f$ is $2$-conditioned on every bounded subset of $\mathbb{R}^N$.
\end{lemma}

\begin{proof}
This function is a separable sum, so, without loss of generality, we can assume from here that $N=1$  (see \cite[Lemma 4.4]{DruLew16}).
Given a real $t \in \mathbb{R}$, we will note its sign with $s(t)$, which is equal to $-1$ (resp. $+1$) if $t<0$ (resp. $t>0$), or $0$ if $t=0$.
Using the convexity, the differentiability of $f$, and the Fermat's rule, we see that $f$ admits a unique minimizer 
$\bar x$, defined by the relations
\begin{equation*}
    0 = s(\bar x) \vert \bar x \vert^{p-1} - u \Leftrightarrow
    \bar x = s(u) \vert u \vert^{\frac{1}{p-1}}
    \Leftrightarrow
    u = s(\bar x) \vert \bar x \vert^{p-1}.
\end{equation*}
If $u=0$, it is immediate to see that $f$ is $2$-conditioned on $]-1,1[$, where the relation $\vert t \vert^2 \leq \vert t \vert^p$ holds.
We therefore assume from now that $u \neq 0$, which also means that $\bar x \neq 0$.
We now compute  (we note $q = p/(p-1)$)
\begin{equation*}
    \inf f 
    = f(\bar x)
    = \frac{1}{p}\vert \bar x \vert^p - u \bar x
    = \frac{1}{p}\vert \bar x \vert^p - s(\bar x) \vert \bar x \vert^{p-1} \bar x
    = \frac{1}{p}\vert \bar x \vert^p - \vert \bar x \vert^{p}
    = - \frac{1}{q} \vert \bar x \vert^p,
\end{equation*}
meaning that we are looking for an inequality like
\begin{equation*}
    \gamma \vert x - \bar x \vert^2
    \leq
    \frac{1}{p}\vert x \vert^p - ux - \inf f 
    = \frac{1}{p}\vert x \vert^p - s(\bar x) \vert \bar x \vert^{p-1} x + \frac{1}{q} \vert \bar x \vert^p.    
\end{equation*}
Using the L'H\^opital rule twice allows us to study the following limit:
\begin{equation*}
    \lim\limits_{x \to \bar x} \ 
    \frac{\frac{1}{p}\vert x \vert^p - s(\bar x) \vert \bar x \vert^{p-1} x + \frac{1}{q} \vert \bar x \vert^p}{\vert x - \bar x \vert^2}
    =
     \lim\limits_{x \to \bar x} \ 
     \frac{s( x) \vert  x \vert^{p-1} - s(\bar x) \vert \bar x \vert^{p-1}}{2(x - \bar x)}
    =
    \lim\limits_{x \to \bar x} \ 
    \frac{(p-1)\vert x \vert^{p-2}}{2}
    =
    \frac{(p-1)}{2}\vert \bar x \vert^{p-2}.
\end{equation*}
Note that our assumption that $\bar x \neq 0$ ensures that we can take the derivative of the second numerator around $\bar x$.
Since this limit is well-defined, and nonnegative, it means that $f$ is $2$-conditioned on a small enough neighbourhood of $\bar x$.
To conclude the proof, it remains to verify that $f$ is $2$-conditioned on \textit{any}  bounded set.
This follows immediately from Proposition \ref{P:trivial conditioning without minimizers} and the fact that $\argmin f = \{\bar x \}$.
\end{proof}
}

{
\begin{lemma}[Kullback-Leibler divergences are $2$-tilt conditioned]\label{L:kullback libler tilt conditioning}
Let $\bar x \in ]0,+\infty[^N$, and $f \in \Gamma_0(\mathbb{R}^N)$ be the Kullback-Leibler divergence to $\bar x$:
\begin{equation*}
    f(x) = KL(\bar x ; x) = \sum_{i=1}^N kl(\bar x_i; x_i)
    \quad \text{ where } \quad 
    kl(\bar t; t) =
    \begin{cases}
    \bar t \log(\frac{\bar t}{t}) - \bar t + t & \text{ if }  t>0, \\
    +\infty & \text{ else.}
    \end{cases}
\end{equation*}
Then $f$ is $2$-tilt-conditioned on every bounded set of $\mathbb{R}^N$.
\end{lemma}

\begin{proof}
Let $d \in \mathbb{R}^N$, and define the tilted function $\tilde f = f + \langle d   \cdot \rangle$.
Using Fermat's rule, we see that ${\rm{argmin}}~f = \partial f^*(-d)$.
It is a simple exercice to verify that  $\dom \partial f^* = ]-\infty, 1[^N$, so ${\rm{argmin}}~\tilde f \neq \emptyset$ if and only if $d \in ]-1,+\infty[^N$.
Let $d$ be such vector, and write, for any $x_i>0$:
\begin{equation*}
    \tilde f_i(x_i) = 
    \bar x_i \log(\frac{\bar x_i}{x_i}) - \bar x_i + x_i + d_i x_i
    =
    (1+d_i) \left( \frac{\bar x_i}{1+d_i} \log(\frac{\bar x_i}{x_i}) - \frac{\bar x_i}{1+d_i} + x_i \right).
\end{equation*}
Let $X_i := \frac{\bar x_i}{1+d_i}$, which is well defined under our assumption that $d_i > -1$. Then 
\begin{equation*}
    \tilde f_i(x_i) = (1+d_i) \left( X_i \log(\frac{X_i}{x_i}) - X_i + x_i + X_i \log(1+d_i) \right) 
    = (1+d_i) kl(X_i;x_i) + a_i ,
\end{equation*}
where $a_i = X_i (1+d_i) \log(1+d_i) >0$.
We then observe that ${\rm{argmin}}~\tilde f_i = \{X_i\}$, from which we deduce that ${\rm{argmin}}~\tilde f = \{X\}$ with $X = (X_i)_{i=1}^N$.

Now, let $\delta >0$ be fixed, and let $x \in \mathbb{B}(X,\delta)$.
Let $\underline{d} := \min_i d_i > -1$, $c := N \Vert X \Vert_\infty$, and 
\begin{equation*}
    C := \frac{1}{\delta^2c} \left( \frac{\delta}{c} - \ln \left( 1+ \frac{\delta}{c} \right) \right)
    \quad \text{ which is nonnegative because $t > \ln(1+t)$ on $]0,+\infty[$}.
\end{equation*}
For each $i \in \{1,\dots,N\}$, we have $\vert x_i - X_i \vert \leq \delta$, so we can use \cite[Lem. A.2]{CalGarRosVil17} on $\tilde f_i$ to write
\begin{eqnarray*}
\tilde f(x) - \inf \tilde f 
& =& 
\sum_{i=1}^N \tilde f_i(x) - \tilde f_i(X_i)
 = 
\sum_{i=1}^N (1+d_i) kl(X_i;x_i) \\
& \geq &
\sum_{i=1}^N (1+d_i) C \vert X_i - x_i \vert^2
\geq 
(1+\underline{d})C \Vert X-x \Vert^2.
\end{eqnarray*}
This proves that $\tilde f$ is $2$-conditioned on $\mathbb{B}(X,\delta)$, which conludes the proof.
\end{proof}
}

\subsection{The Forward-Backward algorithm and proofs of Section \ref{S:CV rates for Forward Backward}}
\label{SS:Annex proof section 4}

\begin{definition}
Given a positive real sequence $(r_n)_{\nin}$ converging to zero, we say that $r_n$ converges:
\begin{itemize}
	\item \textit{sublinearly} (of order $\alpha \in ]0,+\infty[$) if  $\exists C \in ]0,+\infty[$ such that $\forall \nin$, $r_n \leq C n^{-\alpha}$,
	\item \textit{Q-linearly} if $\exists \eps \in ]0,1[$ such that $\forall \nin$, $r_{n+1} \leq \eps r_n$,
	\item \textit{R-linearly} if $\exists (s_n)_\nin$ Q-linearly converging such that $\forall \nin$, $r_n \leq s_n$,
	\item \textit{Q-superlinearly} (of order $\beta \in ]1,+\infty[$) if  $\exists C \in ]0,+\infty[$ such that  $\forall \nin$, $r_{n+1} \leq C r_n^\beta$,
	\item \textit{R-superlinearly}  if $\exists (s_n)_\nin$ Q-superlinearly convergent such that $\forall \nin$, $r_n \leq s_n$.
\end{itemize}
\end{definition}

\noindent It is easy to verify that $r_n$ is R-superlinearly convergent of order $\beta> 1$ if and only if
\begin{equation*}
(\forall \eps \in ]0,1[)(\exists C>0)(\forall \nin) \quad r_n \leq C \eps^{\beta^n}.
\end{equation*}
\noindent {Note that  $R$-linear and $R$-superlinear convergence ensures only  the overall  decrease of the sequence, while  
$Q$-linear and $Q$-superlinear convergence requires the sequence to decrease at a certain speed for each index. 
It is immediate from the definition that $Q$-convergence implies $R$-convergence.}

\begin{lemma}[Estimate for sublinear real sequences]\label{L:sublinear sequence}
Let $(r_n)_\nin$ be a real sequence being strictly positive and satisfying, for some $\kappa > 0$, $\alpha > 1$ and all $\nin$: $r_n - r_{n+1} \geq \kappa r_{n+1}^\alpha.$
Define $\tilde \kappa:= \min\{\kappa,\kappa^\frac{\alpha-1}{\alpha} \}$, and
$\delta:= \max\limits_{s \geq 1} \min \left\{ \frac{ \alpha-1}{s} , \kappa^{-\frac{\alpha - 1}{\alpha}} r_0^{1-\alpha} \left( 1 - s^{-\frac{\alpha - 1}{ \alpha}} \right) \right\} \in \left]0, + \infty\right [.$
Then,  for all $\nin$, $r_n \leq (\tilde \kappa \delta n)^{-1/(\alpha -1)}.$
\end{lemma}

\begin{proof}
It can be found in \cite[Lemma 7.1]{LiMor12}, see also the proofs of \cite[Theorem 2]{AttBol09} or \cite[Theorem 3.4]{FraGarPey15}.
\end{proof}

\begin{lemma}
\label{L:estimates for the forward-backward}
If Assumption \ref{ass:H}  holds, then for all $(x,u)\in X^2$ and all $\lambda >0$:
\begin{enumerate}[i)]
	\item\label{L:estimates for the forward-backward:descent} $\Vert \T x -  u \Vert^2 - \Vert x -  u \Vert^2 \leq \left( {\lambda L} - 1 \right)\Vert \T x - x \Vert^2 + 2\lambda ( f(u) - f(\T x)).$
	\item\label{L:estimates for the forward-backward:gradients} $\Vert \partial f(\T x) \Vert_\_ \leq \lambda^{-1} \Vert \T x - x \Vert \leq \Vert \partial f(x) \Vert_\_.$
\end{enumerate}
\end{lemma}

\begin{proof}[Proof of Lemma \ref{L:estimates for the forward-backward}]
To prove item i), start by writing
\begin{equation*}
\Vert \T x -  u \Vert^2 - \Vert x -  u \Vert^2 = - \Vert \T x - x \Vert^2 + 2\left\langle {x - \T x} , u - \T x \right\rangle.
\end{equation*}
The optimality condition in \eqref{D:prox} gives ${x - \T x}\in \lambda \partial g(\T x) + \lambda \nabla h(x)$ so that, by using the convexity of $g$:
\begin{equation*}
\Vert \T x -  u \Vert^2 - \Vert x -  u \Vert^2 \leq - \Vert \T x - x \Vert^2 + 2\lambda \left( g(u) - g(\T x) + \langle \nabla h(x),u - \T x \rangle \right).
\end{equation*}
Since we can write $\langle \nabla h(x),u -\T x \rangle = \langle \nabla h(x),u -x \rangle + \langle \nabla h(x),x - \T x\rangle$, we deduce from the convexity of $h$ and the Descent Lemma (\cite[Theorem 18.15]{BauCom}) that
\begin{equation*}
\langle \nabla h(x),u -\T x \rangle \leq h(u) -h(x) + h(x) - h(\T x) + \frac{L}{2} \Vert \T x - x \Vert^2 = h(u) - h(\T x) + \frac{L}{2} \Vert \T x - x \Vert^2.
\end{equation*}
Item i) is then proved after combining the two previous inequalities.
For item ii), the optimality condition in \eqref{D:prox}, together with a sum rule (see e.g. \cite[Theorem 3.30]{Pey}), to deduce that
\begin{equation}\label{gve1}
\forall (u,v) \in X^2, \quad v = \prox_{\lambda g}(u) \Leftrightarrow \lambda^{-1}(u-v) + \nabla h(v) \in \partial f(v).
\end{equation}
For the first inequality, use \eqref{gve1} with $(u,v)=(x-\lambda \nabla h(x),\T x)$, together with the contraction property of the gradient map $x \mapsto x - \lambda \nabla h(x)$ when $0<\lambda\leq 2/L$ (see \cite[Cor. 18.17 \& Prop. 4.39 \& Remark 4.34.i]{BauCom}) to obtain:
\begin{equation*}
 \Vert \partial f (\T x) \Vert_\_ \leq \lambda^{-1}\Vert (x - \lambda \nabla h(x) )-(\T x  - \lambda \nabla h(\T x)) \Vert \leq \lambda^{-1} \Vert \T x - x \Vert.
\end{equation*}
For the second inequality, consider $x^* := \proj(-\nabla h(x),\partial g(x))$, and use \eqref{gve1} with $(u,v)=(x + \lambda x^*,x)$, together with the nonexpansiveness of the proximal map (see \cite[Prop. 12.28]{BauCom}):
\begin{equation*}
\Vert \T x - x \Vert = \Vert   \prox_{\lambda g}(x - \lambda \nabla h(x)) - \prox_{\lambda g}(x + \lambda x^*)\Vert \leq \lambda \Vert   \nabla h(x)+ x^*\Vert = \lambda \Vert \partial f(x) \Vert_\_. \tag*{\qedhere}
\end{equation*}
\end{proof}

\begin{lemma}[Descent Lemma for H\"older smooth functions]\label{L:descent lemma}
Let $f : X \longrightarrow \R$ and $C \subset X$ be convex. Assume that $f$ is Gateaux differentiable on $C$, and that there exists $(\alpha,L) \in ]0,+\infty[^2$, such that for all $(x,y) \in C^2$, $\Vert \nabla f(x) - \nabla f(y) \Vert \leq L \Vert x - y \Vert^\alpha$ holds.
Then:
\begin{equation*}
(\forall (x,y) \in C^2) \quad f(y) - f(x) - \langle \nabla f(x) , y - x \rangle \leq \frac{L}{\alpha + 1	} \Vert x - y \Vert^{\alpha + 1 }.
\end{equation*}
\end{lemma}

\begin{proof}
The argument used in \cite[Remark 3.5.1]{Zal} for $C=X$ extends directly to convex sets.
\end{proof}

Now we can prove the convergence rate results of Section \ref{SS:Convergence with Lojasiewicz}:

\begin{proof}[Proof of Theorem \ref{T:CVKL discrete}]
We first show that $(x_n)_{\nin}$ has finite length.
Since $\inf f > - \infty$, $r_n:=f(x_n) - \inf f \in [0, + \infty[$, and it follows from Lemma~\ref{L:estimates for the forward-backward} that
\begin{eqnarray}
a \Vert x_{n+1} - x_n \Vert^2 & \leq & r_n - r_{n+1},  \text{ with } a=\frac{1}{2\lambda}(2-\lambda L) > 0,\label{cflH1} \\
 \Vert \partial f (x_{n+1}) \Vert_\_ & \leq & b \Vert x_n - x_{n+1} \Vert, \text{ with } b= \lambda^{-1}.\label{cflH2}
\end{eqnarray}
If there exists $\nin$ such that $r_n=0$ then the algorithm would stop after a finite number of iterations (see  \eqref{cflH1}), therefore it is not
restrictive to assume that $r_n>0$ for all $\nin$.
We set $\varphi(t):=p t^{1/p}$ and $c:=c_{f,\Omega}$, so that the {\L}ojasiewicz inequality at $x_n \in  \Omega \cap \dom^* f$ can be rewritten as
\begin{equation}\label{cfl2}
(\forall \nin) \quad 1 \leq c \varphi'(r_n) \Vert \partial f(x_n) \Vert_\_.
\end{equation}
Combining \eqref{cflH1}, \eqref{cflH2},  and  \eqref{cfl2}, and using the concavity of $\varphi$, we obtain for all $n \geq 1$:
\begin{eqnarray*}
\Vert x_{n+1} - x_n \Vert^2 & \leq & \frac{bc}{a}\varphi'(r_n)(r_n - r_{n+1}) \Vert x_n - x_{n-1} \Vert  \leq  \frac{bc}{a} (\varphi(r_n) - \varphi(r_{n+1})) \Vert x_n - x_{n-1} \Vert.
\end{eqnarray*}
By taking the square root on both sides, and using  Young's inequality, we obtain
\begin{equation}\label{cfl4}
(\forall n \geq 1) \quad 2 \Vert x_{n+1} - x_n \Vert \leq \frac{bc}{a}(\varphi(r_n) - \varphi(r_{n+1}))+ \Vert x_n - x_{n-1} \Vert.
\end{equation}
Sum  this inequality, and reorder the terms to finally obtain
\begin{equation*}\label{cfl4.5}
(\forall n \geq 1) \quad \sum\limits_{k=1}^{n} \Vert x_{k+1} - x_k \Vert \leq \frac{bc}{a}\varphi(r_1) + \Vert x_1 - x_{0} \Vert.
\end{equation*}
We deduce that $(x_n)_\nin$ has finite length and converges strongly to some $x_\infty$.
Moreover, from \eqref{cflH2} and the strong closedness of $\partial f : X \rightrightarrows X$, 
we conclude  that $0 \in \partial f(x_\infty)$.

Now we prove the convergence rates. 
Let $c=c_{f,\Omega}$ for short. 
We  first derive rates for the sequence of values $r_n:=f(x_n) - \inf f$, from which we will derive the rates for the iterates.
Equations \eqref{cflH1} and \eqref{cflH2} yield
\begin{equation*}\label{cfl5}
r_n-r_{n+1} \geq a \Vert x_{n+1} - x_n \Vert^2 \geq \frac{a}{b^2}\Vert \partial f(x_{n+1}) \Vert_\_^2.
\end{equation*}
The {\L}ojasiwecz inequality at $x_{n+1} \in \Omega \cap \doms$ implies $c^2 r_{n+1}^{2/p}(r_n - r_{n+1}) \geq  ab^{-2}  r_{n+1}^2,$
so we deduce that
\begin{equation}\label{cfl7}
(\forall \nin) \quad
r_{n+1} \neq 0 \ \Rightarrow \
 r_{n+1}^{2/p}(r_n - r_{n+1}) \geq \kappa r_{n+1}^2 , \quad \text{ with }  \kappa:={a}{(bc)^{-2}}.
\end{equation}
The rates for the values are derived from the analysis of the sequences satisfying the inequality in \eqref{cfl7}.
Depending on the value of $p$, we obtain different rates.

\noindent $\bullet$ If $p=1$, then we deduce from \eqref{cfl7} that for all $\nin,  r_{n+1}\neq 0$ implies $r_{n+1} \leq r_n - \kappa.$
Since the sequence $(r_n)_\nin$ is decreasing and positive,  $r_{n+1}\neq 0$ implies $ n\leq r_0\kappa^{-1}$.

For the other values of $p$, we will assume that $r_n >0$.
In particular, we get from \eqref{cfl7} 
\begin{equation}\label{cfl8}
(\forall \nin) \quad r_n - r_{n+1} \geq \kappa r_{n+1}^\alpha , \quad \text{ with }\  \alpha:= {2(p-1)}{p^{-1}}\  \text{ and }\  \kappa:={a}{b^{-2}c^{-2}}.
\end{equation}

\noindent $\bullet$ If $p\in]1,2[$, then $\alpha \in ]0,1[$. The positivity of $r_{n+1}$ and \eqref{cfl8} imply that for all $\nin$, $r_{n+1} \leq \kappa^{-1/\alpha} r_n^{1/\alpha}$,
meaning that $r_n$ converges Q-superlinearly.

\noindent $\bullet$ If $p=2$, then $\alpha=1$ and we deduce from \eqref{cfl8} that for all $\nin$, $r_{n+1} \leq {(1+\kappa)^{-1}} r_n$,
meaning that $r_n$ converges Q-linearly.

\noindent $\bullet$ If $p \in\left ]2,+\infty\right[$, then $\alpha \in \left ]1,2\right[$, and the analysis still relies on studying the asymptotic 
behaviour of a real sequence satisfying \eqref{cfl8}.
Lemma \ref{L:sublinear sequence} in the Annex shows that we have $r_{n+1} \leq (C_p')^{p/(p-2)} n^{-p/(p-2)}$, by taking
\begin{equation}
\label{e:cp}
(C_p')^{-1}:=\min\left\{\kappa,\kappa^{\frac{p-2}{2p-2}} \right\} 
\max\limits_{s \geq 1} \min \left\{ 
\frac{ p-2}{ps} ,
 \kappa^{\frac{2-p}{2p-2}} r_0^{\frac{2-p}{p}} \left( 1 - s^{-{\frac{p-2}{2p-2}}} \right) \right\}.
\end{equation}

To end the proof, we will prove that the rates for $\Vert x_n - x_\infty \Vert$ are governed by the ones of $r_n$.
Let $1\leq n \leq N < +\infty$, and sum the inequality in \eqref{cfl4} between $n$ and $N$ to obtain (remind that $b=\lambda^{-1}$):
\begin{equation*}
\|x_N-x_n\|\leq \sum\limits_{k=n}^{N} \Vert x_{k+1} - x_k \Vert \leq \frac{pc}{a \lambda} r_{n}^{1/p} + \Vert x_n - x_{n-1} \Vert.
\end{equation*}
Next, we pass to the limit for $N \to \infty$, we use \eqref{cflH1}, and the fact that $r_n$ is decreasing to obtain
\begin{equation}\label{cfl10}
(\forall n \geq 1) \quad \Vert x_\infty - x_n \Vert \leq \frac{pc}{a\lambda} r_{n-1}^{1/p} + \frac{1}{\sqrt{a}} \sqrt{r_{n-1}}.
\end{equation}
Note that ${r_{n-1}^{1/2}} \leq r_0^{\frac{1}{2}-\frac{1}{p}} r_{n-1}^{1/p}$ if $p\in\left[2,+\infty\right[$, and $r_{n-1}^{1/p} \leq r_0^{\frac{1}{p}-\frac{1}{2}} {r_{n-1}^{1/2}}$ if $p\in[1,2]$.
So, by defining
\begin{equation}
\label{e:cpp}
 C_p:=
\begin{cases}
2pc(2-\lambda L)^{-1} + (2\lambda r_0)^{1/2} (2-\lambda L)^{-1/2} r_0^{-1/p} & \text{ if } p \geq 2, \\
2pcr_0^{1/p}(2-\lambda L)^{-1}r_0^{-1/2} + (2\lambda )^{1/2} (2-\lambda L)^{-1/2} & \text{ if } p \leq 2,
\end{cases}
\end{equation}
we finally conclude from \eqref{cfl10} that $\Vert x_\infty - x_n \Vert \leq C_p r_{n-1}^{1/\max\{2,p\}}$ when $n \geq 1$. 
\end{proof}

\begin{proof}[Proof of Theorem \ref{T:CVKL discrete rates negative p}]
The proof is as for the case $p \in\left ]2,+\infty\right[$ of Theorem \ref{T:CVKL discrete}:
the $p$-{\L}ojasiewicz property implies \eqref{cfl7}, and the statement follows from  Lemma \ref{L:sublinear sequence} with $\alpha=2(p-1)/p \in \left]2,+\infty\right[$.
\end{proof}

\begin{proof}[Proof of Theorem \ref{T:CV general descent}]
The proofs of Theorems~\ref{T:CVKL discrete} and  \ref{T:CVKL discrete rates negative p}  rely on  the combination of the {\L}ojasiewicz inequality with 
the estimations \eqref{cflH1} and \eqref{cflH2}, which can be replaced by  \eqref{H1} and \eqref{H2}.
\end{proof}

\subsection{Linear inverse problems and proofs of Section \ref{SS:least squares in Hilbert spaces}}
\label{SS:Annex proof section 5.1}

\noindent Here we will make use of is the Moore-Penrose \textit{pseudo-inverse} of $A$.
It is a linear operator (not necessarily bounded), whose domain is $D(A^\dagger):= R(A) + R(A)^\perp$, and satisfying
\begin{equation*}
    (\forall y \in D(A^\dagger)) \quad
    A^\dagger y := \argmin \{ \Vert x \Vert \ | \ A^*Ax = A^*y \}.
\end{equation*}
It is easy to see that, whenever $y \in D(A^\dagger)$, the solution set of \eqref{D:least squares inverse problems} is  $A^\dagger y + \ker A$.

{
\begin{lemma}\label{L:function of operator commutes HELL}
Let $A$ be a bounded linear opertator from $X$ to $Y$.
Then, for every continuous function $\phi : [0,+\infty[ \rightarrow \mathbb{R}$, we have $A \phi(A^*A) = \phi(AA^*) A$.
\end{lemma}

\begin{proof}
A simple induction argument shows that, for every $k \geq 0$, $A(A^*A)^k = (AA^*)^k A$.
Taking linear combinations of this equality allows to see that, for every polynomial $P \in \mathbb{R}[X]$,  $AP(A^*A) = P(AA^*)A$.
Now, if $\phi$ is continuous on $[0,+\infty[$, it is in particular continuous on $[0, \Vert A \Vert^2]$, which is an interval containing the spectrum of both $A^*A$ and $AA^*$. 
Thus, $\phi$ restricted to this interval can be written as the uniform limit of a sequence of polynomials.
Passing to the limit (see \cite[Thm. VI.32.1]{Hel69}) in the last equality  gives the desired result.
\end{proof}
}

\begin{lemma}\label{L:range of A and sqrt AA*}
For all $b \in Y$, $r \in \left]0,+\infty\right[$, the following two properties are equivalent:
\begin{enumerate}
\item $(\exists x \in \ker A^\perp) \quad b=Ax, \quad \Vert x \Vert = r $ 
\item $(\exists y \in \cl R(A)) \quad b=\sqrt{AA^*} y, \quad \Vert y \Vert=r.$
\end{enumerate}
\end{lemma}

\begin{proof}
It is shown in \cite[Proposition 2.18]{EngHanNeu} that $R(A) = R( \sqrt{AA^*})$, so it is enough to verify this implication:
\begin{equation*}
(\forall (x,y)\in \ker A^\perp \times \cl R(A))\quad \ Ax = \sqrt{AA^*} y \ \Rightarrow \ \Vert x \Vert = \Vert y \Vert.
\end{equation*}
Let $(x,y)$ be such a pair.
Since $Ax=\sqrt{AA^*}y$ and $y \in \cl R(A)=\ker \sqrt{AA^*}^\perp$, we deduce that $y=(\sqrt{AA^*})^\dagger Ax$.
Therefore, since $AA^*$ is self-adjoint, $(AA^*)^\dagger Ax = (A^*)^\dagger x$ (see \cite[p.35]{EngHanNeu}), and $A^*(A^*)^\dagger x=\proj(x;\ker A^\perp)$, we get
\begin{equation*}
\Vert y \Vert^2 = \Vert (\sqrt{AA^*})^\dagger Ax \Vert^2 = 
\langle ({AA^*})^\dagger Ax, Ax \rangle  = 
 \langle A^*(A^*)^\dagger x,x \rangle = 
 \Vert x \Vert^2. \tag*{\qedhere}
\end{equation*}
\end{proof}

\begin{proof}[Proof of Lemma \ref{L:regularity trough data space}]
Remind that $y^\dagger=Ax^\dagger=AA^\dagger y$ and let $\nu=\mu+1/2$.
Then, Lemma \ref{L:range of A and sqrt AA*} yields:
$$\begin{array}{rclc}
b \in A^{-1}Y_{\nu,\delta} & \Leftrightarrow & (\exists \omega \in \cl R(A))\quad \Vert \omega \Vert \leq \delta, \quad Ab = y^\dagger + (AA^*)^\mu (AA^*)^\frac{1}{2} \omega & \text{ with $\nu=\mu+1/2$,}\\ 
&\Leftrightarrow & (\exists w \in \ker A^\perp)\quad \Vert w \Vert \leq \delta, \quad Ab = AA^\dagger y + (AA^*)^\mu A w & \text{with Lemma  \ref{L:range of A and sqrt AA*}, } \\
& \Leftrightarrow & (\exists w \in \ker A^\perp)\quad \Vert w \Vert \leq \delta, \quad Ab= AA^\dagger y + A(A^*A)^\mu w &  \text{ with Lemma \ref{L:function of operator commutes HELL}, }\\
& \Leftrightarrow & (\exists w \in \ker A^\perp)\quad \Vert w \Vert \leq \delta, \quad b - x^\dagger -(A^*A)^\mu w \in \ker A &  \\
& \Leftrightarrow & b \in X_{\mu,\delta}. & \hspace{0.28\linewidth}
\end{array}$$
\end{proof}

\begin{lemma}[Interpolation inequality {\cite[p. 55]{EngHanNeu}}]\label{L:interpolation inequality}
For all $x \in X$ and $0\leq \alpha < \beta$, we have
\begin{equation*}
\Vert (\AA)^\alpha x \Vert \leq \Vert (\AA)^\beta x \Vert^{\frac{\alpha}{\beta}} \ \Vert x \Vert^{1- \frac{\alpha}{\beta}}.
\end{equation*}
\end{lemma}

{
\begin{lemma}[Powers of self-adjoint operators]\label{L:powers self adjoint operators}
Let $S$ be a bounded selfadjoint positive linear operator on a Hilbert space.
Then, for all $\alpha >0$, $\ker S = \ker S^\alpha$, and $\cl R(S^\alpha) = \cl R(S)$.
\end{lemma}

\begin{proof}
Given any $0<\alpha<\beta$, we can write $S^\beta = S^{\beta - \alpha} S^\alpha$, from which we deduce that $\ker S^\alpha \subset \ker S^\beta$.
This means that $(\ker S^\alpha)_{\alpha >0}$ is a nondecreasing family.
To prove that this family is constant, it is enough to see that $\ker S^2 \subset \ker S$, which we verify now:
If $x \in \ker S^2$, then $\Vert Sx \Vert^2 = \langle Sx,Sx \rangle = \langle S^2 x,x \rangle = 0$, therefore $x \in \mbox{Ker}~S$.
The conclusion follows from the fact that $\ker S^\perp = \cl R(S)$.
\end{proof}
}

\subsection{Regularized inverse problems and proofs of Section \ref{SS:sparse inverse problems}}
\label{SS:Annex proof section 5.2}

{

\begin{proposition}\label{P: matrix coercive on injective cone}
Let $K \subset \mathbb{R}^N$ be a closed cone and $S\in\mathcal{S}_+(\mathbb{R}^N)$.
Then $S$ is coercive on $K$ if and only if $K \cap \ker S = \{0\}$.
\end{proposition}

\begin{proof}
The direct implication is immediate from Definition \ref{D:coercive operator}.
For the reverse implication, let $K$ be a closed cone such that $K \cap \ker S = \{0\}$.
Since $S$ is linear, we know that $d \mapsto \langle Sd,d \rangle$ is convex and continuous.
So, using the compactness of $K \cap \mathbb{S}_{\mathbb{R}^N}$ we deduce that:
\begin{equation}
\label{e:dbar}
(\exists \bar d \in K \cap \S)\quad \inf\limits_{d \in K \cap \S} \langle S d,d \rangle = \langle S \bar d, \bar d \rangle.
\end{equation}
Because $\bar d \in K$ and $\bar d \neq 0$, we deduce from our assumption  that $\bar d \notin \mbox{Ker}~S$.
Therefore, $ \gamma :=\langle S \bar d, \bar d \rangle > 0$, from which we deduce that $S$ is $\gamma$-coercive on $K$.
\end{proof}

\begin{definition}[Cone enlargement]\label{def:enlcone}
Let $K \subset \mathbb{R}^N$ be a cone, and $ \theta \in [0, \frac{\pi}{2}]$.
We define the $\theta$-enlargement of $K$ as
\begin{equation*}
    K_\theta := \mathbb{R} \left\{ x \in \mathbb{S}_{\mathbb{R}^N} \ | \ (\exists y \in K \cap \mathbb{S}_{\mathbb{R}^N}) \ \arccos \left( {\vert \langle x,y \rangle \vert}\right) \leq \theta \right\}.
\end{equation*}
\end{definition}

\begin{lemma}\label{L:enlarged cone is closed}
If $K$ is a closed cone, then $K_\theta$ is a closed cone containing $K$ for all $\theta \in [0, \frac{\pi}{2}]$.
\end{lemma}

\begin{proof}
By definition, $K_\theta$ is a cone containing $K$ and $\Delta_\theta:=\left\{ x \in \mathbb{S}_{\mathbb{R}^N} \ | \ (\exists y \in K \cap \mathbb{S}_{\mathbb{R}^N}) \ \arccos \left( {\vert \langle x,y \rangle \vert}\right) \leq \theta \right\}$ is compact, due to the 
compactness of $K\cap \S$.
Since $0\not\in \Delta_\theta$, by compactness of $\Delta_\theta$,  we deduce that $K_\theta=\R \Delta_\theta$ is a closed cone (see e.g. \cite[Proposition A.1.1]{Gar}).
\end{proof}

\begin{proposition}\label{P:coercivity on enlarged cones}
Let $S\in\mathcal{S}_+(\mathbb{R}^N)$ which is $\gamma$-coercive on a closed cone $K$.
Then, for every $\gamma' \in ]0,\gamma]$, $S$ is $\gamma'$-coercive on $K_\theta$, with 
$\theta:=\arcsin \left( \frac{\gamma - \gamma'}{\Vert S \Vert} \right) \in [0, \frac{\pi}{2}[$.
\end{proposition}

\begin{proof}
Let $\theta$ and $\gamma$ be as in the statement.
Since $S$ is $\gamma$-coercive on $K$, we see that $\gamma \leq \| S \|$, which guarantees that $\theta \in [0, \frac{\pi}{2}[$.
Now, the fact that $K_\theta$ is closed (Lemma \ref{L:enlarged cone is closed}) implies that $K_\theta \cap \S$ is compact in $X$, so we can use the same arguments as in \eqref{e:dbar} to deduce that there exists $\bar d \in K_\theta \cap \mathbb{S}_{\mathbb{R}^N}$ such that $\langle S\bar d, \bar d \rangle = \inf\limits_{d \in K_\theta\cap \mathbb{S}_{\mathbb{R}^N}} \langle Sd,d \rangle$.
Since $\bar d \in K_\theta$, there exists by definition of $K_\theta$ some $\bar v \in K \cap \S$ such that $\arccos(\vert \langle \bar d, \bar v \rangle \vert) \leq \theta$.
We can use \cite[Theorem 1]{KnyArg06} to write
\begin{equation}\label{elf0}
    \vert \langle S \bar v, \bar v \rangle - \langle S \bar d, \bar d \rangle \vert
    \leq  \|  S \| \sin\arccos(\vert \langle \bar v, \bar d \rangle \vert )\leq\| S \| \sin \theta.
\end{equation}
Since $\bar v \in K \cap \mathbb{S}_{\mathbb{R}^N} \subset K_\theta \cap \mathbb{S}_{\mathbb{R}^N}$, we have $\langle S \bar v, \bar v \rangle \geq \langle S \bar d, \bar d \rangle$. 
Moreover, $\arccos(\vert \langle \bar v, \bar d \rangle \vert) \leq \theta$, so  \eqref{elf0}, implies
\begin{equation*}
 \langle S \bar d , \bar d \rangle \geq \langle S \bar v , \bar v \rangle - \| S \| \sin \theta \geq \gamma - \| S \| \sin \theta = \gamma'.
\end{equation*}
We deduce from the definition of $\bar d$ that $S$ is $\gamma'$-coercive on $K_\theta$.
\end{proof}

\begin{proposition}\label{P:curvature manifold}
Let $C \subset \mathbb{R}^N$, and $\bar x \in C$.
\begin{enumerate}[label=\roman*)]
    \item\label{P:curvature manifold:prox regular} For $\rho >0$, $C$ is $\rho$-prox-regular at $\bar x$ if and only if :
\begin{equation}\label{cmd1}
    (\forall x \in C)(\forall \eta \in N_C(\bar x)) \quad 
    \langle \eta , x - \bar x \rangle \leq \frac{\rho}{2} \Vert \eta \Vert \Vert x - \bar x \Vert^2.
\end{equation}
    \item\label{P:curvature manifold:local} If $C$ is a $C^2$ manifold, then there exists $\delta,\rho >0$ such that $C \cap \mathbb{B}(\bar x, \delta)$ is $\rho$-prox-regular.
\end{enumerate}
\end{proposition}

\begin{proof}
 Item \ref{P:curvature manifold:prox regular} :
Definition \ref{D:prox regular set} can be rewritten as $(\forall \eta \in N_C(\bar x) \cap \mathbb{S}_{\mathbb{R}^N}) (\forall x \in C) \quad x \notin \mathbb{B}(\bar x + \frac{1}{\rho} \eta, \frac{1}{\rho})$, where the condition $x \notin \mathbb{B}(\bar x + \frac{1}{\rho} \eta, \frac{1}{\rho})$ is equivalent to, after developing the square:
\begin{eqnarray*}
\frac{1}{\rho^2} \leq \Vert x - \bar x - \frac{1}{\rho} \eta \Vert^2 
=
\Vert x - \bar x \Vert^2 + \frac{1}{\rho^2} \Vert \eta \Vert^2 - \frac{2}{\rho} \langle x - \bar x , \eta \rangle
=
\Vert x - \bar x \Vert^2 + \frac{1}{\rho^2} - \frac{2}{\rho} \langle x - \bar x , \eta \rangle.
\end{eqnarray*}
The conclusion follows after cancelling and reorganizing the terms. Item \ref{P:curvature manifold:local} : Every $C^2$-manifold is prox-regular in the sense of \cite[Def. 10.23 \& Prop. 13.32]{RocWet}.
Therefore, for every $\bar x \in C$, there exists $\delta, \rho >0$ such that for every $x \in C \cap \mathbb{B}(\bar x,\delta)$, and for every $\eta \in N_C(\bar x) \cap \mathbb{S}_{\mathbb{R}^N}$, the inequality \eqref{cmd1} holds \cite[Exercice 13.31]{RocWet}. Conclusion follows from the fact that $N_C(\bar x) = N_{C\cap \mathbb{B}(\bar x,\delta)}(\bar x)$.

\end{proof}

Here is a needed result estimating locally the coercivity of an operator on a prox-regular set via its coercivity on the tangent cone.
}
{
\begin{proposition}\label{P:ellipticity lifted from tangent to manifold}
Let $C \subset X$ be $\rho$-prox-regular at $\bar x \in C$.
Let $S : X \rightarrow X$ be a bounded positive selfadjoint linear operator, being $\gamma$-coercive on $T_C(\bar x)$.
Then, for all $\gamma' \in \left]0, \gamma\right[$, there exists a cone $K \subset X$
such that $S$ is $\gamma'$-coercive on $K$, and $C \cap \B_X(\bar x,\delta) \subset \bar x + K$, with $\delta = \frac{2(\gamma - \gamma')}{\rho \Vert S \Vert}$.
\end{proposition}

\begin{proof}
Let $\gamma' \in \left ]0, \gamma\right [$ be fixed, and define $\theta:=\arcsin((\gamma - \gamma')\| S \|^{-1}) \in ]0, \frac{\pi}{2}[$.
Let $K_\theta$ be the $\theta$-enlargement of $T_{C}(\bar x)$, then Proposition \ref{P:coercivity on enlarged cones} guarantees that $S$ is $\gamma'$-coercive on $K_\theta$.
It remains to prove that there exists 
$\delta \in \left]0,+\infty\right[$ such that $C  \cap \mathbb{B}(\bar x , \delta) \subset \bar x +  K_\theta$.
Let $x \in C$. Because $C$ is $\rho$-reached at $\bar x$, we know that $T_C(\bar x)$ is a convex cone (use \cite[Thm. 4.8.(12)]{Fed59} and the fact that $C$ is locally closed at $\bar x$), so we can define $y:= \proj(x - \bar x, T_C(\bar x)) $, and $\eta :=  \proj(x - \bar x, N_C(\bar x))$.
Using Moreau's Theorem \cite[Thm. 6.30]{BauCom}, we deduce that $\eta =x - \bar x -y$ with $\langle \eta, y  \rangle = 0$.
We define $\delta := \Vert x- \bar x\Vert$, and look for a condition on it so that $x \in \bar x + K_\theta$.
For this to happen, it is enough to verify that
\begin{equation}\label{ftg1}
    \langle x - \bar x, y \rangle  \geq \cos(\theta) \Vert x - \bar x \Vert \Vert y \Vert.
\end{equation}
Now, use Proposition \ref{P:curvature manifold}.\ref{P:curvature manifold:prox regular} together with the Cauchy-Schwarz inequality, and the polynomial inequality $X^2 - cX \geq c^2/4$, to write
\begin{equation*}
    \Vert y \Vert^2 = \Vert x - \bar x - \eta \Vert^2 \geq 
    \Vert x-\bar x \Vert^2 + \Vert \eta \Vert^2 - \rho \Vert \eta \Vert \Vert x - \bar x \Vert^2
    \geq 
    \delta^2(1 - \rho^2 \delta^2 /4).
\end{equation*}
We can use this inequality, together with the facts that $x - \bar x = y + \eta$ and $\langle y,\eta  \rangle =0$, to write
\begin{equation*}
    \langle x - \bar x,y \rangle^2
    =
    \Vert y \Vert^4
    \geq
    \Vert y \Vert^2 \delta^2(1 - \rho^2 \delta^2 /4).
\end{equation*}
This allows us to conclude that \eqref{ftg1} holds as long as:
\begin{equation*}
    1 - \rho^2 \delta^2 /4 \geq \cos(\theta)^2
    \ \Leftrightarrow  \ 
    \rho^2 \delta^2 /4 \leq 1 -\cos(\theta)^2 
    \ \Leftrightarrow  \  
    \rho \delta /2 \leq \sin(\theta) = \frac{\gamma - \gamma'}{\Vert S \Vert}.
\end{equation*}
\end{proof}
}

\begin{proof}[Proof of Proposition \ref{P:ellipticity implies conditioning}]
\noindent   Let $0<\gamma' < \gamma$, and set $S:=\argmin f$. 
Since $h$ is of class $C^2$ around $\bar x \in S$, there exists some $\delta >0$ such that for all $u \in \delta \B_X$, $\| \nabla^2 h(\bar x + u) - \nabla^2 h(\bar x) \| \leq \gamma - \gamma'.$
Notice that when $\nabla^2 h$ is Lipschitz continuous, we can take $\delta=(\gamma - \gamma')/L$. Also, if it is constant, we can just take $\delta= +\infty$ and $\gamma'= \gamma$. 
Let us show that $f$ is $2$-conditioned on $\Omega:= \bar x + (K \cap \delta \B_X)$ with the constant $\gamma_{f,\Omega}=\gamma'$.
Take $x \in \Omega \cap \dom g$ and use the optimality condition at $\bar x \in S$ and the convexity of $g$ to obtain
\begin{eqnarray*}
f(x) - \inf f 
& = & g(x) - g(\bar x) + \langle \nabla h(\bar x), x - \bar x \rangle + h(x) - h(\bar x) - \langle \nabla h (\bar x), x - \bar x) \rangle 
\geq
h(x) - h(\bar x) - \langle \nabla h (\bar x), x - \bar x \rangle.
\end{eqnarray*}
By {Taylor}'s theorem applied to $h$, we deduce from the inequality above that there exists $y \in [x,\bar x]$ such that:
\begin{eqnarray*}
f(x)- \inf f & \geq & 
(1/2) \langle \nabla^2 h(\bar x) (x-\bar x), x - \bar x \rangle+
(1/2)\langle \left[ \nabla^2 h(y) - \nabla^2 h(\bar x)\right] (x- \bar x), x - \bar x \rangle.
\end{eqnarray*}
On the one hand, since $x\in\Omega$, we have that $x - \bar x \in K$. Thus, from the coercivity of $\nabla^2 h(\bar x)$ we have
\begin{equation*}
\langle \nabla^2 h(\bar x) (x-\bar x), x - \bar x \rangle \geq \gamma \Vert x -\bar x \Vert^2.
\end{equation*}
On the other hand, we  use the Cauchy-Schwarz inequality together with the definition of $\delta$ and the fact that $\Vert y - \bar x \Vert \leq \Vert x - \bar x \Vert < \delta$ to obtain
\begin{equation*}\label{qch3}
\langle \left[ \nabla^2 h(y) - \nabla^2 h(\bar x)\right] (x- \bar x), x - \bar x \rangle \geq - (\gamma - \gamma') \Vert x - \bar x \Vert^2.
\end{equation*}
By combining the three previous inequalities, we deduce that
\begin{equation}\label{qch2}
 f(x) - \inf f \geq (\gamma'/2) \Vert x - \bar x \Vert^2.
\end{equation}
This implies that $(\bar{x}+K)\cap \argmin f=\{\bar{x}\}$,
and the statement follows from $\Vert x - \bar x \Vert \geq  \dist (x;S)$. 
\end{proof}

\begin{footnotesize}

\end{footnotesize}


\begin{thebibliography}{00}
\bibitem{AbsMahAnd05} P.-A. Absil, R. Mahony and B. Andrews, \textit{Convergence of the iterates of descent methods for analytic cost functions}, SIAM Journal on Optimization, \textbf{16}, pp. 531--547, 2005.

\bibitem{AraGeo08} F.J. Arag\'on Artacho and M.H. Geoffroy, \textit{Characterization of metric regularity of subdifferentials}, Journal of Convex Analysis, \textbf{15}(2), pp.365--380, 2008.

\bibitem{AttBol09} H. Attouch and J. Bolte, \textit{On the convergence of the proximal algorithm for nonsmooth functions involving analytic features}, Mathematical Programming, \textbf{116}(1-2), pp. 5--16, 2009.

\bibitem{AttBolRedSou10} H. Attouch, J. Bolte, P. Redont and A. Soubeyran,  \textit{Proximal alternating minimization and projection
methods for nonconvex problems. An approach based on the Kurdyka-\L ojasiewicz inequality}, Mathematics of 
Operations Research, \textbf{35}(2), pp. 438--457, 2010.

\bibitem{AttBolSva13} H. Attouch, J. Bolte and B.F. Svaiter, \textit{Convergence of descent methods for semi-algebraic and tame problems: proximal algorithms, forward-backward splitting, and regularized Gauss-Seidel methods},  Mathematical Programming, \textbf{137}(1-2), pp. 91--129, 2013.

\bibitem{AttWet93} H. Attouch and R. Wets, \textit{Quantitative stability of variational systems II, a framework for nonlinear conditioning}, SIAM Journal on Optimization, \textbf{3}(2), pp. 359--381, 1993.

\bibitem{AzeCor14} D. Az\'e and J.-N. Corvellec, \textit{Nonlinear local error bounds via a change of metric}, Journal of Fixed Point Theory and Applications, \textbf{16}(1), pp. 351--372, 2014.

\bibitem{Bai78} J.-B. Baillon, \textit{Un exemple concernant le comportement asymptotique de la solution du problème $du/dt + \partial \vartheta \ni 0$}, Journal of Functional Analysis, \textbf{28}(3), pp. 369--376, 1978.

\bibitem{BanDobMixSaw13} A.S. Bandeira, E. Dobriban, D.G. Mixon, and W.F. Sawin, \textit{Certifying the restricted isometry property is hard}, IEEE Transactions on Information Theory,  \textbf{59}(6), pp.  3448--3450, 2013. 

\bibitem{BauBor93} H.H. Bauschke and J.M. Borwein, \textit{On the convergence of von Neumann's alternating projection algorithm for two sets}, Set-Valued Analysis, \textbf{1}(2), pp. 185--212, 1993.

\bibitem{BauCom} H.H. Bauschke and P. Combettes, \textit{Convex analysis and monotone operator theory}, 2nd Edition, Springer, 2017.

\bibitem{BecTeb09} A. Beck and M. Teboulle,  {\it  A fast iterative shrinkage-thresholding algorithm for linear inverse problems},  SIAM Journal on Imaging Sciences, \textbf{2}(1), pp. 183--202, 2009.

\bibitem{BegBolJen15} P. B\'egout, J. Bolte and M.A. Jendoubi, \textit{On damped second-order gradient systems}, Journal of Differential Equations, \textbf{259}(7), pp. 3115--3143, 2015.

\bibitem{BolDanLew07} J. Bolte, A. Daniilidis and A. Lewis, \textit{The {\L}ojasiewicz Inequality for Nonsmooth Subanalytic Functions with Applications to Subgradient Dynamical Systems}, SIAM Journal on Optimization, \textbf{17}(4), pp. 1205--1223, 2007.

\bibitem{BolDanLewShi07} J. Bolte, A. Daniilidis, A.S. Lewis and M. Shiota, \textit{Clarke subgradients of stratifiable functions}, SIAM Journal on Optimization, \textbf{18}(2), pp. 556--572, 2007.

\bibitem{BolDanLeyMaz10} J. Bolte, A. Daniilidis, O. Ley and L. Mazet, \textit{Characterizations of \L ojasiewicz inequalities: Subgradient flows, talweg, convexity}, Transactions of the American Mathematical Society, \textbf{362}, pp. 3319--3363, 2010.

\bibitem{BolNguPeySut15} J. Bolte, T.P. Nguyen, J. Peypouquet and B. Suter, \textit{From error bounds to the complexity of first-order descent methods for convex functions}, Mathematical Programming \textbf{165}(2), pp. 471--507, 2017.

\bibitem{BolSabTeb13} J. Bolte, S. Sabach and M.  Teboulle, \textit{Proximal alternating linearized minimization for nonconvex and nonsmooth problems}, Mathematical Programming, \textbf{146}(1-2), pp. 459--494,  2013.

\bibitem{BonSha00} J.F. Bonnans and A. Shapiro, \textit{Perturbation Analysis of Optimization Problems}, Springer-Verlag, New York, 2000.

\bibitem{BreLor08b} K. Bredies, and D.A. Lorenz, \textit{Linear convergence of iterative soft-thresholding. Journal of Fourier Analysis and Applications}, \textbf{14}(5-6), pp. 813--837, 2008.

\bibitem{Bre} H. Br\'ezis, \textit{Op\'erateurs maximaux monotones et semi-groupes de contractions dans les espaces de Hilbert}, North-Holland/Elsevier, New-York, 1973.

\bibitem{Bre70} H. Br\'ezis,\textit{On a characterization of flow-invariant sets}, Communications on Pure and Applied Mathematics, \textbf{23}(2), pp. 261--263, 1970.

\bibitem{BurFer93} J. Burke and M.C. Ferris, \textit{Weak Sharp Minima in Mathematical Programming}, SIAM Journal on Control and Optimization, \textbf{31}(5), pp. 1340--1359, 1993.

\bibitem{CalGarRosVil17} L. Calatroni, G. Garrigos, L. Rosasco, and S. Villa, \textit{Accelerated iterative regularization via dual diagonal descent}, preprint on arXiv:1912.12153, 2019.

\bibitem{Can08} E.J. Cand\`es, \textit{The restricted isometry property and its implications for compressed sensing}, Comptes Rendus Mathematique, \textbf{346}(9-10), pp. 589-592, 2008.

\bibitem{ChaRecParWil12} V. Chandrasekaran, B. Recht, P.A. Parillo and A.S. Willsky, \textit{The Convex Geometry of Linear Inverse Problems}, Foundations of Computational Mathematics, \textbf{12}(6), pp. 805--849, 2012.

\bibitem{ChoPesRep14} E. Chouzenoux, J.-C. Pesquet and A. Repetti, \textit{A block coordinate variable metric forward-backward algorithm}, Journal on Global Optimization, \textbf{66}, pp. 457--485, 2016.

\bibitem{ComPes11} P.L. Combettes and J.-C. Pesquet,  \textit{Proximal splitting methods in signal processing}, in Fixed-point algorithms for inverse problems in science and engineering, Springer New York, 2011.


\bibitem{CorJouZal97} O. Cornejo, A. Jourani and C. Zalinescu,  \textit{Conditioning and Upper-Lipschitz Inverse Subdifferentials in Nonsmooth Optimization Problems}, Journal of Optimization Theory and Applications, \textbf{95}(1), pp. 127--148, 1997.

\bibitem{CraGoc20} D. K. Crane and M. Gockenbach, \textit{The Singular Value Expansion for Arbitrary Bounded Linear Operators}, Mathematics, \textbf{8}(8), 2020.

\bibitem{DauDefDem04} I. Daubechies, M. Defrise and C. De Mol, \textit{An iterative thresholding algorithm for linear inverse problems with a sparsity constraint}, Communications in Pure and Applied Mathematics, \textbf{57}(11), pp. 1413--1457, 2004.

\bibitem{DavYin14} D. Davis and W. Yin, \textit{Convergence rate analysis of several splitting schemes},in: Splitting Methods in Communication, Imaging, Science, and Engineering, Springer International Publishing, 2014.

\bibitem{DevCapRos05} E. De Vito, A. Caponnetto and L. Rosasco, \textit{Model selection for regularized least-squares algorithm in learning theory}, Foundations of Computational Mathematics, \textbf{5}(1), pp. 59--85, 2005.

\bibitem{DevRosVer06} Y. Yao, L. Rosasco and A. Caponnetto \textit{On regularizationearly stopping in gradient descent learning}, Constructive Approximation {\bf 26}, pp. 289--315, 2007.

\bibitem{Dev86} R. DeVore, Approximation of functions, Approximation Theory, Proceedings Symp. Applied Mathematics, AMS 36 (1986) 1-20.

\bibitem{DonLewRoc03} A.L. Dontchev, A.S. Lewis and R.T. Rockafellar, \textit{The Radius of Metric Regularity}, Transactions of the American Mathematical Society, \textbf{355}(2), pp. 493--517, 2003.

\bibitem{DonRoc09} A. Dontchev and T. Rockafellar, \textit{Implicit functions and Solution  Mappings}, Springer, New York, 2009.

\bibitem{DonZol93} A. Dontchev and T. Zolezzi, \textit{Well-posed Optimization Problems}, Springer-Verlag, Berlin,1993.

\bibitem{DruIof15} D. Drusvyatskiy and A.D. Ioffe, \textit{Quadratic growth and critical point stability of semi-algebraic functions}, Mathematical Programming, \textbf{153}(2) Ser. A, pp. 635--653, 2015.

\bibitem{DruLew16} D. Drusvyatskiy and A.D. Lewis, \textit{Error bounds, quadratic growth, and linear convergence of proximal methods}, Mathematics of Operations Research \textbf{43}, pp. 693--1050, 2018.

\bibitem{DruMorNgh14} D. Drusvyatskiy, B.S. Mordukhovich, and T.T.A. Nghia, \textit{Second-order growth, tilt stability, and metric regularity of the subdifferential}, Journal of Convex Analysis, \textbf{21}(4), pp. 1165--1192, 2014.


\bibitem{EngHanNeu} H. Engl, M. Hanke, and A. Neubauer, \textit{Regularization of Inverse Problems}, Kluwer, Dordrecht, 1996.

\bibitem{FadMalPey18} J. Fadili, J. Malick and G. Peyr\'e, \textit{Sensitivity Analysis for Mirror-Stratifiable Convex Functions}, SIAM Journal on Optimization, \textbf{28}(4), pp. 2975--3000, 2018.

\bibitem{Fed59} H. Federer, \textit{Curvature Measures}, Transactions of the American Mathematical Society, \textbf{93}(3), pp. 418--491, 1959.

\bibitem{Fer91} M.C. Ferris, \textit{Finite termination of the proximal point algorithm}, Mathematical Programming, \textbf{50}, pp. 359--366, 1991.

\bibitem{FraGarPey15} P. Frankel, G. Garrigos and J. Peypouquet, \textit{Splitting methods with variable metric for Kurdyka-\L ojasiewicz functions and general convergence rates}, Journal of Optimization Theory and Applications, \textbf{165}(3), pp. 874--900, 2015.

\bibitem{FouRau} S. Foucart and H. Rauhut, \textit{A mathematical introduction to compressive sensing}, Springer, 2013.

\bibitem{Gar} G. Garrigos, \textit{Descent dynamical systems and algorithms for tame optimization and multi-objective problems}, Ph.D. thesis, 2015. Available on \url{https://tel.archives-ouvertes.fr/tel-01245406}

\bibitem{GarRosVil20} G. Garrigos, L. Rosasco and S. Villa, \textit{Thresholding gradient methods in Hilbert spaces: support identification and linear convergence}, ESAIM Control Optimization and Calculus of Variations, \textbf{26}, 28 (20 pages) 2020. 

\bibitem{Gol62} A.A. Goldstein, \textit{Cauchy's method of minimization}, Numerische Mathematik, \textbf{4}(1), pp. 146--150, 1962.

\bibitem{Gro77} C. W. Groetsch, \textit{Generalized inverses of linear operators: representation and approximation}, Dekker, 1977.


\bibitem{Gul91} O. G\"uler, \textit{On the Convergence of the Proximal Point Algorithm for Convex Minimization}, SIAM Journal on Control and Optimization, \textbf{29}(2), pp. 403--419, 1991.

\bibitem{HarJen11} A. Haraux and M.A. Jendoubi, \textit{The \L ojasiewicz gradient inequality in the infinite dimensional Hilbert space framework}, Journal of Functional Analysis, \textbf{260}(9), pp. 2826--2842, 2011.

\bibitem{HarLew04} W.L. Hare and A.S. Lewis, \textit{Identifying Active Constraints via Partial Smoothness and Prox-Regularity}, Journal of Convex Analysis, \textbf{11}(2), pp. 251--266, 2004.

\bibitem{HarLew07} W.L. Hare and A.S. Lewis, \textit{Identifying active manifolds}, Algorithmic Operations Research, \textbf{2}(2), pp.  75--82, 2007.

\bibitem{Hel69} G. Helmberg, \textit{Introduction to Spectral Theory in Hilbert Space}, Elsevier, 1969.

\bibitem{HirLem} J.-B. Hiriart-Urruty and C. Lemar\'echal, \textit{Convex analysis and minimization algorithms I:  Fundamentals}, Springer Science \& Business Media, 1993.

\bibitem{Hof52} A.J. Hoffman, \textit{On approximate solutions of systems of linear inequalities}, Journal of Research of the National Bureau of Standards, \textbf{49}(4), pp. 263--265, 1952.

\bibitem{Hoh02} T. Hohage, \textit{Lecture notes on inverse problems}, Vorlesungskript, University of G\"ottingen, Germany, 2002.

\bibitem{HouZhoSoLuo13} K. Hou, Z. Zhou, A. M.-C. So and Z.-Q. Luo, \textit{On the Linear Convergence of the Proximal Gradient Method for Trace Norm Regularization}, in: Advances in Neural Information Processing Systems, pp. 710--718, 2013.

\bibitem{KarNutSch16} H. Karimi, J. Nutini and M. Schmidt, \textit{Linear Convergence of Gradient and Proximal-Gradient Methods Under the Polyak-\L ojasiewicz Condition}, in: Machine Learning and Knowledge Discovery in Databases (ECML PKDD). Lecture Notes in Computer Science, vol 9851. Springer, 2016.  

\bibitem{KnyArg06} A.V. Knyazev and M.E. Argentati, \textit{On proximity of Rayleigh quotients for different vectors and Ritz values generated by different trial subspaces}, Linear algebra and its applications, \textbf{415}(1), pp. 82--95, 2006.

\bibitem{LadLak74} G.S. Ladde and V. Lakshmikantham, \textit{On flow-invariant sets}, Pacific Journal of Mathematics, \textbf{51}(1), pp. 215--220, 1974.

\bibitem{Lem92} B. Lemaire, \textit{About the Convergence of the Proximal Method}, in Advances in Optimization, Lecture Notes in Economics and Mathematical Systems, \textbf{382}, pp. 39--51, 1992.

\bibitem{Lem96} B. Lemaire, \textit{Stability of the iteration method for non expansive mappings}, Serdica Mathematical Journal, \textbf{22}(3), pp. 331--340, 1996.

\bibitem{Lem98} B. Lemaire, \textit{Well-posedness, conditioning and regularization of minimization, inclusion and fixed-point problems}, Pliska Studia Mathematica Bulgarica, \textbf{12}(1), pp. 71--84, 1998.

\bibitem{Lev09} D. Leventhal, \textit{Metric subregularity and the proximal point method}, Journal of Mathematical Analysis and Applications, \textbf{360}(2), pp. 681--688, 2009.

\bibitem{Lew02} A.S. Lewis, \textit{Active sets, nonsmoothness, and sensitivity}, SIAM Journal on Optimization, \textbf{13}(3), pp. 702--725, 2002.

\bibitem{LewMal08} A. Lewis and J. Malick, \textit{Alternating projections on manifolds}, Mathematics of Operations Research, \textbf{33}(1), pp. 216--234, 2008.

\bibitem{Li95} W. Li, \textit{Error Bounds for Piecewise Convex Quadratic Programs and Applications}, SIAM Journal on Control and Optimization, \textbf{33}(5), pp. 1510--1529, 1995.

\bibitem{Li13} G. Li,  \textit{Global error bounds for piecewise convex polynomials},  Mathematical Programming, \textbf{137}(1-2), Ser. A, pp.
37--64, 2013.

\bibitem{LiMor12} G. Li and B. Mordukhovich, \textit{H\"older Metric Subregularity with Applications to Proximal Point Method}, SIAM Journal on Optimization, \textbf{22}(4), pp. 1655--1684, 2012.

\bibitem{LiMorPha15}
G. Li, B. S. Mordukhovich and T. S. Pham, \textit{New fractional error bounds for polynomial systems with
applications to Holderian stability in optimization and spectral theory of tensors}, Mathematical Programming, \textbf{153}, pp. 333--362, 2015.

\bibitem{LiPon17} G. Li and T.K. Pong, \textit{Calculus of the exponent of Kurdika-{\L}ojasiewicz inequality and its applications to linear convergence of first-order methods}, Foundations of Computational Mathematics, \textbf{18}, pp.1199-1232, 2018.

\bibitem{LiaFadPey14} J. Liang, J. Fadili and G. Peyr\'e, \textit{Local linear convergence of Forward--Backward under partial smoothness}, in: Advances in Neural Information Processing Systems, pp. 1970--1978, 2014.

\bibitem{LiaFadPey17} J. Liang, J. Fadili and G. Peyr\'e, \textit{Activity identification and local linear convergence of Forward-Backward-type methods}, SIAM Journal on Optimization {\bf 27}, pp. 408--437, 2017.

\bibitem{LiaFadPey16} J. Liang, J. Fadili and G. Peyr\'e, \textit{A Multi-step Inertial Forward--Backward Splitting Method for Non-convex Optimization},  in: Advances in Neural Information Processing Systems, pp. 4042--4050, 2016.

\bibitem{LiuWriReBitSri15} J. Liu, S.J. Wright, C. R\'e, V. Bittorf and S. Sridhar, \textit{An Asynchronous Parallel Stochastic Coordinate Descent Algorithm}, Journal of Machine Learning Research, \textbf{16}(1), pp. 285--322, 2015.

\bibitem{Loj63} S. {\L}ojasiewicz, \textit{Une propri\'et\'e topologique des sous-ensembles analytiques r\'eels}, in: Les \'Equations aux
D\'eriv\'ees Partielles, \'Editions du centre National de la Recherche Scientifique, Paris, pp. 87--89, 1963.

\bibitem{Luk13} R. Luke, \textit{Prox-Regularity of Rank Constraint Sets and Implications for Algorithms}, Journal of Mathematical Imaging and Vision, \textbf{47}(3), pp. 231--238, 2013.

\bibitem{LuoTse93} Z. Q. Luo and P. Tseng, \textit{Error bounds and convergence analysis of feasible descent methods: a general approach}. Annals of Operations Research, \textbf{46}(1), pp. 157--178, 1993.

\bibitem{Luq84} F. Luque, \textit{Asymptotic Convergence Analysis of the Proximal Point Algorithm}, SIAM Journal on Control and Optimization, \textbf{22}(2), pp. 277--293, 1984.

\bibitem{MerPie10} B. Merlet, M. Pierre, \textit{Convergence to equilibrium for the backward Euler scheme and applications}, Commun. Pure Appl. Anal, \textbf{9}, pp. 685--702, 2010.

\bibitem{MesRosSan10}  S.  Mosci, L. Rosasco,M.  Santoro, A. Verri,and S.  Villa, {\em  Solving structured sparsity regularization with proximal methods},  in Joint European Conference on Machine Learning and Knowledge Discovery in Databases (pp. 418-433). Springer Berlin Heidelberg.

\bibitem{NecNesGli15} I. Necoara, Y. Nesterov and F. Glineur, \textit{Linear convergence of first order methods for non-strongly convex optimization}, Mathematical Programming, \textbf{175}, pp. 69--107,2019.

\bibitem{Pen96} J.-P. Penot, \textit{Conditioning convex and nonconvex problems}, Journal of Optimization Theory and Applications, \textbf{93}(3), pp. 535--554, 1996.

\bibitem{Pey} J. Peypouquet, \textit{Convex optimization in normed spaces. Theory, methods and examples.}, Springer Science \& Business media, 2015.

\bibitem{PolRoc96} R.A. Poliquin and R.T. Rockafellar, \textit{Prox-regular functions in variational analysis}, Transactions of the American Mathematical Society, \textbf{348}(5), pp. 1805--1838, 1996.

\bibitem{Pol63} B.T. Polyak, \textit{Gradient methods for minimizing functionals}, Zh. Vychisl. Mat. Mat. Fiz.,  \textbf{3}(4), pp. 643--653, 1963.

\bibitem{Pol} B.T. Polyak, \textit{Introduction to Optimization}, Optimization Software, New York, 1987.

\bibitem{Roc76} R.T. Rockafellar, \textit{Monotone Operators and the Proximal Point Algorithm}, SIAM Journal on Control and Optimization, \textbf{14}(5), pp. 877--898, 1976.

\bibitem{Roc} R.T. Rockafellar, \textit{Convex Analysis}, Princeton University Press, 1996.

\bibitem{RocWet} R.T. Rockafellar and R. J.-B. Wets \textit{Variational Analysis}, Springer Science \& Business Media, 2009.

\bibitem{Sal16} S. Salzo, \textit{The variable metric forward-backward splitting algorithm under mild differentiability assumptions}, SIAM Journal on Optimization, \textbf{27}(4), pp. 2153--2181, 2017.

\bibitem{SchLerBac11} M. Schmidt, N. Le Roux and F. Bach, \textit{Convergence Rates of Inexact Proximal-Gradient Methods for Convex Optimization}, in Advances in neural information processing systems, pp. 1458--1466, 2011.

\bibitem{Spi85} J.E. Spingarn, \textit{Applications of the method of partial inverses to convex programming: Decomposition}, Mathematical Programming, \textbf{32}(2), pp. 199--223, 1985.

\bibitem{Spi87} J.E. Spingarn, \textit{A projection method for least-squares solutions to overdetermined systems of linear inequalities}, Linear Algebra and its Applications, \textbf{86}, pp. 211--236, 1987.

\bibitem{Vai70} M. M. Vainberg, \textit{Le probl\`eme de la minimisation des fonctionelles non lin\'eaires}, C.I.M.E. IV
ciclo (1970).

\bibitem{VaiPeyFad14} S. Vaiter, G. Peyr\'e, and J.M. Fadili, \textit{ Model consistency of partly smooth regularizers}, IEEE Transactions on Information Theory, \textbf{64}(3), pp. 1725--1737, 2017.

\bibitem{Wri93} S. Wright, \textit{Identifiable Surfaces in Constrained Optimization}, SIAM Journal on Control and Optimization, \textbf{31}(4), pp. 1063--1079, 1993.

\bibitem{Zal} C. Zalinescu, \textit{Convex Analysis in General Vector Spaces}, Singapore: World Scientific, 2002.

\bibitem{ZhaTre95} R. Zhang and J. Treiman, \textit{Upper-Lipschitz Multifunction and Inverse Subdifferentials}, Nonlinear Analysis: Theory, Methods, and Applications, \textbf{24},
pp. 273--286, 1995.

\bibitem{ZhoSo15} Z. Zhou, and A.M.-C. So, \textit{A unified approach to error bounds for structured convex optimization problems}, Mathematical Programming {\bf 165}, pp. 689--728, 2017.

\bibitem{ZhoZhaSo15} Z. Zhou, Q. Zhang, and A.M.-C. So, \textit{$\ell_{1,p}$-Norm Regularization: Error Bounds and Convergence Rate Analysis of First-Order Methods.}, in Proceedings of the 32nd International Conference on Machine Learning, pp. 1501--1510, 2015.

\bibitem{Zol78} T. Zolezzi, \textit{On equiwellset minimum problems}, Appl. Math. Optim, \textbf{4}, pp. 209--223, 1978.
\end{thebibliography}
\end{document}